\documentclass[a4paper,oneside,10pt]{article}

\usepackage{a4}
\usepackage{amsmath}
\usepackage{amssymb}
\usepackage{amsthm}
\usepackage{amscd}
\usepackage{epsfig}

\parskip\smallskipamount

\swapnumbers
\newtheorem{thm}{Theorem}
\newtheorem{prop}[thm]{Proposition}
\newtheorem{lem}[thm]{Lemma}
\newtheorem{cor}[thm]{Corollary}
\newtheorem*{intro}{Theorem}
\newtheorem*{intro2}{Corollary}

\theoremstyle{definition}
\newtheorem{dfn}[thm]{Definition}
\newtheorem{rem}[thm]{Remark}
\numberwithin{thm}{section}

\def\IN{\mathbb N}
\def\IZ{\mathbb Z}
\def\IR{\mathbb R}
\def\IC{\mathbb C}
\def\dx{\mathrm d}
\def\lap{\Delta}
\def\H{\mathrm H} 
\def\C{\mathrm C} 
\def\L{\mathrm L}
\def\vol{\,\mathrm {vol}}
\def\Ord{\mathrm O}
\def\ord{\mathrm o}
\def\total{\mathcal E}
\def\base{\mathcal B}
\def\inbase{\mathcal U}
\def\Zeta{\mathcal Z}

\providecommand{\abs}[1]{\left|#1\right|}
\providecommand{\norm}[1]{\left\lVert#1\right\rVert}
\providecommand{\scal}[2]{\langle#1\, ,\,#2\rangle}
\providecommand{\st}{\,|\,}

\providecommand{\dist}[2]{\mathrm d(#1,#2)}
\DeclareMathOperator{\arsinh}{Arsinh}
\DeclareMathOperator{\arcosh}{Arcosh}
\DeclareMathOperator{\supp}{supp}
\DeclareMathOperator{\tr}{tr}

\DeclareMathOperator{\re}{Re}
\DeclareMathOperator{\im}{Im}
\DeclareMathOperator{\interior}{int}
\DeclareMathOperator{\Li}{Li_2}

\def\isom{\mathrm{isom}}
\def\sig{\sigma} 

\providecommand{\inv}[2]{(#1,#2)}
\providecommand{\graph}[1]{{#1}}

\begin{document}
\pagestyle{empty}
\begin{titlepage}
\end{titlepage}
\pagestyle{plain}
\title{\vskip-5mm On the resolvent of the Laplacian on functions\\ for degenerating surfaces of finite geometry}
\author{Michael Schulze}
\date{August 29, 2004}
\maketitle
\begin{abstract}
  We consider families $(Y_n)$ of degenerating hyperbolic surfaces.
  The surfaces are geometrically finite of fixed topological type. Let
  $\Zeta_n$ be the Selberg Zeta function of $Y_n$, and let $\Zeta^d_n$
  be the contribution of the pinched geodesics to $\Zeta_n$. Extending
  a result of Wolpert's, we prove that $\Zeta_n(s)/\Zeta^d_n(s)$
  converges to the Zeta function of the limit surface for all $s$ with
  $\re(s)>1/2$. The technique is an examination of resolvent of the
  Laplacian, which is composed from that for elementary surfaces via
  meromorphic Fredholm theory. The resolvent $(\lap_n-t)^{-1}$ is
  shown to converge for all $t\notin[1/4,\infty)$. We also use this
 property to define approximate Eisenstein functions and scattering matrices.
\end{abstract}
\tableofcontents
\setcounter{section}{-1}
\section{Introduction}
A family of \emph{degenerating Riemannian manifolds} consists of a
manifold $M$ and a family $(g_\ell)_{\ell>0}$ of Riemannian metrics on
$M$ that meet the following assumptions:
\begin{itemize}
\item There are finitely many disjoint open subsets $Z_i\subset M$
  that are diffeomorphic to cylinders $F_i\times J_i$. The fibre
  $F_i$ is a compact manifold and $J_i\subset\IR$ is a neighbourhood
  of $0$.
\item The restriction of each metric $g_\ell$ to $Z_i=F_i\times J_i$
  is a product metric
  \begin{equation*}
    (x,a)\longmapsto \nu_{i,\ell}(a)\cdot g_{F_i}(x)+ \mu_{i,\ell}(a)\cdot \dx a^2
  \end{equation*}
  such that  $\nu_{i,\ell}(0)\to 0$ and $\mu_{i,\ell}(0)\to \infty$ as
  $\ell\to 0.$
\item On the complement of $\bigcup_i F_i\times\{0\}$ in $M$, the
  metrics $g_\ell$ converge to a Riemannian metric $g_0$.
\end{itemize}
Spectral geometric properties of certain types of degenerating
manifolds have been examined by several authors, let us mention
Colbois-Courtois \cite{colbois-courtois2}, Chavel-Dodziuk
\cite{chavel-dodziuk} and Judge \cite{judge1}. 

We consider the case of hyperbolic surfaces, which is the fundamental
example of such a degeneration. Here $M$ is an oriented surface of
negative Euler characteristic, and the metrics $g_\ell$ are
hyperbolic, chosen in such a way that there are finitely
many closed curves $c_i$, geodesic with respect to all metrics, with
the length $\ell_i$ of each curve converging to 0 as $\ell$ decreases.
On the complement of the distinguished curves, the sequence of metrics
is required to converge to a hyperbolic metric.

In the description above, the geodesics $c_i$ correspond to the
central fibres $F_i\times\{0\}$. The collar lemma of hyperbolic
geometry ensures that each $c_i$ has got a collar neighbourhood
$Z_i=\IR/\IZ\times (-\epsilon,\epsilon)$, with the Riemannian metric
on $Z_i$ being given by
\begin{equation*} 
    (x,a)\longmapsto (\ell_i^2+a^2)\dx x^2+(\ell_i^2+a^2)^{-1}\dx
    a^2.
\end{equation*}
Let $M_{\ell}$ denote the surface $M$ equipped with the metric
$g_\ell$ if $\ell>0$, and let $M_0=M\setminus\bigcup_i c_i$ carry the
limit metric $\lim_{\ell\to 0}g_{\ell}$. Note that $M_0$ is a complete
hyperbolic surface by definition, which contains a pair of cusps for
each $i$.

From the point of view of spectral theory, this example was initially
studied by Schoen-Wolpert-Yau \cite{schoen-wolpert-yau},
Colbois-Courtois \cite{colbois-courtois1} and by Hejhal \cite{hej3},
Ji \cite{ji1} and Wol\-pert \cite{wolpert}. To exemplify how the
spectrum may behave during this process, assume that $M$ is compact
for the moment.  Then $M_0$ is of finite area but not compact. The
spectrum of compact manifolds is purely discrete, whereas that of
$M_0$ is the union of finitely many eigenvalues in $[0,1/4)$, and the
essential spectrum is $[1/4,\infty)$. It was observed that the small
eigenvalues of $M_0$ are limits of eigenvalues of $M_\ell$ as $\ell\to
0$, and the eigenvalues of $M_\ell$ accumulate at each point in
$[1/4,\infty)$.  One also tries to obtain information on
embedded eigenvalues in the essential spectrum of $M_0$ by means of
such an approximation \cite{wolpert2}.

Hejhal and Wolpert proved results on the behaviour of the Selberg Zeta
function. Our motivation for this work was to extend one of these
results.

Let us recall the definition of the Selberg Zeta function. It is a
meromorphic function $\Zeta$ on the complex plane, associated with a
hyperbolic surface. In the domain $\{s\in\IC\st \re(s)>1\}$ it is
given by an absolutely convergent product
\begin{equation*}
  \Zeta\colon\ s\longmapsto \prod_c \Zeta_{c}(s),\quad\text{where}\quad \Zeta_{c}(s):=\prod_{k=0}^\infty \bigl(1-e^{-(s+k)\ell(c)}\bigr)^2.
\end{equation*}
The product ranges over the set of all unoriented, simple, closed
geodesics $c$ of the surface, and $\ell(c)$ denotes the length of $c$.
Now if $\ell(c)$ decreases as the metric changes, we have (the precise
asymptotics are given in lemma \ref{asymptotic})
\begin{equation*}
  \Zeta_c(s)=\Ord\bigl(\ell(c)^{1-2\re(s)}e^{-\pi^2/3\ell(c)}\bigr),\quad\ell(c)\to 0.
\end{equation*}
Therefore, one cannot expect to Selberg Zeta function of $M_\ell$ to
converge to that of the limit $M_0$.  In section \ref{zeta function}
we prove
\begin{intro} Consider a family $M_\ell$ of degenerating surfaces, and let $\{c_i\}$ be the set of distinguished geodesics that are pinched. Then $\Zeta(s)/\prod_i \Zeta_{c_i}(s)$ converges to the Zeta function of the limit surface $M_0$ if\/ $\re(s)>1/2$.
\end{intro}
This is Wolpert's \emph{Conjecture 2} \cite{wolpert}. Hejhal proved
that it holds in the domain of convergence $\re(s)>1$, and this was
extended by Wolpert to a neighbourhood of $s=1$.  He also concluded
from the functional equation of $\Zeta$ that the same cannot be true
in any domain that intersects $\{s\in\IC\st \re(s)<1/2\}$, at least if
the surfaces $M_\ell$ are compact. 

Each zero $s$ of the Zeta function with $\re(s)>1/2$ corrensponds to
an eigenvalue $s(1-s)\in [0,1/4)$ of the Laplacian.  We also prove
Wolpert's \emph{Conjecture 1}, which states that one may divide out
these zeroes from the quotient above to obtain uniform boundedness
from below.  It applies to geometrically finite surfaces of
both finite and infinite area.

Let us elaborate a little on the proof of these statements. Analysis
of the Selberg Zeta function is based on the trace formula, which
relates the logarithmic derivative of $\Zeta$ with the resolvent
kernel of the Laplacian. Our primary object of investigation is the
resolvent operator. The following theorem, proved throughout section
\ref{resolvent}, extends a result of Jorgenson and Lundelius
\cite{jorgenson-lundelius-heat} that arose as a consequence of their
discussion of the heat kernel.
\begin{intro} If\/ $t\in\IC\setminus[1/4,\infty)$, let $R_\ell(t)$ denote the pull-back of the resolvent $(\lap-t)^{-1}$ from $M_\ell$ to $M_0$. Then $R_\ell(t)$ converges to $R_0(t)$ in the topology of continuous linear maps $\L^2(M_0)\to \L^2_\mathrm{loc}(M_0)$.
\end{intro}
Restriction to the vector space $\L^2_\mathrm{loc}(M_0)$ in the image
of the resolvent means that we prove convergence of the resolvent
kernel on the complement of small neighbourhoods of the pinched
geodesics. This restriction may be dropped if $(\lap-t)^{-1}$ is
replaced with
\begin{equation}\label{difference}
  (\lap-t)^{-1}-\sum_i \psi_i (\lap_{\bar Z_i}-t)^{-1}\phi_i,
\end{equation}
where $(\lap_{\bar Z_i}-t)^{-1}$ denotes the resolvent on an infinite
cylinder $\bar Z_i$ that admits an isometric embedding $Z_i\to\bar
Z_i$. The functions $\psi_i$, $\phi_i$ are suitable cut-off functions
such that the operator is defined. We prove that if two operators like
\eqref{difference} with different values of $t$ are considered, their
difference converges in the trace class topology. This immediately
implies our results on the Zeta function.

The notion of convergence in the previous theorem is to be understood
in the following sense: The map $t\mapsto R_\ell(t)$ is meromorphic on
$\IC\setminus[1/4,\infty)$ with possibly finitely many poles of finite
rank in $[0,1/4)$. If $t_0$ is not a pole of $R_{0}$, then there exist
neighbourhoods $V$ of $t_0$ and $\inbase$ of $0$ such that $R_\ell$
has no poles in $V$ for all $\ell\in\inbase$, and convergence holds
uniformly in $V$. So if $\cal C$ is a closed curve in the complement
of the discrete spectrum of $R_{0}$, then the Riesz-projector
\begin{equation*}
  -\frac 1{2\pi i}\oint_{\cal C} R_\ell(t)\,\dx t
\end{equation*}
converges to that of $R_0$ as $\ell\to 0$. In particular, we obtain an alternative proof of the convergence of small eigenvalues for degenerating surfaces \cite{colbois-courtois1}.

To prove the previous theorem, we apply a technique that is typically
used to establish a continuation of the resolvent kernel for $M_0$,
across the essential spectrum $[1/4,\infty)$, to a branched cover of
$\IC$. More precisely, we use meromorphic Fredholm theory to compose
the resolvent of a geometrically finite surface from those for a
compact surface and for elementary quotients of the hyperbolic plane.
If $\re(s)>1/2$, and $s(1-s)$ is not an eigenvalue von $\lap$, then
the following equation holds, where $\lap_i$ denotes the Laplacian of
an auxiliary surface:
\begin{equation}\label{int1}
  (\lap-s(1-s))^{-1}=\Bigl(\sum_i\psi_i(\lap_i-s(1-s))^{-1} \phi_i\Bigr)(1+K(s))^{-1}.
\end{equation}
Here $(1+K(s))$ denotes a meromorphically invertible family of
Fredholm operators.  This formula gives the information needed to
deduce convergence of the resolvent from that of its constituent
parts.

Now the right-hand side of equation \eqref{int1} is known
to have a continuation in $s$ to the complex plane. In view of this
analytic continuation, it would be interesting to obtain similar
results on the left of the critical axis $\{s\in\IC\st \re(s)=1/2\}$.
But obviously the continuation is symmetric in $s$ and $1-s$ for compact
surfaces, while there is no such trivial relation in the other cases.
To overcome this symmetry, we introduce \emph{approximate
  Eisenstein functions}\/ for $M_\ell$ as follows.

Recall that the definition of Eisenstein functions is based on certain
eigenfunctions of the Laplacian on a cusp $Z_k^-\subset M_0$,
where
\begin{equation*}
  Z_k^-:=\IR/\IZ\times (-\epsilon,0)\subset Z_k
\end{equation*}
carries the metric $(x,a)\mapsto a^2\dx x^2+a^{-2}\dx a^2$. In these
coordinates, the functions are given by
\begin{equation*}
  h(0,s)\colon Z_k^-\longrightarrow\IC,\quad (x,a)\longmapsto\abs{a}^{-s}.
\end{equation*}
The number $0$ in $h(0,s)$ refers to the \lq diameter\rq{} of a cusp.
In section \ref{eisenstein} we define functions $h(\ell(k),s)\colon
Z_k^-\to\IC$ that depend on the length $\ell(k)>0$ of the closed
geodesic in $Z_k$. The corresponding notion of approximate Eisenstein
functions on $M_\ell$ then consists of a meromorphic family of
functions associated with each half-cylinder $Z_k^{\pm}\subset Z_k$.
Up to a jump discontinuity on the respective closed geodesic, they are
eigenfunctions of the Laplacian on $M_\ell$.  In this introduction
they shall be denoted by $E_i(s)$, where $i$ runs through the set
$\tilde S$ of half-cylinders.

In this context, the classical scattering matrix is replaced with a
pair of matrices $(C_{ij}(s))_{i,j\in\tilde S}$ and
$(D_{ij}(s))_{i,j\in\tilde S}$ with the following properties: On each
half-cylinder $Z_j^-$ we may define fibrewise Fourier coefficients
$F_j^n$, and the approximate Eisenstein functions $E_i(s)$ satisfy
\begin{equation*}
  F_j^0E_i(s)= D_{ij}(s)\cdot h(\ell(j),s)+ C_{ij}(s)\cdot h(\ell(j),1-s).
\end{equation*}
Then, in consequence of convergence of the resolvent, we obtain
\begin{intro}
  \begin{enumerate}
  \item As $\ell\to 0$, the approximate Eisenstein function $E_i(s)$
    on $M_\ell$ converge to Eisenstein functions for $M_0$ if\/
    $\re(s)>1/2$.
  \item If\/ $\re(s)>1/2$, the matrix $(C_{ij}(s))_{ij}$\/ converges to
    the scattering matrix of $M_0$.
  \end{enumerate}
\end{intro}
Let us mention that the first part of this theorem in particular
implies that the Eisenstein functions $E_i(s)$ on $M_\ell$, defined if
the latter has cusps, converge to Eisenstein functions of the limit if
$\re(s)>1/2$. Results of this kind were used by Obitsu \cite{obitsu}
to study the geometry of Teichm\"uller spaces.

But the approximate Eisenstein functions do not immediately
accomplish the problem of extending the convergence results to
$\{s\in\IC\st \re(s)<1/2\}$. Rather, they satisfy a Mass-Selberg
relation that gives rise to functional equations for surfaces of
finite area.  These can be used to prove the following assertion. Here
$E(s)$ denotes the column vector that has the approximate Eisenstein
functions as entries.
\begin{intro2} Assume that the surfaces $M_\ell$ are of finite area. For $\ell$ near $0$, the meromorphic family of matrices $D\colon s\mapsto (D_{ij}(s))$ is meromorphically invertible. Then $D(s)^{-1}\cdot E(s)$ and $D(s)^{-1}\cdot C(s)$ converge on $\{s\in\IC\st \re(s)\ne 1/2\}$ to the Eisenstein functions and the scattering matrix of $M_0$, respectively.
\end{intro2}

The corollary suggests a replacement for
the quotient $\Zeta(s)/\prod_i\Zeta_{c_i}(s)$ in the first theorem, namely
\begin{equation*}
  \det D(s)\cdot \Zeta(s)/\textstyle\prod_i\Zeta_{c_i}(s).
\end{equation*}
The additional factor does not alter the limit if $\re(s)>1/2$, and we
prove in Theorem \ref{zeta everywhere} that this expression also
converges to the Zeta function of $M_0$ if $\re(s)<1/2$ for a
degenerating family of compact surfaces. It is in this sense, that the
approximate scattering data admit an extension of the first theorem to
the left of the critical axis. Unfortunately we do not know how this
term behaves on the critical axis itself.

The text is arranged as follows:
In section \ref{surface geometry}, we give an explicit description of
the metric degeneration in terms of Fenchel-Nielsen coordinates. It
comes with a convenient choice of coordinates for elementary cylinders
that are embedded in the surfaces, and these coordinates are the basis
for our comparison of integral kernels later on.
Section \ref{resolvent} is divided into two parts. After a few remarks
on compact surfaces, the first part uses the resolvent kernel of the
Laplacian on the hyperbolic plane to compare the kernels for elementary
surfaces of different diameters. The main difficulty is to obtain
trace class estimates that carry over to general surfaces of finite geometry. The second
part applies meromorphic Fredholm theory to examine these surfaces.
Approximate Eisenstein series are described in section
\ref{eisenstein}, and in section \ref{zeta function} we apply our results to
the Selberg Zeta function.

From the year 2000 on I was a member of the DFG research group
``Zetafunktionen und lokalsymmetrische R\"aume'', located at
Clausthal-Zellerfeld and G\"ottingen, and of the Graduate Program
``Gruppen und Geometrie''. I am very grateful to the members of these
groups for their support. But first and foremost I want to express my
thanks to Prof.~Dr.~Ulrich~Bunke. He raised my interest in this
subject, and he was of great influence on me during the past years. I
also want to thank Dr.~Martin Olbrich and Dr.~Margit R\"osler for
their help.

\section{Surface geometry}\label{surface geometry}
A collection of $3 p-3$ incontractible, simple, closed curves that are
disjoint and homotopically distinct may be used to dissect a closed, oriented
surface of genus $p\ge 2$ into pairs of pants. Fenchel-Nielsen
coordinates as gluing data for hyperbolic pairs of pants determine a
hyperbolic surface of the same genus. Such a decomposition of hyperbolic
surfaces will be applied to examine their degeneration in
the sections to come.

To introduce some notation, we begin with giving an explicit
construction for a family of extended pairs of pants. These are
extended in the sense that they are not compact, but the removal of
infinite cylinders will result in the common pairs or pants, which have
a boundary consisting of three closed geodesics. The family provides
the building blocks for general surfaces of finite geometry, and it
also serves as the fundamental example for the degeneration process.

\subsection{A degenerating family of pairs of pants}\label{subs pairs}
Topologically, the
surfaces are thrice punctured spheres, so their fundamental group is the free
group on two generators $\IZ *\IZ$. To provide such a one carrying a
hyperbolic structure, we specify the generators of a Fuchsian group that is
isomorphic to $\IZ *\IZ$. Up to conjugation within the orientation preserving
isometries of the hyperbolic plane, the Fuchsian group is determined by a
triple of three nonnegative numbers, of which every positive stands for the
length of a closed geodesic. In case at least one of the lengths vanishes, we
speak of a degeneration.

For the time being, we use the unit disc model $D$ of the hyperbolic
plane. Hyperbolic distance between $z_1, z_2\in D$ is denoted by
$\dist{z_1}{z_2}$. To determine a conjugacy class for the groups
to be constructed, fix three distinct points that occur in the order $v_1$,
$v_2$, $v_3$ on $\partial D$, the indices will be taken as elements of $\IZ
/ 3\IZ$. Let $L_i$ denote the geodesic line that joins $v_{i+1}$
with $v_{i+2}$.
\begin{figure}
\begin{center}
\begin{minipage}{6cm}
\begin{picture}(0,0)%
\includegraphics{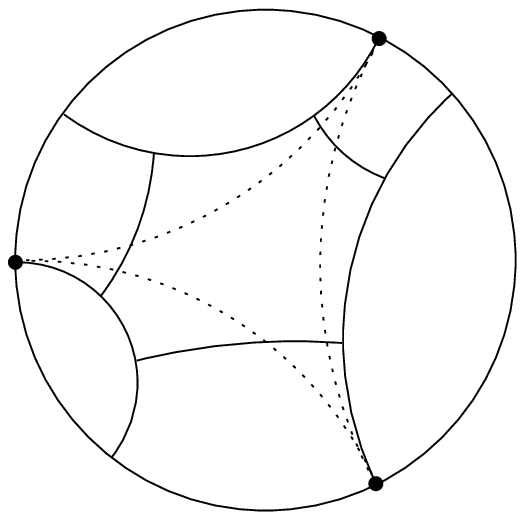}%
\end{picture}%
\setlength{\unitlength}{4144sp}%
\begingroup\makeatletter\ifx\SetFigFont\undefined
\def\x#1#2#3#4#5#6#7\relax{\def\x{#1#2#3#4#5#6}}%
\expandafter\x\fmtname xxxxxx\relax \def\y{splain}%
\ifx\x\y   
\gdef\SetFigFont#1#2#3{%
  \ifnum #1<17\tiny\else \ifnum #1<20\small\else
  \ifnum #1<24\normalsize\else \ifnum #1<29\large\else
  \ifnum #1<34\Large\else \ifnum #1<41\LARGE\else
     \huge\fi\fi\fi\fi\fi\fi
  \csname #3\endcsname}%
\else
\gdef\SetFigFont#1#2#3{\begingroup
  \count@#1\relax \ifnum 25<\count@\count@25\fi
  \def\x{\endgroup\@setsize\SetFigFont{#2pt}}%
  \expandafter\x
    \csname \romannumeral\the\count@ pt\expandafter\endcsname
    \csname @\romannumeral\the\count@ pt\endcsname
  \csname #3\endcsname}%
\fi
\fi\endgroup
\begin{picture}(2517,2492)(21,-1801)
\put(1827,-717){\makebox(0,0)[lb]{\smash{\SetFigFont{12}{14.4}{rm}{$\scriptstyle T_1$}%
}}}
\put(1203,-773){\makebox(0,0)[lb]{\smash{\SetFigFont{12}{14.4}{rm}{$\scriptstyle L_2$}%
}}}
\put(1090, 20){\makebox(0,0)[lb]{\smash{\SetFigFont{12}{14.4}{rm}{$\scriptstyle T_2$}%
}}}
\put(543,-1005){\makebox(0,0)[lb]{\smash{\SetFigFont{12}{14.4}{rm}{$\scriptstyle T_3$}%
}}}
\put(561,-409){\makebox(0,0)[lb]{\smash{\SetFigFont{12}{14.4}{rm}{$\scriptscriptstyle\frac 12 \ell_1$}%
}}}
\put(1189,-1110){\makebox(0,0)[lb]{\smash{\SetFigFont{12}{14.4}{rm}{$\scriptscriptstyle\frac 12\ell_2$}%
}}}
\put( 21,-613){\makebox(0,0)[lb]{\smash{\SetFigFont{12}{14.4}{rm}{$\scriptstyle v_3$}%
}}}
\put(1868,-1743){\makebox(0,0)[lb]{\smash{\SetFigFont{12}{14.4}{rm}{$\scriptstyle v_1$}%
}}}
\put(1927,535){\makebox(0,0)[lb]{\smash{\SetFigFont{12}{14.4}{rm}{$\scriptstyle v_2$}%
}}}
\put(1740,  1){\makebox(0,0)[lb]{\smash{\SetFigFont{12}{14.4}{rm}{$\scriptscriptstyle\frac 12\ell_3$}%
}}}
\put(1225,-442){\makebox(0,0)[lb]{\smash{\SetFigFont{12}{14.4}{rm}{$\scriptstyle L_1$}%
}}}
\put(1448,-583){\makebox(0,0)[lb]{\smash{\SetFigFont{12}{14.4}{rm}{$\scriptstyle L_3$}%
}}}
\end{picture}

\end{minipage}
\hspace{0.5cm}
\begin{minipage}{5cm}
\begin{picture}(0,0)%
\includegraphics{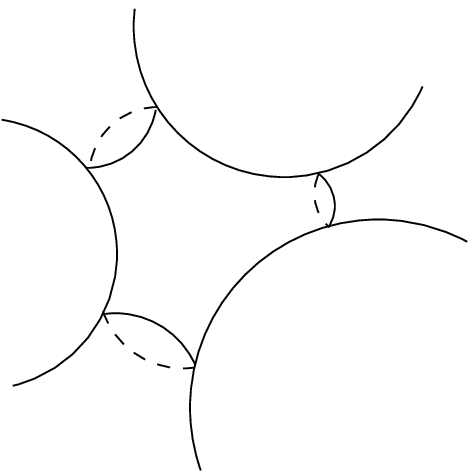}%
\end{picture}%
\setlength{\unitlength}{4144sp}%
\begingroup\makeatletter\ifx\SetFigFont\undefined
\def\x#1#2#3#4#5#6#7\relax{\def\x{#1#2#3#4#5#6}}%
\expandafter\x\fmtname xxxxxx\relax \def\y{splain}%
\ifx\x\y   
\gdef\SetFigFont#1#2#3{%
  \ifnum #1<17\tiny\else \ifnum #1<20\small\else
  \ifnum #1<24\normalsize\else \ifnum #1<29\large\else
  \ifnum #1<34\Large\else \ifnum #1<41\LARGE\else
     \huge\fi\fi\fi\fi\fi\fi
  \csname #3\endcsname}%
\else
\gdef\SetFigFont#1#2#3{\begingroup
  \count@#1\relax \ifnum 25<\count@\count@25\fi
  \def\x{\endgroup\@setsize\SetFigFont{#2pt}}%
  \expandafter\x
    \csname \romannumeral\the\count@ pt\expandafter\endcsname
    \csname @\romannumeral\the\count@ pt\endcsname
  \csname #3\endcsname}%
\fi
\fi\endgroup
\begin{picture}(2143,2126)(12,-1635)
\put(397,  7){\makebox(0,0)[lb]{\smash{\SetFigFont{12}{14.4}{rm}{$\scriptstyle \ell_1$}%
}}}
\put(545,-1243){\makebox(0,0)[lb]{\smash{\SetFigFont{12}{14.4}{rm}{$\scriptstyle\ell_2$}%
}}}
\put(1595,-413){\makebox(0,0)[lb]{\smash{\SetFigFont{12}{14.4}{rm}{$\scriptstyle\ell_3$}%
}}}
\end{picture}

\end{minipage}
\end{center}
\caption{Construction of a pair of pants with prescribed circumference of funnels.}\label{pant}
\end{figure}
\begin{lem}\label{hex} Any triple $(\ell_1,\ell_2,\ell_3)$ of non-negative reals determines disjoint geo\-de\-sics $T_1, T_2, T_3$ $\subset D$ with the following properties  (cf.~figure \ref{pant}).
  \begin{enumerate}
  \item \label{order} Each $T_i$ meets $\partial D$ at $v_i$ and the geodesic $L_{i-1}$ separates $T_i$ from $L_i$.
  \item The number $\ell_i/2$ is equal to the hyperbolic distance $\dist{T_{i+1}}{T_{i+2}}$.
  \end{enumerate}
  In fact, $T_i$ is the unique geodesic with the first property that satisfies
  \begin{equation}\label{distance hexagon}
    \cosh (\dist{T_i}{L_i}) = \left( \cosh (\ell_i/2) +1\right)^{-1} \left[ m(\ell_1,\ell_2,\ell_3) +\cosh (\ell_{i+2}/2) -\cosh(\ell_{i+1}/2)\right]
  \end{equation}
  where
  \begin{equation*}\begin{split}
      m(\ell_1,\ell_2,\ell_3)&=\left( \cosh^2(\ell_1/2)+ \cosh^2(\ell_2/2)+
        \cosh^2(\ell_3/2)\right.\\
      &\qquad \bigl.+2 \cosh( \ell_1/2) \cosh(\ell_2/2) \cosh(\ell_3/2)
      -1\bigr)^{1/2}.
    \end{split}\end{equation*}
\end{lem}
\begin{proof} It is convenient to consider this as a statement on the
  inversive product of spheres in the extended plane $\hat\IR^2$. The set
  of spheres is identified with a subset of $\IR P^3$ by mapping $\bigl\{
  x\in\IR^2 \st a_0 \norm x^2 -2 \scal x{(a_1,a_2)} + a_3=0\bigr\}$ to the
  equivalence class of $(a_0,a_1,a_2,a_3)$. With the bilinear form
  $q(a,b)= 2(a_1b_1+a_2b_2)- a_0b_3- a_3 b_0$, the inversive product of two
  spheres $S$ and $T$, represented by $a,b\in \IR P^3$, is defined to be
  (cf.~Beardon \cite[p.~28{\em ff\ }]{beardon}) 
  \begin{equation*} \inv ST = \frac{\abs{q(a,b)}}{\abs{q(a,a)}^{1/2} \abs{q(b,b)}^{1/2}}.
  \end{equation*}
  Now, given spheres $L_1$, $L_2$, $L_3$ and $T_1$, $T_2$, $T_3$ as in the proposition with property \ref{order}, one can compute $\inv{T_i}{L_i}$ from $\inv{T_i}{T_{i+1}}$ and $\inv{T_{i+1}}{L_{i+1}}$:
  \begin{equation*}
    \inv{T_i}{L_i} = \frac{ 2\left( \inv{T_i}{T_{i+1}} +\inv{T_{i+1}}{L_{i+1}}\right)}{ 1+\inv{T_{i+1}}{L_{i+1}}} -1.
  \end{equation*}
  The formula can be verified conveniently in the upper half plane with
  $\left\{ v_1,v_2,v_3\right\}= \left\{0,1,\infty\right\}$ since the
  inversive product is invariant under M\"obius transformations.
  
  This gives a system of three equations, solved by the (non-negative)
  inversive products $\inv{T_i}{L_i}$ if and only if
  \begin{equation*}\begin{split}
      \inv{T_i}{L_i} &= \left(\inv{T_{i+1}}{T_{i+2}}+1\right)^{-1}\cdot\Bigl[\inv{T_i}{T_{i+1}}-\inv{T_i}{T_{i+2}}\Bigr.\\
        &\quad\Bigl. +\sqrt{ \inv{T_1}{T_2}^2+ \inv{T_2}{T_3}^2+ \inv{T_3}{T_1}^2 + 2\inv{T_1}{T_2} \inv{T_2}{T_3} \inv{T_3}{T_1} -1}\Bigr].
  \end{split}\end{equation*}
  Finally, if the intersections of two spheres $S$,$T$ with $D$ are disjoint geodesics, we have $\inv ST=\cosh \dist{S\cap D}{T\cap D}$. So the previous equation calculates $\inv{T_i}{L_i}$ from $\ell_1$, $\ell_2$, $\ell_3$, and this in turn determines $T_i$.
\end{proof}
In view of this lemma, the set $\base_0:=[0,\infty)^3$ now serves as
parameter space for Fuchsian groups. Each point
$(\ell_1,\ell_2,\ell_3)\in \base_0$ gives rise to geodesics $T_i$. Let
$\sigma_i$ denote reflection in $T_i$ and
$\gamma_i=\sigma_{i+2}\sigma_{i+1}$. Then $\gamma_i$, $\gamma_{i+1}$
act as side-pairing transformations on the convex polygon bounded by
$T_i$, $T_{i+1}$, $\sigma_{i+2}T_i$ and $\sigma_{i+2}T_{i+1}$\label{fdom}.
Observe that $\gamma_i$ is a parabolic isometry of $D$ if $T_{i+2}$
and $T_{i+1}$ meet in $v_{i+2}$, i.e.~if $\ell_i=0$. Otherwise it is a
hyperbolic transformation of translation length $\ell_i$, its axis
being the common orthogonal of $T_{i+2}$ and $T_{i+1}$. 

The group generated by $\gamma_1$, $\gamma_2$, $\gamma_3$ has a
presentation $\langle\gamma_1,\gamma_2,\gamma_3\st
1=\gamma_3\gamma_2\gamma_1\rangle$, and Poincar\'e's theorem implies
that it is a discrete group of isometries. We fix the isomorphism of
$\IZ*\IZ$ onto this group that maps the natural generators to
$\gamma_1$, $\gamma_2$. Equation \eqref{distance hexagon} shows that
we constructed a family
\begin{equation*}
  \phi\colon\quad \base_0\longrightarrow \hom(\IZ *\IZ,\isom(D))
\end{equation*}
that is smooth even on the boundary of $\base_0$.

Each of the generators $\gamma_1,\gamma_2,\gamma_3$ corresponds to a
cylinder embedded in the quotient $Y_0(\ell):=\phi(\ell)\backslash D$.
General surfaces will be defined by a gluing procedure along the
cylinders, and we need to choose suitable coordinates. For this
purpose, the following models of the hyperbolic plane are used. The
Riemannian metrics are special cases of those employed by Judge
\cite{judge1} in his study of the spectrum of degenerating manifolds
of more general geometry.
\begin{dfn}\label{model}
  For each positive real number $t$, let $X_t$ denote the manifold $\IR^2$
  endowed with the Riemannian metric
  \begin{equation*} 
    (x,a)\longmapsto (t^2+a^2)\dx x^2+(t^2+a^2)^{-1}\dx
    a^2.\end{equation*}
  For $t=0$ we consider the disconnected manifold
  \begin{equation*} X_t=X_0=\left\{(x,a)\in\IR^2\st a\ne 0\right\},\quad
    (x,a)\longmapsto a^2\dx x^2+ a^{-2}\dx a^2.
  \end{equation*}
  We shall also use the notation $X_t^{\pm}=\left\{ (x,a)\in X_t\st
    \pm a > 0\right\}$, and more generally, for arbitrary subsets
  $A\subset\IR$, the set $\left\{(x,a)\in X_t\st a\in A\right\}$ is
  denoted by $X_t^A$.
\end{dfn}
The plane $X_t$, if $t\ne 0$, is mapped isometrically onto the upper half-plane model for the hyperbolic plane by 
\begin{equation*}
  (x,a)\longmapsto e^{t x} \left(\frac a{\sqrt{t^2+a^2}},\, \frac{t}{\sqrt{t^2+a^2}}\right),
\end{equation*}
and so is $X_0^{\pm}$ by $(x,a)\mapsto (x,\pm a^{-1})$. The map 
\begin{equation*}
\gamma\colon\  X_t\longrightarrow X_t,\quad (x,a)\longmapsto (x+1,a)
\end{equation*}
corresponds to a hyperbolic isometry of translation $t$ in the
first case and to a parabolic one in the second. So the quotient
$Z_{t}:=\langle\gamma\rangle\backslash X_{t}$ is either an
elementary hyperbolic cylinder or the disjoint union of two cusps.

It is convenient to denote subsets of $Z_t$ by $Z_t^\pm$ and $Z_t^A$
as in definition \ref{model}. The following lemma is essentially the
collar lemma from hyperbolic geometry. The definition of $A(t)$ therein is only a preliminary, cf.~\eqref{A(t)}.
\begin{lem}\label{identif} Let $\ell=(\ell_1,\ell_2,\ell_3)\in\base_0$, and for all $t\ge 0$ we put
  \begin{equation*}
    A(t)=\begin{cases}
      \bigl( -\tfrac{t}{2 \sinh(t/2)}, \infty\bigr)& \text{if
        $t> 0$},\\
      \left(-1,0\right)& \text{if $t=0$.}\end{cases}
  \end{equation*}
  Then the construction of lemma \ref{hex} gives rise to canonical embeddings of cylinders
  \begin{equation*}
    \Psi_i\colon Z_{\ell_i}^{A(\ell_i)}\longrightarrow Y_0(\ell)
  \end{equation*}
  such that the images $C_i(\ell):=\Psi_i\bigl(Z_{\ell_i}^{A(\ell_i)}\bigr)$ are
  mutually disjoint.
\end{lem}
\begin{proof} Any pair $(S_1,S_2)$ of disjoint geodesics with positive distance $t/2$ determines an orientation-preserving isometry of $X_t$ onto $D$. It maps $\left\{ (x,a)\st x=0\right\}$ to $S_1$ and $\left\{(x,a)\st x=1/2\right\}$ to $S_2$. In case the geodesics are disjoint with distance $t/2=0$, it may be necessary to interchange $S_1$ and $S_2$, but then there is a unique isometry $X_0^-\to D$ with the same property. We apply this to the geodesics that are associated with $(\ell_1,\ell_2,\ell_3)$ by Lemma \ref{hex}: Each pair $(T_{i+1},T_{i+2})$ of geodesics gives rise to $\tilde\Psi_i\colon X_{\ell_i}\rightarrow D$ (or $\tilde\Psi_i\colon X_{0}^-\rightarrow D$, respectively). The set $\tilde\Psi_i\bigl(X_{\ell_i}^{A(\ell_i)}\bigr)\subset D$ is precisely invariant under the action of $\langle\gamma_i\rangle\subset \IZ *\IZ$ (this is one way of proving the collar lemma). Then the $\tilde\Psi_i$ induce maps $\Psi_i$ between the quotients as proposed.
\end{proof}
So we defined a connected hyperbolic surface $Y_0(\ell)$ for each
$\ell\in\base_0$ with distinguished subsets $C_1,C_2,C_3$. Each $C_i$
is diffeomorphic to a cylinder $Z_{\ell_i}^{A(\ell_i)}$, and a point
$p\in C_i$ ist given, via the map $\Psi_i$, by coordinates $p=(x,a)$
where $x\in\IR$ is well-defined modulo $\IZ$. The image of
$Z_{\ell_i}^+\cap Z_{\ell_i}^{A(\ell_i)}$ under $\Psi_i$ is
\begin{equation*}
  C_i^+=\{(x,a)\in C_i\st a>0\},
\end{equation*}
and the interior $P(\ell)$ of the set $Y_0(\ell)\setminus
\bigcup_{i=1}^3 C_i^+$ is an open pair of pants in the customary sense
of this word. Note that $C_i^+$ is empty by definition if
$\ell_i=0$.

\subsection{Assembling a surface from pairs of pants}\label{assembling} 
\begin{figure}
  \begin{center}
\begin{picture}(0,0)%
\includegraphics{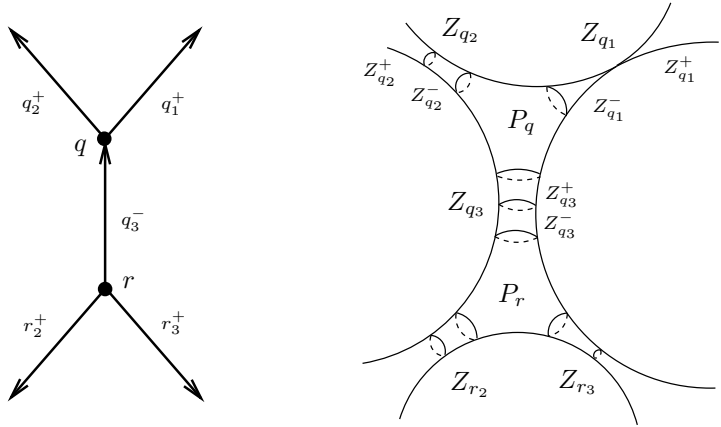}%
\end{picture}%
\setlength{\unitlength}{4144sp}%
\begingroup\makeatletter\ifx\SetFigFont\undefined
\def\x#1#2#3#4#5#6#7\relax{\def\x{#1#2#3#4#5#6}}%
\expandafter\x\fmtname xxxxxx\relax \def\y{splain}%
\ifx\x\y   
\gdef\SetFigFont#1#2#3{%
  \ifnum #1<17\tiny\else \ifnum #1<20\small\else
  \ifnum #1<24\normalsize\else \ifnum #1<29\large\else
  \ifnum #1<34\Large\else \ifnum #1<41\LARGE\else
     \huge\fi\fi\fi\fi\fi\fi
  \csname #3\endcsname}%
\else
\gdef\SetFigFont#1#2#3{\begingroup
  \count@#1\relax \ifnum 25<\count@\count@25\fi
  \def\x{\endgroup\@setsize\SetFigFont{#2pt}}%
  \expandafter\x
    \csname \romannumeral\the\count@ pt\expandafter\endcsname
    \csname @\romannumeral\the\count@ pt\endcsname
  \csname #3\endcsname}%
\fi
\fi\endgroup
\begin{picture}(4286,2542)(335,-1762)
\put(3348, 25){\makebox(0,0)[lb]{\smash{\SetFigFont{10}{12.0}{rm}{$P_q$}%
}}}
\put(4283,350){\makebox(0,0)[lb]{\smash{\SetFigFont{10}{12.0}{rm}{$\scriptstyle Z_{q_1}^+$}%
}}}
\put(2953,575){\makebox(0,0)[lb]{\smash{\SetFigFont{10}{12.0}{rm}{$Z_{q_2}$}%
}}}
\put(3783,540){\makebox(0,0)[lb]{\smash{\SetFigFont{10}{12.0}{rm}{$Z_{q_1}$}%
}}}
\put(3573,-400){\makebox(0,0)[lb]{\smash{\SetFigFont{10}{12.0}{rm}{$\scriptstyle Z_{q_3}^+$}%
}}}
\put(3563,-585){\makebox(0,0)[lb]{\smash{\SetFigFont{10}{12.0}{rm}{$\scriptstyle Z_{q_3}^-$}%
}}}
\put(3653,-1535){\makebox(0,0)[lb]{\smash{\SetFigFont{10}{12.0}{rm}{$Z_{r_3}$}%
}}}
\put(3858,120){\makebox(0,0)[lb]{\smash{\SetFigFont{10}{12.0}{rm}{$\scriptstyle Z_{q_1}^-$}%
}}}
\put(2768,185){\makebox(0,0)[lb]{\smash{\SetFigFont{10}{12.0}{rm}{$\scriptstyle Z_{q_2}^-$}%
}}}
\put(2483,330){\makebox(0,0)[lb]{\smash{\SetFigFont{10}{12.0}{rm}{$\scriptstyle Z_{q_2}^+$}%
}}}
\put(2988,-470){\makebox(0,0)[lb]{\smash{\SetFigFont{10}{12.0}{rm}{$Z_{q_3}$}%
}}}
\put(3293,-1020){\makebox(0,0)[lb]{\smash{\SetFigFont{10}{12.0}{rm}{$P_r$}%
}}}
\put(3008,-1550){\makebox(0,0)[lb]{\smash{\SetFigFont{10}{12.0}{rm}{$Z_{r_2}$}%
}}}
\put(758,-115){\makebox(0,0)[lb]{\smash{\SetFigFont{10}{12.0}{rm}{$q$}%
}}}
\put(1038,-560){\makebox(0,0)[lb]{\smash{\SetFigFont{10}{12.0}{rm}{$\scriptstyle q_3^-$}%
}}}
\put(1048,-935){\makebox(0,0)[lb]{\smash{\SetFigFont{10}{12.0}{rm}{$r$}%
}}}
\put(453,-1195){\makebox(0,0)[lb]{\smash{\SetFigFont{10}{12.0}{rm}{$\scriptstyle r_2^+$}%
}}}
\put(1278,-1180){\makebox(0,0)[lb]{\smash{\SetFigFont{10}{12.0}{rm}{$\scriptstyle r_3^+$}%
}}}
\put(1278,135){\makebox(0,0)[lb]{\smash{\SetFigFont{10}{12.0}{rm}{$\scriptstyle q_1^+$}%
}}}
\put(443,140){\makebox(0,0)[lb]{\smash{\SetFigFont{10}{12.0}{rm}{$\scriptstyle q_2^+$}%
}}}
\end{picture}
  \end{center}
  \caption{A surface for the graph of type (0,4), $\ell(q_1)=0$.}\label{(0,4)}
\end{figure}
We describe the gluing procedure that composes a geometrically
finite surface from pairs of pants. Our presentation
follows that of Kra's article on horocyclic coordinates~\cite{kra},
with focus on Riemannian geometry rather than complex structure. It
incorporates a choice of local coordinates for embedded cylinders,
which are modelled after lemma \ref{identif}.

Starting point is an \emph{admissible graph $\graph G$ of type
  $(p,n)$}. Basically, it determines the topological
structure of the surface to be constructed, which will be of genus $p$
with $n$ deleted discs in the non-degenerate case. The graph consists
of a set $\graph G_0$ of vertices, a set $\tilde {\graph G}_1$ of
oriented edges, an involution $\iota$ of $\tilde{\graph G}_1$ without
fixed elements, and a map $s\colon\tilde{\graph G}_1\to\graph G_0$
that associates the source to every edge. The set of unoriented edges
$\tilde{\graph G}_1/\iota$ is denoted by $\graph G_1$. The numbers
$p,n\in\IN_0$ satisfy $2p-2+n>0$ and $3p-3+n\ge 0$ by definition. The
number of vertices is $2p-2+n$, the number of unoriented edges is
$3p-3+n$, and the preimage of any $q\in\graph G_0$ under $s$ consists
of three edges at most.

It is convenient to introduce additional edges that do not belong to
the domain of definition of the involution $\iota$. They correspond to
infinite funnels, i.e.~cylinders with only one end attached to a pair
of pants.
\begin{dfn}\label{augmented} An \emph{augmented admissible graph} $\graph G^*$ consists of a set $\graph G_0$ of vertices, a set  $\tilde{\graph G}_1^*$ of oriented edges, a source map $s:\tilde{\graph G}_1^*\to \graph G_0$, and a distinguished subset $\tilde{\graph G}_1\subset\tilde{\graph G}_1^*$ with an involution $\iota$ of $\tilde{\graph G}_1$ such that
  \begin{itemize}
  \item the tuple $\bigl(\tilde{\graph G}_1,\graph
    G_0,s\vert_{\tilde{\graph G}_1},\iota\bigr)$ is an admissible
    graph,
  \item for each vertex $q\in\graph G_0$, the preimage $s^{-1}(q)\subset\tilde{\graph G}_1^*$
    consists of exactly three elements.
  \end{itemize}
  The set $\graph G_1^*$ of \emph{unoriented edges}\/ is the quotient
  of $\tilde{\graph G}_1^*$ under the equivalence relation generated
  by $e\sim\iota(e)$, and the elements of $\tilde{\graph
    G}_1^*\setminus\tilde{\graph G}_1$ are called \emph{phantom
    edges}.
\end{dfn}
If $\graph G^*$ is an augmented admissible graph, we may choose a map
\begin{equation*}
  \graph G_0\longrightarrow \tilde{\graph G}_1^*\times \tilde{\graph G}_1^*\times
  \tilde{\graph G}_1^*,\quad q\longmapsto (q_1,q_2,q_3),
\end{equation*}
that assigns all three adjacent edges to each vertex. Now a surface
$Y_{\graph G}(\lambda)$ with a hyperbolic metric is determined by a
labelling $\lambda$ of the edges of $\graph G^*$ with length and twist
parameters as follows. Let
\begin{equation*}
  \lambda\colon\ \graph G_1^*\longrightarrow \left[0,\infty\right)\times \IR ,\quad
  d\longmapsto(\ell(d),\tau(d))
\end{equation*}
be an arbitrary map. The lift of $\lambda$ to $\tilde {\graph G}_1^*$
is also denoted by $\lambda$. We consider a pair of pants
\begin{equation*}
  P_{q}:=\left\{q\right\}\times P(\ell(q_{1}),\ell(q_{2}),\ell(q_{3}))
\end{equation*}
for each vertex $q\in\graph G_0$ as defined at the end of section \ref{subs pairs}, and an infinite
elementary half-cylinder or a cusp
\begin{equation*}
  Z_{e}:=\left\{e\right\}\times \left(Z_{\ell(e)}^{A(\ell(e))}\cup
  Z_{\ell(e)}^+ \right)
\end{equation*}
for each \emph{phantom} edge $e\in\tilde{\graph
  G}_{1}^{*}\setminus\tilde {\graph G}_{1}$.  These half-cylinders
will be attached to ends of the $P_q$, while additional collars will
serve as connectors between the remaining ends of pairs of pants.  To
define these connectors, we alter the definition of the interval
$A(t)$ from previous subsection into a more symmetric one,
\begin{equation}\label{A(t)}
  A(t):=
  \begin{cases}
    \bigl( -\frac{t}{2\sinh(t/2)},
    \frac{t}{2\sinh(t/2)}\bigr)& \text{if $t\ne 0$},\\
    (-1,1)&\text{if $t=0$}.
  \end{cases}
\end{equation}
This allows for the definition of a collar $Z_e$, for each
\emph{proper}\/ edge $e\in\tilde{\graph G}_1$, by
\begin{equation*}
  Z_e:=\{e\}\times Z_{\ell(e)}^{A(\ell(e))},
\end{equation*}
which is indeed a collar around a closed geodesic if $\ell(e)\ne 0$.

Now we define $Y_{\graph G}(\lambda)$ by an equivalence relation,
where the coordinates we use on $Z_e$ and on $C_i^-\subset P_q$ are
the canonical ones induced by definition~\ref{model} and lemma~\ref{identif}.
\begin{dfn} Let $\base$ be the set of maps $\graph G_1^*\longrightarrow \left[0,\infty\right)\times \IR$. For each $\lambda=(\ell,\tau)\in\base$, let $Y_{\graph G}(\lambda)$ be the surface defined by
  \begin{equation}\label{surface}
    Y_{\graph G}(\lambda) = \biggl( \bigcup_{q\in\graph G_{0}}P_{q}\ \cup
    \bigcup_{e\in\tilde{\graph G}_{1}^{*}}Z_{e}\biggr)/\sim,
  \end{equation}
  with the equivalence relation being generated as follows:
  \begin{align*}
    (x,a)\sim (x-\tau(q_i)/2,a)\qquad
    \text{where}&\quad (x,a)\in C_i^-\subset
    P_q,\\
    &\quad(x-\tau(q_i)/2,a)\in Z_{q_i}^{-}\subset Z_{q_i};\\
    (x,a)\sim (-x,-a)\qquad
    \text{where}&\quad (x,a)\in Z_e,\ (-x,-a)\in Z_{\iota(e)}.
  \end{align*}
\end{dfn}
The first relation glues a cylinder to each end $C_i^-$ of a pair of
pants $P_q$, and the second identifies two such cylinders $Z_e$,
$Z_{\iota(e)}$ according to the structure of the graph. Recall that
the involution $\iota$ is defined on the \emph{proper}\/ edges in
$\tilde{\graph G}_1$ only, so the condition for the second generator
is never satisfied if $e$ is a phantom edge.

From now on, the notation $P_q$ and $Z_e$ will refer to the respective
images in the quotient $Y_{\graph G}(\lambda)$. Then
$Z_e=Z_{\iota(e)}$ for each proper edge $e$, and we may speak of a
subset $Z_d\subset Y_\graph G(\lambda)$ for each unoriented edge
$d\in\graph G_1^*$. Note the following: The choice of an orientation
$e\in\tilde {\graph G}_1^*$ of $d$ endows the cylinder $Z_d$ with a
canonical quotient structure $\langle\gamma\rangle\backslash
X_{\ell(e)}^A$. Replacing $e$ with $\iota(e)$ corresponds to the
coordinate change $(x,a)\mapsto (-x,-a)$. In particular, we may use
oriented edges to specify half-cylinders $Z_e^\pm\subset Z_d$, and
$Z_e^\pm=Z_{\iota(e)}^\mp$ holds. This will be of importance in the
definition of approximate Eisenstein series in section
\ref{eisenstein}.

The $P_q$ and $Z_d$ provide an open cover of $Y_{\graph G}(\lambda)$
as illustrated in figure \ref{(0,4)}. If the admissible graph is of
type $(0,3)$, then the extended pairs of pants $Y_0(\ell)$, as defined
in the previous paragraph, occur as the non-elementary component of
$Y_{\graph G}(\lambda)$, and $Y_{\graph G}(\lambda)$ is not connected
if a length $\ell(d)$ vanishes. In this instance, the edge $d$ gives
rise to an isolated cusp $Z_d^{+}\subset Y_{\graph G}$. The reader is
advised to have a look at Kra's text \cite{kra} for a discussion of
several graph types.

The quotient $Y_{\graph G}(\lambda)$ inherits a hyperbolic metric.
Obviously, there exists an isometry between the emerging surfaces if a
twist $\tau(d)$ is replaced with $\tau(d)+1$ (or with $\tau(d)+\delta$
for arbitrary $\delta$ if $d\in\graph G_{1}^{*}\setminus\graph
G_{1}$).

\subsection{Choice of a trivialisation}
In section \ref{resolvent} we will need maps between surfaces
$Y_{\graph G}(\lambda)$ for different values of $\lambda$, and these
will be constructed here. In definition \ref{model}, the identity
mapping on $\IR^{2}$ descends to diffeomorphisms between the surfaces
$Z_t^+\cup Z_t^-$ as $t$ varies in $\left[0,\infty\right)$. We want
our maps to coincide with these canonical maps on the cylinders
embedded in $Y_{\graph G}(\lambda)$. What we do is to define such maps
for each pair of pants $P_q$ (and possibly its adjacent infinite
half-cylinders) separately.  They do not necessarily extend to a
diffeomorphism near a closed geodesic where two pairs of pants meet.
But this is sufficient for our needs as it yields operators on
$\L^2$-spaces that are strongly continous with respect to the
geometry.

So we start with choosing suitable diffeomorphisms for the extended
pairs of pants $Y_0(\ell)$ of subsection \ref{subs pairs}. Recall that
we have $\base_0=[0,\infty)^3$, and $\phi$ maps $\base_0$ to the
monomorphisms of $\IZ*\IZ$ into the isometries of the hyperbolic plane
$D$. There is a fibre space in analogy with the Bers fibre space in
Teichm\"uller theory:
\begin{dfn} Let $\tilde\total_0:=\base_0\times D$ and $\tilde\rho\colon \tilde\total_0\to\base_0$ be the canonical map. Then we define $\total_0:=(\IZ *\IZ)\backslash\tilde \total_0$ and $\rho\colon\total_0\to\base_0$ to be the map induced by $\tilde\rho$, where $\IZ *\IZ$ acts on $\tilde\total_0$ by $\gamma(\ell,z)=(\ell,\phi(\ell)(\gamma) z)$.
\end{dfn}
The $\IZ *\IZ$-action on $\tilde\total_0$ is smooth, and we refer to
lemma \ref{discont} for a proof that it is freely discontinuous. So
$\total_0$ is a smooth manifold. The group acts on each fibre of
$\tilde\rho$, so $\rho$ is well-defined, and $\rho^{-1}(\ell)\cong Y_0(\ell)$.

The aim is to equip $\rho$ locally with a trivialisation that exhibits
a specific behaviour on the infinite ends of the fibres. This
behaviour is modelled on the canonical map induced by the identity on
$Z_t^+\cup Z_t^-$ by means of lemma \ref{identif}:

Let us fix $j\in \{1,2,3\}$ for a moment. Differentiation of each
map $\Psi_i$ with respect to $\ell_j$ defines a smooth vector field
$s_i^j$ on the open subset
\begin{equation*}
  \bigcup_{\ell}\{\ell\}\times C_i(\ell)\subset \total_0
\end{equation*}
such that $s_i^j-\partial_{\ell_j}$ is tangent to the fibre of $\rho$.
Then discontinuity of the $\IZ*\IZ$-action on $\tilde\total_0$ implies
the existence of a smooth vector field $s^j$ on $\total_0$ such that
$s^j-\partial_{\ell_j}$ is again vertical with respect to $\rho$, the
restriction of $s^j$ to each subset $\{\ell\}\times C_i(\ell)$ being
equal to $s_i^j$.

If we agree on the order, in which the integral flows of $s^1$, $s^2$
and $s^3$ are to be applied, they can be used to ``trivialise''
$\rho$. With reference to appendix \ref{families} for the details, we
only state the outcome in the proposition below. Here
$P(\ell)\subset\rho^{-1}(\ell)$ again denotes the embedded pair of
pants.
\begin{prop}\label{triv} For each $\ell\in\base_0$ there
  exists a neighbourhood $\inbase$ of $\ell$ and a smooth map
\begin{equation}\label{psi}
  \Psi\colon\ \inbase\times P(\ell)\longrightarrow\total_0
\end{equation}
with the following properties:
\begin{itemize}
\item For all $\ell'\in\inbase$, the set $P(\ell')\subset\rho^{-1}(\ell')$ is
  the diffeomorphic image of $\{\ell'\}\times P(\ell)$.
\item Via the identification of $C_i(\ell')$ with
  $Z_{\ell'_i}^{A(\ell'_i)}$ in lemma \ref{identif}, the map $$\Psi(\ell',\cdot)\colon
  \ C_i(\ell)\longrightarrow P(\ell')$$ corresponds to the canonical identification of
  $Z_{\ell_i}^I$ with $Z_{\ell'_i}^I$, where $I$ is the intersection of
  $A(\ell_i)$ with $A(\ell'_i)$. In particular, this restriction to $Z_{\ell_i}^I$ is area-preserving,
  since the volume form on $X_{\ell'}$ is the Euclidean one for all $\ell'$.
\item If the degeneracies of $\rho^{-1}(\ell')$ coincide with those of $\rho^{-1}(\ell)$, that is to say, if $\ell'_i=0$ implies $\ell_i=0$ and vice versa, then the
  integral flows extend this map to a diffeomorphism $\Psi(\ell')\colon \rho^{-1}(\ell)\to \rho^{-1}(\ell')$.
\end{itemize}
\end{prop}
Now let $\graph G$ be an arbitrary admissible graph. To compare
$Y_{\graph G}(\lambda')$ with $Y_{\graph G}(\lambda)$, we apply the
trivialising maps from proposition \ref{triv} to each $P_q$ with its
adjacent cylinders $Z_{q_i}^-$ separately. Thus, if $U(\lambda)\subset
Y_{\graph G}(\lambda)$ is the complement of all closed geodesics that
are associated with elements of $\graph G_1^*$, there is a
distinguished map $U(\lambda)\rightarrow Y_{\graph G}(\lambda')$. This
map is a diffeomorphism onto the corresponding set $U(\lambda')$. The
pull-back metric has a unique extension to $Y_{\graph G}(\lambda)$ if
there are no additional degeneracies in $Y_{\graph G}(\lambda')$, and
these metrics match up to form a continuous family with respect to
$\lambda'$.

\subsection{The Selberg Zeta function in its domain of convergence}\label{zetaconv}
The result proved here is a prerequisite for section \ref{zeta
  function}. It only requires the construction from the first two
sections, so we include this assertion here. It makes use of a
standard procedure to estimate the number of prime closed geodesics in
a hyperbolic surface that also proves convergence of the infinite
product in the definition of the Zeta function. The argument is
essentially due to Hejhal \cite{hej3}, but we skip the notion of
regular b-groups and consider arbitrary surfaces of finite geometry.

Let $\graph G$ be an admissible graph of type $(p,n)$ and $\graph G^*$
its augmentation.  The graph will remain fixed. Let $\Zeta_\lambda$
denote the Selberg Zeta function of the surface $Y(\lambda)=Y_\graph
G(\lambda)$, which is the product of the Zeta functions of all
connected components. The contribution of an elementary cusp to
$\Zeta_\lambda$ is trivial by definition.

For each $d\in \graph G_1^*$ and $\lambda=(\ell,\tau) \in\base$ with
$\ell(d)\ne 0$, an entire function $\Zeta_{d,\lambda}$ is defined by
the infinite product
\begin{equation}\label{zetaquot}
  \Zeta_{d,\lambda}\colon\ s\longmapsto\prod_{k=0}^\infty\left(
  1-e^{-(s+k)\ell(d)}\right)^2.
\end{equation}
This product is absolutely convergent. We extend this definition to
$\ell(d)=0$ by putting $\Zeta_{d,\lambda}= 1$ in that case. The
purpose of this section is to prove the following assertion.
\begin{prop}\label{prop zeta}  With respect to the topology of locally uniform
  convergence of analytic functions on $\{s\in\IC\st \re(s)> 1\}$, the
  map $\lambda\mapsto\Zeta_\lambda/\prod_{d\in\graph
    G_1^*}\Zeta_{d,\lambda}$ is continuous on $\base$.
\end{prop}
The basic idea is the following elementary estimate: If $Y$ is a
compact, connected hyperbolic surface, let $N(r)$ be the number of
closed geodesics on $Y$ of length at most $r$. Then there exists $C>0$
such that $N(r)\le Ce^r$ holds for all $r>0$. The constant $C$ is
given explicitly in terms of area and diameter of a fundamental
domain for a uniformising group.

In our situation we must not assume that a family of Fuchsian groups
has fundamental domains of uniformly bounded diameter as cusps emerge.
But each closed geodesic that is not entirely contained in a cylinder
must meet a certain fixed subset of the surface, and it is sufficient
that this subset has uniformly bounded diameter. To make this precise,
we choose uniformising groups and fundamental domains for all
$Y_{\graph G}(\lambda)$.

As $Y_{\graph G}(\lambda)$ may have several non-elementary components,
it is convenient to have one group $\Gamma_q (\lambda)$ at hand, for
each vertex $q\in \graph G_0$, that is a uniformizing group for the
component that contain a pair of pants $P_q$. So let $q\in\graph G_0$
be a vertex with the adjacent edges $q_1,q_2,q_3\in\tilde{\graph
  G}_1^*$. For each $\lambda=(\ell,\tau)\in\base$ we defined in
\ref{subs pairs} geodesics $T_1, T_2, T_3\subset D$ such that their
distances satisfy
\begin{equation*}
  \dist{T_{i+1}}{T_{i+2}}=\ell(q_i)/2,
\end{equation*}
and a Fuchsian group $\Delta_q=\phi(\ell(q_1),\ell(q_2),\ell(q_3))$
with generators $\gamma_1,\gamma_2,\gamma_3$.
\begin{figure}
\begin{center}
\begin{minipage}{6cm}
\begin{picture}(0,0)%
\includegraphics{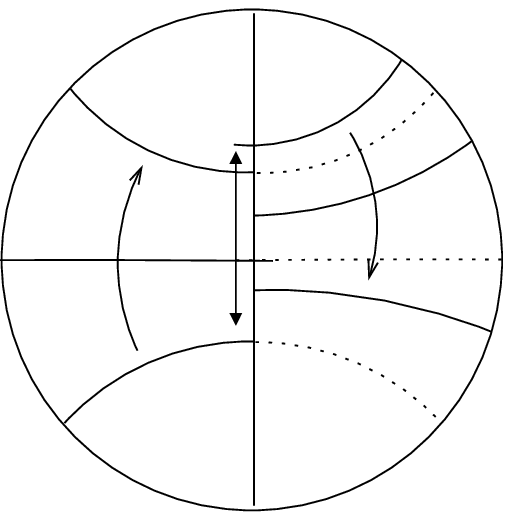}%
\end{picture}%
\setlength{\unitlength}{4144sp}%
\begingroup\makeatletter\ifx\SetFigFont\undefined
\def\x#1#2#3#4#5#6#7\relax{\def\x{#1#2#3#4#5#6}}%
\expandafter\x\fmtname xxxxxx\relax \def\y{splain}%
\ifx\x\y   
\gdef\SetFigFont#1#2#3{%
  \ifnum #1<17\tiny\else \ifnum #1<20\small\else
  \ifnum #1<24\normalsize\else \ifnum #1<29\large\else
  \ifnum #1<34\Large\else \ifnum #1<41\LARGE\else
     \huge\fi\fi\fi\fi\fi\fi
  \csname #3\endcsname}%
\else
\gdef\SetFigFont#1#2#3{\begingroup
  \count@#1\relax \ifnum 25<\count@\count@25\fi
  \def\x{\endgroup\@setsize\SetFigFont{#2pt}}%
  \expandafter\x
    \csname \romannumeral\the\count@ pt\expandafter\endcsname
    \csname @\romannumeral\the\count@ pt\endcsname
  \csname #3\endcsname}%
\fi
\fi\endgroup
\begin{picture}(2305,2306)(234,-1711)
\put(1601,263){\makebox(0,0)[lb]{\smash{\SetFigFont{12}{14.4}{rm}{$\scriptstyle gT_{j+1}'$}%
}}}
\put(1235,-1537){\makebox(0,0)[lb]{\smash{\SetFigFont{12}{14.4}{rm}{$\scriptstyle S_i$}%
}}}
\put(812,-1228){\makebox(0,0)[lb]{\smash{\SetFigFont{12}{14.4}{rm}{$\scriptstyle T_{i+1}$}%
}}}
\put(353,-706){\makebox(0,0)[lb]{\smash{\SetFigFont{12}{14.4}{rm}{$\scriptstyle T_{i+2}$}%
}}}
\put(2084,-280){\makebox(0,0)[lb]{\smash{\SetFigFont{12}{14.4}{rm}{$\scriptstyle gT_{j+2}'$}%
}}}
\put(623,-421){\makebox(0,0)[lb]{\smash{\SetFigFont{12}{14.4}{rm}{$\scriptstyle\gamma_i$}%
}}}
\put(1094,-418){\makebox(0,0)[lb]{\smash{\SetFigFont{12}{14.4}{rm}{$\scriptstyle \tau\cdot\ell$}%
}}}
\end{picture}
\end{minipage}
\begin{minipage}{5.5cm}
\begin{picture}(0,0)%
\includegraphics{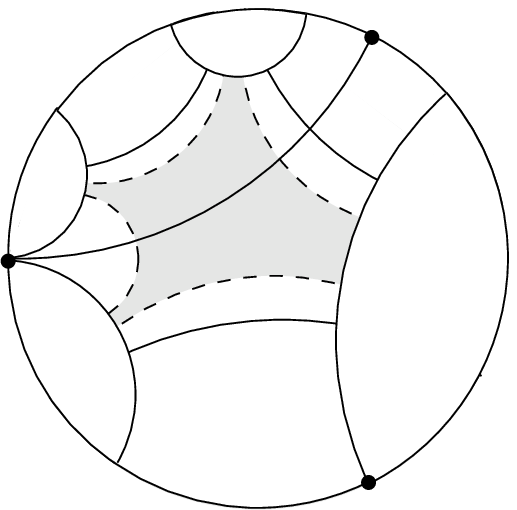}%
\end{picture}%
\setlength{\unitlength}{4144sp}%
\begingroup\makeatletter\ifx\SetFigFont\undefined
\def\x#1#2#3#4#5#6#7\relax{\def\x{#1#2#3#4#5#6}}%
\expandafter\x\fmtname xxxxxx\relax \def\y{splain}%
\ifx\x\y   
\gdef\SetFigFont#1#2#3{%
  \ifnum #1<17\tiny\else \ifnum #1<20\small\else
  \ifnum #1<24\normalsize\else \ifnum #1<29\large\else
  \ifnum #1<34\Large\else \ifnum #1<41\LARGE\else
     \huge\fi\fi\fi\fi\fi\fi
  \csname #3\endcsname}%
\else
\gdef\SetFigFont#1#2#3{\begingroup
  \count@#1\relax \ifnum 25<\count@\count@25\fi
  \def\x{\endgroup\@setsize\SetFigFont{#2pt}}%
  \expandafter\x
    \csname \romannumeral\the\count@ pt\expandafter\endcsname
    \csname @\romannumeral\the\count@ pt\endcsname
  \csname #3\endcsname}%
\fi
\fi\endgroup
\begin{picture}(2332,2300)(205,-1705)
\put(1853,221){\makebox(0,0)[lb]{\smash{\SetFigFont{12}{14.4}{rm}{$\scriptstyle T_3$}%
}}}
\put(1133,-479){\makebox(0,0)[lb]{\smash{\SetFigFont{12}{14.4}{rm}{$\scriptstyle F_q^A$}%
}}}
\put(510,-879){\makebox(0,0)[lb]{\smash{\SetFigFont{12}{14.4}{rm}{$\scriptstyle T_1$}%
}}}
\put(566, 31){\makebox(0,0)[lb]{\smash{\SetFigFont{12}{14.4}{rm}{$\scriptstyle \sigma_3T_1$}%
}}}
\put(1853,-591){\makebox(0,0)[lb]{\smash{\SetFigFont{12}{14.4}{rm}{$\scriptstyle T_2$}%
}}}
\put(784,284){\makebox(0,0)[lb]{\smash{\SetFigFont{12}{14.4}{rm}{$\scriptstyle\sigma_3T_2$}%
}}}
\end{picture}
\end{minipage}
\end{center}
\caption{Amalgamated free product/HNN-extension and the set $F_q^A$.}\label{amal}
\end{figure}
If $\ell(q_i)\ne 0$ and $\iota(q_i)=q_j'$ for some $q'\in\graph G_0$,
we use the combination theorems to form a new Fuchsian group
\cite[ch.~7]{maskit}. This group is either an amalgamated free product
of $\Delta_q$ and $\Delta_{q'}$, where $\gamma_i$ is identified with
$\gamma_j'^{-1}$ if $q\ne q'$, or an HNN-extension of $\Delta_q$ if
$q=q'$. An iteration of this procedure will yield the group
$\Gamma_q(\lambda)$. For the definition consider
\begin{equation*}
  X_{\graph G}:=\coprod_{q\in\graph G_0} \{q\}\times D.
\end{equation*}
We define an equivalence relation on $X_{\graph G}$ such that the
quotient is canonically isometric to $Y_{\graph
  G}(\lambda)$. The equivalence relation has the following generators:
\begin{enumerate}
\item $(q,z)\sim(q',z')$ if $q=q'$ and $z'=\delta z$ for some $\delta\in \Delta_q$.
\item Let $e\in\tilde{\graph G}_1$ with $\ell(e)\ne 0$. Let $q=s(e)$,
  so $e=q_i$ for some $i\in\IZ/3\IZ$ in the notation of section
  \ref{assembling}. The common orthogonal of $T_{i+1},T_{i+2}$ is
  denoted by $S_i$. Then $S_i$ is the axis of $\gamma_i\in\Delta_q$.
  We have $\iota(e)=q_j'$ for some $q'\in \graph G_0$, and if $S_j'$
  is the corresponding axis of $\gamma_j'\in\Delta_{q'}$, there exists
  a unique isometry $g$ of the hyperbolic plane such that (see figure \ref{amal})
  \begin{itemize}
  \item $g(S_j')=S_i$,
  \item $g^{-1}\gamma_j' g=\gamma_i^{-1}$,
  \item $\gamma_i^{\tau(d)}$ maps $T_{i+1}$ to $g(T_{j+1}')$.
  \end{itemize}
  Then we put $(q,gz)\sim (q',z)$.
\end{enumerate}
This gives $Y_{\graph G}(\lambda)\cong X_{\graph G}/\!\!\sim$ and
\begin{equation*}
  \Gamma_q(\lambda)=\{\gamma\in\isom(D)\st (q,z)\sim (q,\gamma z)\ \text{for all $z\in D$}\}.
\end{equation*}
Now the crucial property is that we may choose a fundamental polygon
$F_q$ in $\{q\}\times D$ for the action of $\Delta_q$ on its Nielsen
domain. Then the union of these polygons is a fundamental domain for
the equivalence relation on $X_\graph G$ up to a set of measure zero.
In particular, let $F_q$ be the set bounded by
$T_1,T_2,\sigma_3T_1,\sigma_3T_2$ and the common orthogonal of
$T_1,T_2$ and $\sigma_3T_1, \sigma_3T_2$ as on page \pageref{fdom}. If
$\inbase\subset \base$ is relatively compact, there exists $A>0$ such
that the cylinders $Z_{d}^A\subset Y_{\graph G}(\lambda)$ of area $2A$
exist for all $d\in\graph G_1^*$ and all $\lambda\in\inbase$. Then the
subset $F_q^A\subset F_q$ that is mapped onto the complement of
$\bigcup_d Z_{d}^A$ in $Y_{\graph G}(\lambda)$ has uniformly bounded
diameter for $\lambda\in\inbase$. (Note the cusp on the left of
$F_q^A$ in figure \ref{amal}. It belongs to $F_q$, so $F_q$ is of
infinite diameter.)  This observation allows to give a uniform
estimate on the number of closed geodesics in $Y_{\graph G}(\lambda)$.
\begin{lem}\label{bound geodesics} Let $N(\lambda;r)$ denote the number of unoriented, primitive, closed geo\-desics in $Y_{\graph G}(\lambda)$ of length at most $r$. Then there exist a neighbourhood $\inbase\subset\base$ for each $\lambda_0\in\base$ and $C>0$ such that
    $N(\lambda;r)\le C\cdot e^r$
  holds for all $\lambda\in\inbase$ and all $r>0$.
\end{lem}
\begin{proof}
  Fix $\inbase$ and $A>0$ as above. Let $\lambda\in\inbase$. If
  $c\subset Y_{\graph G}(\lambda)$ is a closed geodesic that is not
  the central geodesic in one of the cylinders $Z_{d}$, then its
  intersection with $Y_{\graph G}(\lambda)\setminus\bigcup_d Z_{d}^A$
  is nonempty. So there exists $q\in\graph G_0$ and a geodesic in
  $\{q\}\times D$ that is mapped onto $c$ and has nonempty
  intersection with the set $F_q^A$. Moreover, this geodesic is the
  axis of some hyperbolic element in $\Gamma_q(\lambda)$ with
  translation length equal to the length of $c$. So $N(\lambda;
  r)-(3p-3+2n)$ is bounded from above by the number of hyperbolic
  isometries $\gamma$ such that $\gamma\in\Gamma_q(\lambda)$ for some
  $q$ and $\dist{F_q^A}{\gamma F_q^A}\le r$ holds. $3p-3+2n$ is the
  number of edges. As explained above, there is an upper bound $D$ for
  the diameter of $F_q^A$, so
  \begin{equation*}
    N(\lambda;r)-(3p-3+2n)\le \#\{\gamma\st \gamma F_q^A\subset B_{r+D}(F_q^A)\,\text{for some $q$}\}.
  \end{equation*}
  If $v>0$ is a uniform lower bound for the area of $F_q^A$ for all
  $q$ and all $\lambda\in\inbase$, this implies
  \begin{equation*}
    N(\lambda;r)-(3p-3+2n)\le (2p-2+n)\cdot \max\{k\in\IN\st kv\le\vol\, B_{r+2D}\},
  \end{equation*}
  where $B_{r}$ denotes a hyperbolic ball of radius $r$, and $2p-2+n$ is
  the number of vertices. This shows
  \begin{align*}
    N(\lambda;r)&\le 3p-3+2n+ (2p-2+n)\cdot\max\{ k\st kv\le 4\pi\sinh^2(D+r/2)\}\\
    &\le 3p-3+2n+4\pi(2p-2+n)\frac{\sinh^2(D+r/2)}{v}\\
    &\le 3p-3+2n+\frac{\pi(2p-2+n)}{v} e^{2D}\cdot e^r.
  \end{align*}
\end{proof}
Now we prove proposition \ref{prop zeta}. Let $\lambda_0=\{\ell_0,\tau_0)\in\base$. Let $\inbase$ be a relatively compact neighbourhood of $\lambda_0$ such that $\ell(d)=0$ implies $\ell_0(d)=0$ for all $(\ell,\tau)\in\inbase$. There are two kinds of geodesics that contribute to the quotient \eqref{zetaquot} for $\lambda \in\inbase$:
\begin{enumerate}
\item Those that cross a cylinder $Z_d\subset Y_{\graph G}(\lambda)$
  with $\ell_0(d)=0$.  As $\lambda$ converges to $\lambda_0$, their
  length goes to infinity, so the contribution of each of these to the
  zeta function converges to $1$.
\item Those that do not cross such a cylinder. They are given by
  conjugacy classes of isometries in $\Gamma_q(\lambda)$ that converge
  to elements of $\Gamma_q(\lambda_0)$, and each hyperbolic isometry
  in $\Gamma_q(\lambda_0)$ arises as such a limit.
\end{enumerate}
So if the Zeta function is replaced with a finite product over all
geodesics with length less than a fixed constant, the resulting
quotient in \eqref{zetaquot} is continuous at $\lambda_0$. Lemma
\ref{bound geodesics} implies a uniform estimate in $\inbase$ for the
remainder.

\section{The resolvent for geometrically finite surfaces}\label{resolvent}In
this section, the resolvent of the Laplacian on functions defined on
$Y_{\graph G}(\lambda)$ is examined for varying $\lambda$. In
particular, we address its behaviour as $\lambda\in\base$ approaches
the boundary. The main assertion, theorem \ref{continuity resolvent},
states that the resolvent is continuous in $\lambda$ if the spectral
parameter does not belong to $[1/4,\infty)$.

The proof uses meromorphic Fredholm theory to compose the resolvent
from those of auxiliary surfaces, where the latter are rather easy to
describe. So the first subsection is concerned with the auxiliary
surfaces, and its results are combined in the second subsection to
examine general surfaces of finite geometry.

This technique of applying meromorphic Fredholm theory has often been
applied to prove the existence of a continuation of the resolvent
across the essential spectrum of the Laplacian. More precisely, the
map $s\mapsto (\lap-s(1-s))^{-1}$, defined for all $s\in\IC$ such that
$s(1-s)$ belongs to the resolvent set and $\re(s)>1/2$, has a
meromorphic continuation to $\IC$ as a family in, say,
$B(\L^2_{\mathrm c},\H^2_{\mathrm {loc}})$. We will make use of the
continuation in section \ref{eisenstein}. Yet the continuity results
of this section only address the \emph{physical domain} $\re(s)>1/2$,
and similar results fail to hold for the continued resolvent even in
the case of elementary surfaces. We refer to Guillop\'e's text
\cite{guillope1} for further details on the continuation, our
description of the resolvent closely follows his presentation.

\subsection{Auxiliary surfaces}\label{auxiliary surfaces}
We provide all statements on the resolvent operators for model
surfaces that are needed to deduce similar results for a geometrically
finite surface. The models are either compact or quotients of the
hyperbolic plane by an elementary group of isometries. In the compact
case, Hilbert-Schmidt and trace class properties of the resolvent are
immediate, and so is its continuity with respect to an arbitrary
family of Riemannian metrics. We will give the arguments below. In the
case of elementary quotients of the hyperbolic plane, continuous
dependency on the metric refers to the identification of a pair of
cusps with the complement of a closed geodesic in a hyperbolic
cylinder of arbitrary diameter, as it is induced by the identity map
of $\IR^2$ via definition $\ref{model}$. Note that this identification
preserves the Riemannian volume, so it gives an isometry of
$\L^2$-spaces. Continuity of the resolvent as a bounded map into
(local) Sobolev spaces is then obtained from inspection of its
integral kernel (proposition \ref{convergence operator elementary}).
The major part of this section is aimed at a trace class property for
a truncated resolvent and, in particular, uniform boundedness of its trace
norm for cylinders of small diameter (proposition
\ref{unicauchyprop}).

\subsubsection{Compact surfaces}\label{compact surfaces}  Let $(Y,g)$ be a
compact Riemannian manifold of
dimension 2. The spectrum of the Laplacian $\lap$ is purely discrete, and
there exists a complete set of orthonormal eigenfunctions for $\lap$ in
$\L^2(Y,\vol_g)$. Weyl's asymptotic law states that the spectral counting
function $N(\lambda)$, i.e.~the number of eigenvalues below $\lambda$
according to their multiplicities, is asymptotic to $\lambda\vol(Y)/4\pi$ as
$\lambda$ increases. This implies that the resolvent $(\lap+1)^{-1}$ is of
Hilbert-Schmidt class, and $(\lap+1)^{-2}$ is of trace class.

We use the notion of trace class mappings between different Hilbert spaces to
reformulate this observation. A bounded linear map $T\colon E\to F$ of separable
Hilbert spaces is of trace class if the supremum over all sums
\begin{equation*}
  \sum_{i=1}^\infty \abs{\scal{T e_i}{f_i}}
\end{equation*}
is finite, where $(e_i)_{i\in\IN}$ runs through the complete orthonormal
systems of $E$, and $(f_i)$ through those of $F$. The supremum is denoted by
$\norm{T}_1$ or $\norm{T}_{B_1(E,F)}$. The vector space $B_1(E,F)$ of trace class
operators is complete with respect to this norm. If $T\in B_1(E,E)$, the trace
of $T$ is defined by the series
\begin{equation*}
  \tr T:=\sum_i\scal{Te_i}{e_i},
\end{equation*}
which is independent of the chosen orthonormal system.

There is a related concept of Hilbert-Schmidt operators, we refer to
H\"or\-man\-der \cite[pp. 185-193]{hoermander3} for the particulars.

In the present situation, the Hilbert spaces we are interested in are Sobolev
spaces of functions on $Y$. As a vector subspace of $\L^2(Y,\vol_g)$, for any
$k\ge 0$, the Sobolev space $\H^k(Y,g)$ is the image of $(\lap+1)^{-k/2}$. A
Hilbert space structure is defined by the requirement that
\begin{equation*}
  (\lap+1)^{-k/2}\colon\ \L^2(Y,\vol_g)\longrightarrow \H^k(Y,g)
\end{equation*}
be isometric. Then $(\lap+1)^{-2}$, if considered as a bounded operator on
$\L^2(Y,\vol_g)$, is equal to the composition of bounded maps
\begin{equation*}
  \L^2(Y,\vol_g)\xrightarrow{(\lap+1)^{-2}}\H^4(Y,g)\overset{\iota}{\longrightarrow}
  \L^2(Y,\vol_g).
\end{equation*}
Since the first of these is isometric, trace class property for the square of
the resolvent on $\L^2(Y,\vol_g)$ means that the inclusion $\iota$ is of trace class and
\begin{equation*}
  \norm{\iota}_{B_1(\H^4,\L^2)}=\norm{(\lap+1)^{-2}}_{B_1(\L^2,\L^2)}.
\end{equation*}
The following lemma is an immediate consequence.
\begin{lem} Let $A(t)\colon\L^2(Y,\vol_g)\rightarrow \H^4(Y,g)$ be a continuous
  family of bounded linear mappings. Then $\iota\circ A(t)\in
  B_1(\L^2(Y,\vol_g),\L^2(Y,\vol_g))$ also depends con\-tinuous\-ly on $t$.
\end{lem}
\begin{proof}
  We have
  \begin{equation*}
    \norm{\iota\circ \bigl(A(t)-A(t_0)\bigr)}_{B_1(\L^2,\L^2)}
    \le \norm{\iota}_{B_1(\H^4,\L^2)}\cdot\norm{A(t)-A(t_0)}.
  \end{equation*}
\end{proof}
\begin{cor}\label{corollary compact}
  Let $\lap_h$ denote the Laplacian on functions that is associated with a
  Riemannian metric $h$ on $Y$. Then
  \begin{equation*}
    h\longmapsto\tr\left((\lap_h-\lambda)^{-1}-(\lap_h-\lambda_0)^{-1}\right)
  \end{equation*}
  is continuous on the open set of metrics where it is defined.
\end{cor}
\begin{proof}
  Let $g$ be a Riemannian metric on $Y$. A linear isometry $\delta_h\colon\L^2(Y,\vol_g)\to L^2(Y,\vol_h)$ is defined by multiplication of
  functions with the square root of an appropriate Radon-Nikodym
  derivative. The above trace is equal to the trace of the following operator
  on $\L^2(Y,\vol_g)$:
  \begin{multline*}
    \delta_h^{-1}\circ\left((\lap_h-\lambda)^{-1}-(\lap_h-\lambda_0)^{-1}\right)
    \circ\delta_h\\
    =(\lambda-\lambda_0)\cdot\delta_h^{-1}\circ
    (\lap_h-\lambda)^{-1}(\lap_h-\lambda_0)^{-1}\circ\delta_h.
  \end{multline*}
  The coefficients of the differential operator
  $(\lap_h-\lambda)(\lap_h-\lambda_0)$ depend continuously on $h$, so it
  defines a continuous family of bounded maps $\H^4(Y,g)\to
  \L^2(Y,\vol_g)$. Therefore, the inverse maps also depend continuously on $h$,
  and the previous lemma is applicable.
\end{proof}
\subsubsection{Elementary surfaces}
An elementary hyperbolic surface is a quotient $Z=\Gamma\backslash D$ of the hyperbolic plane, where the group $\Gamma$ of isometries has a cyclic subgroup of finite index. We only consider elementary groups that are cyclic, generated by a hyperbolic or a parabolic isometry of $X$.

Our aim is to compare the resolvent of the Laplacians for different surfaces. Define
\begin{gather*}
  X:=\{(x,a)\in\IR^2\st a\ne 0\},
  \intertext{and}
  Z:=\langle\gamma\rangle\backslash X,\quad\text{where}\quad \gamma\colon (x,a)\longmapsto (x+1,a).
\end{gather*}
Recall that, for each $\ell\ge 0$, the split plane $X$ may be identified with a subset of the model $X_\ell$ for the hyperbolic plane (definition \ref{model}). If we fix the Euclidean volume on $X$, this gives rise to an isometry of $\L^2$-spaces. So there is a corresponding family of self-adjoint Laplacians $\lap_\ell$ acting on $\L^2(X)$ and on $L^2(Z)$. In each case, the spectrum of $\lap_\ell$ is $[1/4,\infty)$. We may thus consider the resolvent operators as an analytic family of bounded operators
\begin{equation*}
  s\longmapsto \bigl(\lap_\ell-s(1-s)\bigr)^{-1},\quad\re(s)>1/2,
\end{equation*}
and the half-plane $\{s\in\IC\st \re(s)>1/2\}$ is a maximal domain where it is defined.

We shall prove the following statement.
\begin{prop}\label{elementary resolvent}
  \begin{enumerate}
  \item Let $\psi$ be a smooth function of compact support on $Z$. Then
    \begin{equation*}
      s\longmapsto \psi\bigl(\lap_\ell-s(1-s)\bigr)^{-1}
    \end{equation*}
    depends continuously on $\ell\in [0,\infty)$ in the sense of locally uniform convergence of analytic families of bounded linear maps $\L^2(Z)\rightarrow \H^2(Z)$. If\/ $\re(s_0)>1/2$, the same holds true for
    \begin{equation*}
      \psi\Bigl(\bigl(\lap_\ell-s(1-s)\bigr)^{-1}-\bigl(\lap_\ell-s_0(1-s_0)\bigr)^{-1} \Bigr)\colon\L^2(Z)\longrightarrow \H^4(Z).
    \end{equation*}
  \item Let $\psi$ be bounded and supported in a cylindrical subset $Z^S\subset Z$ of finite area, and let $\chi$ be of compact support in $Z$. Then
    \begin{equation*}
      \psi\Bigl(\bigl(\lap_\ell-s(1-s)\bigr)^{-1}-\bigl(\lap_\ell-s_0(1-s_0)\bigr)^{-1}\Bigr)\chi
    \end{equation*}
    is of trace class, and it depends continuously on $\ell\in[0,\infty)$ as an analytic family of trace class operators.
  \end{enumerate}
\end{prop}
\begin{rem} The integral kernel of the resolvent $(\lap_\ell-s(1-s))^{-1}$ can be continued meromorphically in $s$ to
  the complex plane. The poles, which are called \emph{resonances}, are
  \begin{equation*}
    \begin{split}
      \left\{1/2\right\}\quad &\text{for a cusp $Z_0^+$},\\
      \left\{-n+{2\pi i m}/{\ell}\st n\in\IN_0,\ 
        m\in\IZ\right\}\quad&\text{for a cylinder $Z_\ell,$ $\ell\ne
        0$}.
    \end{split}
  \end{equation*}
  This indicates that one cannot expect a result like proposition
  \ref{elementary resolvent} to hold for the continuation on
  the left of the critical axis $\left\{s\in\IC\st \re(s)=1/2\right\}$.
\end{rem}
The Sobolev spaces in part 1 of proposition \ref{elementary resolvent}
don't need to be specified because the support of $\psi$ is compact.
In the second part, the support of $\psi$ is not necessarily compact.
We note that the methods of section \ref{compact surfaces} apply to
deduce the second part from the first if the degenerate case $\ell=0$
is excluded. The proof for small values of $\ell$ is given in the
forthcoming subsection.
\subsubsection{Proof of proposition \ref{elementary resolvent}}\label{proof elementary}
\paragraph{A continuity statement.}
Let $K_s^\ell$ denote the Schwartz kernel of $\bigl(\lap_\ell-s(1-s)\bigr)^{-1}$, acting on $\L^2(Z)$. It is a smooth function on the complement of the diagonal in $Z\times Z$.

The following is a well-known property of $K_s^\ell$.
\begin{lem}\label{kernel} Let
  \begin{equation*}
    k_s(t):= \frac{4^{s-1}}\pi \int_0^1 (x(1-x))^{s-1}(4x+t)^{-s},\quad t>0,
  \end{equation*}
  and
  \begin{align*}
    \sigma_\ell\colon X\times X&\longrightarrow \IR_+\cup\{\infty\},\qquad
    (z_1,z_2)\longmapsto 2\bigl(\cosh(\dist{z_1}{z_2})-1\bigr),
  \end{align*}
  where $\mathrm d$ denotes hyperbolic distance if $X$ is identified
  with $X_\ell$. Then, for all $s$ with $\re(s)>1/2$, the kernel
  $K_s^\ell$ is given by the function on $X\times X$
  \begin{equation*}
    (z_1,z_2)\longmapsto \sum_{n\in\IZ} k_s\circ\sigma_\ell(z_1,\gamma^n z_2).
  \end{equation*}
\end{lem}
By this formulation we mean that the series is absolutely convergent
on $\{(z_1,z_2)\in X\times X\st z_1\ne\gamma^n z_2\text{ for all
  $n\in\IZ$ } \}$, it is $\gamma$-invariant in each argument, and the
induced function on $Z\times Z\setminus \{(z,z)\st z\in Z\}$ coincides
with the kernel $K_s^\ell$.
\begin{proof}
  The analogous statement holds for an arbitrary hyperbolic surface if
  $\re(s)>1/2$ holds and if the critical exponent of a uniformising
  group is less that $\re(s)$. We refer to Elstrodt \cite{elstrodt1}
  for a proof.  In the present situation, the critical exponent is
  $1/2$ if $\ell=0$, and $0$ otherwise.
\end{proof}
Now we use lemma \ref{kernel} to deduce the first statement of proposition \ref{elementary resolvent}.
\begin{prop}\label{convergence operator elementary} Let $\psi$ be a smooth
  function on $Z$ of compact support. Then $(\ell,s)\mapsto \psi\bigl(\lap_\ell
  -s(1-s)\bigr)^{-1}$ is a continuous mapping of\/ $[0,\infty)\times \{s\in\IC\st
  \re(s)>1/2\}$ into the normed vector space of bounded linear maps
  $\L^2(Z)\to \H^2(Z)$.
\end{prop}
\begin{proof}
  First of all, we observe that $k_s\circ\sigma_\ell$ converges to
  $k_s\circ\sigma_{\ell_0}$ as $\ell\to\ell_0$, uniformly on every
  compact subset of $X\times X$. This even holds on the diagonal in
  the sense that the difference $k_s\circ\sigma_\ell-
  k_s\circ\sigma_{\ell_0}$ has a continuous extension, which converges
  locally uniformly to 0. Lemma \ref{uniform convergence quotient}
  below shows that the same holds true for the kernel
  $K_s^\ell$ on the quotient.
  
  Now let $B$ denote a compact neighbourhood of $\supp\psi$. There exists a
  positive number $C$ such that
  \begin{equation}\label{c e 1}
    \norm{\psi f}_{\H^2}\le C\left( \norm{f}_{\L^2(B)}+
    \norm{\bigl(\lap_{\ell_0}-s(1-s)\bigr)f}_{\L^2(B)}\right) 
  \end{equation}
  for any smooth function $f$, locally uniform in $s$. Since the metric on $B$
  converges as $\ell\to\ell_0$, there exist $C_\ell$ with $C_\ell\to 0$ such that
  \begin{equation*}
    \norm{(\lap_{\ell_0}-\lap_\ell) f}_{\L^2(B)} \le C_\ell \left(
      \norm{\bigl(\lap_\ell-s(1-s)\bigr) f}_{\L^2(B)} +\norm{ f}_{\L^2(B)}\right).
  \end{equation*}
  Therefore, if $f$ is a smooth function of compact support,
  \begin{multline*}
    \lefteqn{\norm{\psi\,\bigl(\lap_\ell-s(1-s)\bigr)^{-1}
        f-\psi\,\bigl(\lap_{\ell_0}-s(1-s)\bigr)^{-1} f}_{\H^2}}\\
    \begin{aligned}
      &\le C\left( \norm{\left(
          \bigl(\lap_\ell-s(1-s)\bigr)^{-1}-\bigl(\lap_{\ell_0}-s(1-s)\bigr)^{-1}\right)
        f}_{\L^2(B)}\right.\\&\left.\phantom{C\Bigl(\quad} +\norm{(\lap_{\ell_0}-\lap_\ell)
        \bigl(\lap_\ell-s(1-s)\bigr)^{-1}f}_{\L^2(B)}\right)\\
      &\le
      C\,\norm{\left(\bigl(\lap_\ell-s(1-s)\bigr)^{-1}-
          \bigl(\lap_{\ell_0}-s(1-s)\bigr)^{-1}\right) f}_{\L^2(B)}\\
      &\quad + CC_\ell\left( \norm{f}_{\L^2(B)}+
        \norm{\bigl(\lap_\ell-s(1-s)\bigr)^{-1}f}_{\L^2(B)}\right).
    \end{aligned}
  \end{multline*}
  The last summand can be estimated in terms of the Euclidean distance between
  $s(1-s)$ and the spectrum of $\lap_\ell$. 
  
  To estimate the first term, we exhaust $Z$ by subsets $\supp\chi_n$
  for compactly supported, smooth functions $\chi_n$. The
  above-mentioned locally uniform convergence of the kernel $K_s^\ell$
  proves for each $n$
  \begin{equation*}
    \norm{\left( \bigl(\lap_\ell-s(1-s)\bigr)^{-1}-\bigl(\lap_{\ell_0}-s(1-s)\bigr)^{-1}\right)\chi_n}_{\L^2(B)}\longrightarrow
    0,\qquad \ell\rightarrow \ell_0.
  \end{equation*}
  So we need to find an upper bound for
  \begin{equation*}
    \norm{\psi\left(\bigl(\lap_\ell-s(1-s)\bigr)^{-1}- \bigl(\lap_{\ell_0}-s(1-s)\bigr)^{-1}
      \right)(1-\chi_n)}
  \end{equation*}
  that decreases with growing $n$. We may require that the distance of
  $B$ from the support of $1-\chi_n$ be greater than $n$ when measured
  with respect to each of the given metrics for $\ell$ near $\ell_0$. Now
  the spectral theorem admits to define the resolvent by
  \begin{equation*}
    \bigl(\lap_\ell-s(1-s)\bigr)^{-1}=\int_0^\infty \frac{\sin(
    t\sqrt{\smash[b]{\lap_\ell-1/4}})} 
    {\sqrt{\smash[b]{\lap_\ell-1/4}}}\, e^{-(s-1/2) t}\dx t.
  \end{equation*}
  The operator
  $(\lap_\ell-1/4)^{-1/2}\sin(t\sqrt{\smash[b]{\lap_\ell-1/4}})$
  solves the wave equation
  \begin{equation*}
    \left[ \frac{\partial^2}{\partial t^2}+(\lap_\ell-1/4)\right]
    \frac{\sin(t\sqrt{\smash[b]{\lap_\ell-1/4}})}
    {\sqrt{\smash[b]{\lap_\ell-1/4}}}=0,
  \end{equation*}
  so it has unit propagation speed. Again by the spectral theorem,
  \begin{multline*}
    \norm{\bigl(\lap_\ell-s(1-s)\bigr)^{-1}(1-\chi_n) f}_{\L^2(B)}\\
    \begin{aligned}
      &= \biggl\lVert{\int_n^\infty \frac{\sin(t
        \sqrt{\smash[b]{\lap_\ell-1/4}})}
        {\sqrt{\smash[b]{\lap_\ell-1/4}}}\, e^{-(s-1/2) t}\dx t\,
        (1-\chi_n)f}\biggr\rVert_{\L^2(B)}\\
      &\le \frac{n-1}{\re(s)-1/2}\, e^{-(\re(s)-1/2)
        n}\norm{(1-\chi_n)f}_{\L^2(X_0)},
    \end{aligned}
  \end{multline*}
  and this vanishes as $n\rightarrow\infty$.
\end{proof}
The proof of proposition \ref{convergence operator
  elementary} is completed with the following observation.
\begin{lem} \label{uniform convergence quotient} The kernel $K_s^\ell$, considered as a smooth function on the complement of the diagonal in $Z\times Z$, depends continuously on $(s,\ell)\in\{s\in\IC\st\re(s)>1/2\}\times [0,\infty)$.
\end{lem}
\begin{proof} It is easy to see that a finite sum $(z_1,z_2)\mapsto \sum_{\abs{n}\le N} k_s(\sig_\ell(z_1,\gamma^n z_2))$ depends continuously on $s$ and $\ell$. So we only need to give a bound for the derivatives of the remainder in the series. We will show below
  \begin{equation}\label{lala}
    \frac{\mathrm{d}^n}{\mathrm{d}t^n} k_s(t) = (4+t)^{-s-n} f_n(s,t)
  \end{equation}
  with a function $f_n$ that is bounded on any compact subset of $\{
  s\in\IC\st \re(s)>1/2\}\times (0,\infty]$. Therefore it is
  sufficient to replace $k_s(t)$ with $(4+t)^{-s}$ and to give a bound
  in the $\C^0$-topology. Now
  \begin{align*}
    \sig_\ell(x_1,a_1,x_2,a_2)& = 2\biggl[ \frac{\cosh(\ell(x_1-x_2))-1}{\ell^2}\,
    \sqrt{\smash[b]{\ell^2+a_1^2}}\, \sqrt{\smash[b]{\ell^2+a_2^2}}\\
    &\quad\phantom{2\biggl[} +\frac{1}{\ell^2} \left( \sqrt{
        \smash[b]{\ell^2+a_1^2}}\,\sqrt{\smash[b]{\ell^2+ a_2^2}}-a_1a_2 \right)-1
    \biggr]\qquad(\ell>0),
    \intertext{and if $(x_1,a_1)$, $(x_2,a_2)$ belong to the same component of $X$}
    \sig_0(x_1,a_1,x_2,a_2)&=\frac{(a_1-a_2)^2}{a_1a_2}+ a_1a_2(x_1-x_2)^2.
  \end{align*}  
  In particular, the inequality
  \begin{equation}\label{simple est}
    \sig_\ell(x_1,a_1,x_2,a_2)\ge \abs{a_1a_2}(x_1-x_2)^2
  \end{equation}
  holds for all $\ell\ge 0$, so
  \begin{equation*}
    \biggl|\sum_{\abs n>N} \bigl(4+\sig_\ell(x,a_1,x+n,a_2)\bigr)^{-s}\biggr| \le
    \sum_{\abs n>N} \bigl( 4+\abs{a_1a_2} n^2\bigr)^{-\re(s)}.
  \end{equation*}
  The last expression is of order $N^{1-2\re(s)}$ as $N\to\infty$, uniformly
  if $\abs{a_1a_2}$ is bounded from below by a positive number.

  We still need to verify \eqref{lala}. This follows from
  \begin{equation*}
    k_s(t)= \frac{4^{s-1}}{\pi}\sum_{j=0}^{\infty}
    \frac{4^j}{j!}\frac{\Gamma(s+j)^2}{ \Gamma(2s+j)}(4+t)^{-s-j},
  \end{equation*}
  which is proved by substituting
  \begin{equation*}
    (4x+t)^{-s}=(4+t)^{-s}\left( 1+\frac{4x-4}{4+t}\right)^{-s}
  \end{equation*}
  in the definition of $k_s$ and expanding the second factor into a
  power series in $\frac{4x-4}{4+t}$.
\end{proof}
The proof of proposition \ref{convergence operator elementary} can
easily be modified to prove the first statement in proposition
\ref{elementary resolvent} concerning the difference
\begin{equation*}
  \psi\Bigl(\bigl(\lap_\ell-s(1-s)\bigr)^{-1}-\bigl(\lap_\ell-s_0(1-s_0)\bigr)^{-1}\Bigr).
\end{equation*}
We shall not do this here. Rather, before we proceed with the second
part of proposition \ref{elementary resolvent}, we want to include an
observation that will not be needed until the Selberg trace formula is
applied in section \ref{zeta function}. The expression arises as a
remainder therein that has to be examined separately.
\begin{lem}\label{convergence remainder trace}
  For all $A>0$ the following map is continuous:
  \begin{equation*}
    (\ell,s)\longmapsto \int_A^\infty \sum_{\substack{n\in\IZ\\ n\ne 0}}
    k_s(\sig_\ell (0,a,n,a))\, \dx a.
  \end{equation*}
\end{lem}
\begin{proof} Consider the two estimates
  \begin{align*}
    \int_B^\infty \sum_{n\ne 0} \sig_\ell(0,a,n,a)^{-\re(s)}\dx a
    &\le \int_B^\infty a^{-2\re(s)}\dx a \sum_{n\ne 0} \abs{n}^{-2\re(s)} \\
    &=\frac{\zeta(2 \re(s))}{ \re(s)-1/2}\, B^{1-2\re(s)}
  \end{align*}
  and
  \begin{align*}
    \int_A^B \sum_{\abs n>N} \sig_{\ell}(0,a,n,a)^{-\re(s)}\dx a
    &\le \int_A^B a^{-2\re(s)}\dx a \sum_{\abs n>N} \abs n^{-2\re(s)} \\
    &\le \frac{ A^{1-2\re(s)}-B^{1-2\re(s)}}{(\re(s)-1/2)^2} N^{1-2\re(s)}.
  \end{align*}
  They show that it is sufficient to replace the series
  with a finite sum and the domain of integration with a finite
  interval. Then the statement follows from the continuity of
  $\ell\mapsto \sigma_\ell$.
\end{proof}
\paragraph{Trace class estimates.} We prove the second part of proposition \ref{elementary resolvent}. Essentially the proof is split into two parts: Trace class property is immediate if $\ell>0$, for in this case $\psi$ corresponds to a compactly supported function on $X_\ell$. The first point is to show that the operator under consideration is of trace class for $\ell=0$. Then proposition \ref{unicauchyprop} states that the trace norm is uniformly bounded in $\ell$ if the support of $\psi$ is small, which reduces the proof to the case where $\psi$ is compactly supported, i.e.~part 1 of proposition \ref{elementary resolvent} applies.

As the elementary cylinder $Z$ is a quotient of $X=\{(x,a)\in\IR^2\st
a\ne 0\}$, the support condition on $\chi$ implies that there is a gap
between the support and the deleted closed curve in $Z$. The
support of $\psi$ in contrast may include a neighbourhood of this
curve, so $\psi$ corresponds to a function of unbounded support on the
disjoint union $Z_0$ of two cusps. Thus our first issue is to prove
that in this situation
\begin{equation*}
  \psi\Bigl(\bigl(\lap_0-s(1-s)\bigr)^{-1}-\bigl(\lap_0-s_0(1-s_0)\bigr)^{-1}\Bigr)\chi
\end{equation*}
is of trace class. Let us explain how the propositions below are
applied to prove this fact. We may assume that $\psi$ and $\chi$ are
supported in the component $Z^+$. In this paragraph we will use $Z^A$
synonymously with $Z^{[-A,A]}$ for numbers $A>0$. The first resolvent
formula implies that the operator is the strong limit of
\begin{equation*}
   \bigl(s(1-s)-s_0(1-s_0)\bigr)\psi\Bigl(\bigl(\lap_0-s(1-s)\bigr)^{-1}\kappa_A\bigl(\lap_0-s_0(1-s_0)\bigr)^{-1}\Bigr)\chi
 \end{equation*}
 as $A\rightarrow 0$, where $\kappa_A$ denotes the characteristic
 function of $Z^+\setminus Z^A$. Lemma \ref{hscusp} below states that
 these operators are indeed of trace class if $A>0$.
Then proposition \ref{cauchyprop} implies that, if $A$ converges to 0,
the sequence of operators is a Cauchy sequence in the trace
class topology. The limit of this sequence coincides with the
strong limit, so the latter is of trace class.

The crucial observation in this argument will be a growth estimate for the
kernel that implies proposition \ref{cauchyprop} and convergence in
the trace class topology.  Essentially the same estimate will be used
again, in proposition \ref{unicauchyprop}, to examine the trace norm
of
\begin{equation*}
  \psi\Bigl(\bigl(\lap_\ell-s(1-s)\bigr)^{-1}\Bigr)-\bigl(\lap_\ell-s_0(1-s_0)\bigr)^{-1}\Bigr)\chi
\end{equation*}
locally uniformly in $\ell$. There it proves that the norm
\begin{equation*}
  \norm{\psi_\epsilon \Bigl(\bigl(\lap_\ell-s(1-s)\bigr)^{-1}-\bigl(
\lap_\ell-s_0(1-s_0)\bigr)^{-1}\Bigr)\chi}_1
\end{equation*}
becomes uniformly small as $\epsilon\to 0$, where $\psi_\epsilon$ is the characteristic
function of $Z^\epsilon$. Then we may apply the first part of proposition \ref{elementary resolvent} to conclude the proof of the second.

We proceed with the technical assertions that contribute to the
proof. They rely on Hilbert-Schmidt properties for $\lap_\ell$, which
are proved using the description of the resolvent kernel in lemma
\ref{kernel}. As it turns out, the logarithmic singularity of the
resolvent kernel is irrelevant. The decisive property is its
asymptotic
\begin{equation*}
  k_s(t)=\Ord\bigl(t^{-\re(s)}\bigr),\quad t\rightarrow\infty,
\end{equation*}
and therefore the following lemma will be useful. It is purely
technical but provides a convenient formula for our estimates.
\begin{lem}\label{formula}
  Let $\ell\ge 0$ and $r>1/2$. Consider the function $h_\ell$ on
  $X\times X$ defined by
  \begin{equation*}
    h_\ell(z_1,z_2)=\sum_{m\in\IZ}\bigl( 1+\sigma_\ell(z_1,\gamma^m z_2)\bigr)^{-r}.
  \end{equation*}
  Then
  \begin{multline}\label{complic}
    \int_0^1\int_0^1 \abs{h_\ell(x_1,a_1,x_2,a_2)}^2\dx{x_2}\dx{x_1}\\
    \le \sqrt\pi\, \frac{\Gamma(2r-1/2)}{\Gamma(2r)}\, g_\ell(a_1,a_2,2r)+\pi\, \frac{\Gamma(r-1/2)^2}{\Gamma(r)^2}\, g_\ell(a_1,a_2,r)^2
  \end{multline}
  holds, where $g_\ell$ is defined for all positive $\ell$ by
  \begin{align*}
    g_\ell(a_1,a_2,r)&= \biggl[ 1+{\frac 2{\ell^2}}\Bigl( \sqrt{(\ell^2+a_1a_2)^2+\ell^2(a_1-a_2)^2}-(\ell^2+a_1a_2)\Bigr)\biggr]^{1/2-r}\\
    &\phantom{-\Bigl[}\cdot \Bigl((\ell^2+a_1a_2)^2+\ell^2(a_1-a_2)^2\Bigr)^{-1/4},
    \intertext{and $g_0$ by}
    g_0(a_1,a_2,r) &= \begin{cases}
      \bigl[ a_1a_2+(a_1-a_2)^2\bigr]^{1/2-r}(a_1a_2)^{-1+r}&\text{if $a_1a_2>0$,}\\
      0&\text{otherwise}.
    \end{cases}
  \end{align*}
\end{lem}
\begin{proof}
  If the product $\abs{h_\ell(x_1,a_1,x_2,a_2)}^2$ of the two series
  in \eqref{complic} is expanded, the integral becomes
  \begin{multline}\label{expand}
    \sum_{m,n}\int_0^1\!\!\int_0^1\bigl( 1+\sigma_\ell(x_1,a_1,x_2+m+n,a_2)\bigr)^{-r} \bigl(1+\sigma_\ell(x_1,a_1,x_2+m,a_2)\bigr)^{-r}\dx{x_2}\dx{x_1}\\
    =\sum_{n\in\IZ}\int_{-\infty}^\infty \bigl( 1+\sigma_\ell(0,a_1,x+n,a_2)\bigr)^{-r} \bigl(1+\sigma_\ell(0,a_1,x,a_2)\bigr)^{-r}\dx x.
  \end{multline}
  The $n=0$ contribution is immediately estimated by substituting the definition of $\sigma_\ell$. This gives for $\ell>0$
  \begin{align*}
    \sig_\ell(0,a_1,x,a_2) &= 2\,\frac{\cosh(\ell x)-1}{\ell^2}
    \sqrt{(\ell^2+a_1a_2)^2+  \ell^2(a_1-a_2)^2} \\
    &\quad +\frac 2{\ell^2}\left( \sqrt
      {(\ell^2+a_1a_2)^2+\ell^2(a_1-a_2)^2}-\bigl(\ell^2+a_1a_2\bigr)\right)\\
    &\ge x^2 \sqrt{(\ell^2+a_1a_2)^2+\ell^2(a_1-a_2)^2}\\
    &\quad +\frac 2{\ell^2} \left(
      \sqrt{(\ell^2+a_1a_2)^2+\ell^2(a_1-a_2)^2}-\bigl(\ell^2+a_1a_2\bigr)\right),
  \end{align*}
  so the formula
  \begin{equation}\label{int}
    \int_{-\infty}^\infty (\alpha+x^2)^{-r}\dx x=\sqrt\pi\, \frac{\Gamma(r-1/2)}{\Gamma(r)}\,\alpha^{1/2-r}
  \end{equation}
  yields a contribution of
  \begin{equation*}
    \int_{-\infty}^\infty \bigl(1+\sig_\ell(0,a_1,x,a_2)\bigr)^{-2r}\dx x \le
    \sqrt\pi\,\frac{\Gamma(2r-1/2)}{\Gamma(2r)}\, g_\ell(a_1,a_2,2r),
  \end{equation*} 
  where the function $g_\ell$ is as proposed. If $\ell=0$, by definition of $\sigma_0$ we obtain from equation \eqref{int}
  \begin{equation*}
    \int_{\infty}^\infty \bigl(1+\sig_0(0,a_1,x,a_2)\bigr)^{-2r}=\sqrt\pi\,\frac{\Gamma(2r-1/2)}{\Gamma(2r)}\, g_0(a_1,a_2,2r).
  \end{equation*}
  The previous two formulae give the first summand in \eqref{complic}.
  With respect to the remaining part of the series in equation
  \eqref{expand}, we observe that the integral decreases as $\abs n$
  increases, and therefore
  \begin{multline*}
    \sum_{n\ne 0}\int_{-\infty}^\infty \bigl( 1+\sigma_\ell(0,a_1,x+n,a_2)\bigr)^{-r}\bigl( 1+\sigma_\ell(0,a_1,x,a_2)\bigr)^{-r}\dx x\\
    \begin{aligned}
      &\le \int_{-\infty}^\infty\int_{-\infty}^\infty \bigl( 1+\sigma_\ell (0,a_1,x+y,a_2)\bigr)^{-r} \bigl(1+\sigma_\ell(0,a_1,x,a_2)^{-r}\dx x\dx y\\
      &=\biggl(\int_{-\infty}^\infty \bigl( 1+\sigma_\ell(0,a_1,x,a_2)\bigr)^{-r}\dx x\biggr)^2.
    \end{aligned}
  \end{multline*}
  Another application of \eqref{int} completes the proof.
\end{proof}
We continue with the first trace class property.
\begin{lem}\label{hscusp}
  Let $s,s_0\in\{w\in\IC\st\re(w)>1/2\}$. If $\kappa$ is a bounded function supported off the deleted curve, then
  \begin{equation*}
    \psi\bigl(\lap_0-s(1-s)\bigr)^{-1}\kappa\bigl(\lap_0-s_0(1-s_0)\bigr)^{-1}\chi
  \end{equation*}
  is of trace class.
\end{lem}
\begin{proof}
  We show that $\psi\bigl(\lap_0-s(1-s)\bigr)^{-1}\kappa$ is a Hilbert-Schmidt operator, the same argument proves that this also holds true for $\kappa\bigl(\lap_0-s_0(1-s_0)\bigr)^{-1}\chi$. We may further assume that
  \begin{equation*}
    \supp\psi\subset Z^S\cap Z^+,\qquad\supp\kappa\subset Z^+\setminus Z^A
  \end{equation*}
  for some $A>0$, so the aim is to see
  \begin{equation*}
    \int_{Z^S\cap Z^+}\int_{Z^+\setminus Z^A} \abs{K_s^0(z_1,z_2)}^2\,\vol(z_1)\,\vol(z_2)<\infty.
  \end{equation*}
  According to lemma \ref{kernel}, the kernel is given by
  \begin{equation*}
    (z_1,z_2)\longmapsto \sum_{m\in\IZ} k_s\circ\sigma_0(\tilde z_1,\gamma^m\tilde z_2),
  \end{equation*}
  where the function $k_s$ satisfies
  \begin{equation*}
    k_s(t)=\Ord(t^{-\re(s)}),\ t\rightarrow \infty,\qquad\text{and}\qquad k_s(t)\sim -\frac 1{4\pi}\log t,\ t\rightarrow 0.
  \end{equation*}
  Integration of the logarithmic singularity over a two-dimensional
  domain causes no difficulties, therefore $k_s\circ\sigma_0$ may be
  replaced with $(1+\sigma_0)^{-\re(s)}$. By lemma \ref{formula} the proof
  is reduced to
  \begin{equation*}
    \int_0^S\int_A^\infty g_0(a_1,a_2,2r)\,\dx {a_2}\,\dx{a_1}<\infty
    \quad\text{and}\quad
    \int_0^S\int_A^\infty g_0(a_1,a_2,r)^2\dx{a_2}\,\dx{a_1}<\infty.
  \end{equation*}
  The explicit formula for $g_0$ shows that the integrands here are of order
  \begin{equation*}
    \Ord(a_1^{-2+2r}),\ a_1\rightarrow 0,\qquad \Ord(a_2^{-2r}),\ a_2\rightarrow\infty;
  \end{equation*}
  and therefore the integrals do exist.
\end{proof}
As explained at the beginning, the next statement
establishes the analogue of lemma \ref{hscusp} if the cut-off
function $\kappa$ is replaced with the constant function 1.
\begin{prop}\label{cauchyprop}
  For each $r\in\bigl(0,\re(s_0)-1/2\bigr)$ there exists $C>0$ such that
  \begin{multline*}
    \norm{\psi\bigl(\lap_0-s(1-s)\bigr)^{-1}(\kappa_B-\kappa_A)\bigl(\lap_0-s_0(1-s_0)\bigr)^{-1}\chi}_1
    \le C\cdot \int_B^A a^{-1+r}\dx a
  \end{multline*}
  holds for all $A>0$ sufficiently small and $A>B>0$. In particular, if\/
  $A-B$ converges to\/ $0$, then the norm converges to\/ $0$ locally
  uniformly in $A$. 
\end{prop}
\begin{proof}
  We begin with a consideration that is local in $A$ and $B$. The
  operators are of trace class by lemma \ref{hscusp}, and the trace norm
  can be estimated by a product of Hilbert-Schmidt norms
  \begin{equation}\label{normprod}
    \norm{\psi\bigl(\lap_0-s(1-s)\bigr)^{-1}(\kappa_B-\kappa_A)}_2\cdot\,\norm{(\kappa_B-\kappa_A) \bigl(\lap_0-s_0(1-s_0)\bigr)^{-1}\chi}_2.
  \end{equation}
  The first factor is the square root of
  \begin{equation*}
    \int_{Z^A\setminus Z^B}\int_{\supp\psi} \abs{\psi(z_1) K_s^0(z_1,z_2)}^2\vol(z_1)\vol(z_2).
  \end{equation*}
  Observe that, for each $B>0$, the function
  \begin{equation*}
    f\colon (0,\infty)\longrightarrow \IR,\quad A\longmapsto \int_{Z^A\setminus Z^B} \int_{\supp\psi} \abs{\psi(z_1) K_s^0(z_1,z_2)}^2 \vol(z_1)\vol(z_2)
  \end{equation*}
  is continuously differentiable. This implies
  \begin{equation*}
    \norm{\psi\bigl(\lap_0-s(1-s)\bigr)^{-1}(\kappa_B-\kappa_A)}_2^2\le (A-B)(f'(A)+\epsilon)
  \end{equation*}
  for arbitrarily small positive $\epsilon$ if $A$ belongs to a fixed
  compact subset of $(0,\infty)$, and if $B$ is sufficiently close to $A$.
  The same argument is applied to the second Hilbert-Schmidt norm
  in \eqref{normprod}, so
  \begin{multline*}
    \norm{\psi\bigl(\lap_0-s(1-s)\bigr)^{-1}(\kappa_B-\kappa_A)\bigl(\lap_0-s_0(1-s_0)\bigr)^{-1}\chi}_1\\
    \le (A-B)\sqrt{f'(A)+\epsilon}\,\sqrt{ g'(A)+\epsilon},
  \end{multline*}
  where
  \begin{equation*}
    g(A)=\int_{Z^A\setminus Z^B}\int_{\supp\chi} \abs{K_{s_0}^0(z_1,z_2)\chi(z_2)}^2\vol(z_1)\vol(z_2).
  \end{equation*}
  This formula leads to the desired `global' estimate as follows: For
  fixed numbers $A$ and $B$, we choose a subdivision
  $B=C_0<C_1<\cdots<C_k=A$, and if this subdivision is sufficiently
  small, we get
  \begin{multline*}
    \norm{\psi\bigl(\lap_0-s(1-s)\bigr)^{-1}(\kappa_B-\kappa_A)\bigl(\lap_0-s_0(1-s_0)\bigr)^{-1}\chi}_1\\
    \begin{aligned}
      &\le \sum_{j=0}^{k-1} \norm{\psi\bigl(\lap_0-s(1-s)\bigr)^{-1}(\kappa_{C_j}-\kappa_{C_{j+1}})\bigl(\lap_0-s_0(1-s_0)\bigr)^{-1}\chi}_1\\
      &\le \sum_{j=0}^{k-1} (C_{j+1}-C_j)\sqrt{f'(C_{j+1})+\epsilon}\,\sqrt{g'(C_{j+1})+\epsilon}.
    \end{aligned}
  \end{multline*}
  The last sum can be considered as a Riemann sum, and by refining the subdivision, we see
  \begin{equation}\label{crux}
    \norm{\psi\bigl(\lap_0-s(1-s)\bigr)^{-1}(\kappa_B-\kappa_A)\bigl(\lap_0-s_0(1-s_0)\bigr)^{-1}\chi}_1\\
    \le \int_A^B \sqrt{f'(a)}\sqrt{g'(a)}\,\dx a.
  \end{equation}
  We finally need to estimate the integrand in \eqref{crux}. The
  fundamental theorem of calculus gives
  \begin{equation*}
    f'(a)= \int_{\supp\psi} \abs{\psi(z) K_s^0(z,(x,a))}^2 \vol(z),
  \end{equation*}
  where $(x,a)$ denotes an arbitrary point on the boundary of $Z^a$.
  The kernel $K_s^0$ is given by an infinite summation of $k_s$. As $a\rightarrow 0$, the logarithmic singularity of $k_s$ can
  only contribute a term of order $\Ord(a^{-1})$ to this integral,
  because the Euclidean diameter of hyperbolic balls in $X_0$ grows
  linearly as the centre approaches the line $a=0$. The remaining
  contribution of $k_s$ is bounded by $(1+\sigma_0)^{-\re(s)}$, and
  lemma \ref{formula} yields
  \begin{align}
    f'(a)&\le\Ord\bigl(a^{-1}\bigr)+ C_1 a^{-1+2\re(s)}\int_0^S \bigl(a_1a+(a_1-a)^2\bigr)^{1/2-2\re(s)}{a_1}^{-1+2\re(s)}\dx a_1\nonumber\\
    &\phantom{\le} +C_2 a^{-2+2\re(s)}\int_0^S \Bigl(\bigl(a_1 a+(a_1-a)^2\bigr)^{1/2-\re(s)}a_1^{-1+\re(s)}\Bigr)^{\!2}\dx a_1\label{splitint}
  \end{align}
  The second summand is of order $a^{-1/2}$. In view of the third, we
  assume $a<\min(1,S)$ and split the integral into two parts to see
  \begin{multline*}
    \int_0^S\left(\left( a_1a+(a_1-a)^2\right)^{1/2-\re(s)}
      a_1^{-2+2\re(s)}\right)^2\dx
    a_1\\
    \begin{aligned}
      &\le \int_0^{a^2}(a_1-a)^{2-4\re(s)} a_1^{-2+2\re(s)}\dx a_1+
      a^{1-2\re(s)}\int_{a^2}^S a_1^{-1}\dx a_1\\
      &\le a^{2-4\re(s)}\abs{1-a}^{2-4\re(s)}\int_0^{a^2}
      a_1^{-2+2\re(s)}\dx a_1+ a^{1-2\re(s)}\log\left( Sa^{-2}\right)\\
      &=\Ord\bigl(a^{1-2\re(s)-2\epsilon}\bigr)
    \end{aligned}
  \end{multline*}
  for each $\epsilon>0$, so $f'(a)=\Ord\bigl(a^{-1-2\epsilon}\bigr)$.

  The estimate for $g'$ is analogous, but the compact support of
  $\chi$ implies that the corresponding integrands in \eqref{splitint}
  are actually bounded in this case, so
  $g'(a)=\Ord\bigl(a^{-2+2\re(s_0)}\bigr)$. This proves
  $\sqrt{f'(a)g'(a)}=\Ord\bigl(a^{-3/2+\re(s_0)-\epsilon}\bigr)$.
\end{proof}
To complete the second part of proposition \ref{elementary resolvent},
we need to examine the dependency of the operators on $\ell$.  In case
$\psi$ vanishes in a neighbourhood of the deleted curve, continuity in
$\ell$ as a family of bounded maps $\L^2(Z)\to \H^4(\supp\psi)$ was
proved in the first part of proposition \ref{elementary resolvent}.
This in turn implies continuity in the trace class topology under
this assumption on $\psi$, so it is sufficient to prove the following.
\begin{prop}\label{unicauchyprop} For each $\epsilon>0$ let $\psi_\epsilon$ denote the characteristic function of the cylindrical set $Z^\epsilon$ of area $2\epsilon$. Given $\delta>0$ and an arbitrary compact subset $R$ of $\{s\in\IC\st\re(s)>1/2\}$, there exists $\epsilon$ such that
  \begin{equation*}
    \norm{\psi_\epsilon\Bigl(\bigl(\lap_\ell-s(1-s)\bigr)^{-1}-\bigl(\lap_\ell-s_0(1-s_0)\bigr)^{-1}\Bigr)\chi}_1\le \delta
  \end{equation*}
  holds for sufficiently small $\ell$ and for all $s,s_0\in R$.
\end{prop}
\begin{proof}
  The proof is similar to the trace class estimates for $\ell=0$, i.e. the
  operator is approximated by variation of the cut-off parameter $A$ in
  \begin{equation*}
    \psi_\epsilon \bigl(\lap_\ell-s(1-s)\bigr)^{-1}\kappa_A\bigl(\lap_\ell-s_0(1-s_0)\bigr)^{-1}\chi.
  \end{equation*}
  The next two equations follow from the explicit formulae in lemma
  \ref{formula}, we will go into further detail at the end of
  this proof. Firstly the resolvent has the following property for
  fixed $A>0$:
  \begin{equation}\label{hsf1}
    \norm{\psi_\epsilon\bigl(\lap_\ell-s(1-s)\bigr)^{-1}\kappa_A}_2= \Ord\bigl(\epsilon^{\re(s)-1/2}\bigr),
  \end{equation}
  and the implicit constant may be chosen uniformly for small $\ell$.
  Secondly, if we assume that $\chi$ is supported in $Z^+$, we observe that
  \begin{equation}\label{hsf2}
    \norm{(1-\kappa_{-A})\bigl(\lap_\ell-s_0(1-s_0)\bigr)^{-1}\chi}_2 =\Ord\bigl(\ell^{\re(s_0)-1/2}\bigr),\quad\ell\rightarrow 0
  \end{equation}
  holds for fixed $A>0$.

  Taking these properties for granted, we need to examine the behaviour of
  \begin{equation*}
    \psi_\epsilon\bigl(\lap_\ell-s(1-s)\bigr)^{-1}\kappa\bigl(\lap_\ell-s_0(1-s_0)\bigr)^{-1}\chi
  \end{equation*}
  if $\kappa$ is supported in a small neighbourhood of the deleted curve. Here we proceed exactly as in the proof of proposition \ref{cauchyprop} to see
  \begin{equation*}
    \norm{\psi_\epsilon\bigl(\lap_\ell-s(1-s)\bigr)^{-1}(\kappa_B-\kappa_A\bigr)\bigl(\lap_\ell-s_0(1-s_0)\bigr)^{-1}\chi}_1\\
    \le C \int_B^A a^{-1+r}\dx a
  \end{equation*}
  for small $A$ with $B<A<0$, and $C$ is locally independent of $\ell$
  and $\epsilon$. Note that, although the function $g_0$ that appeared
  in the proof has to be replaced with $g_\ell$ a priori, monotonicity
  of $g_\ell(a_1,a_2,r)$ in $\ell$ if $a_1a_2>0$ and $\lim_{\ell\to 0}
  g_\ell=g_0$ imply that the original estimates remain valid. 
  
  Now given $\delta$ and $R$ as proposed, the last estimate together
  with a symmetry consideration imply that we can choose $A>0$ such
  that
  \begin{equation*}
    \norm{\psi_\epsilon\bigl(\lap_\ell-s(1-s)\bigr)^{-1}(\kappa_{-A}-\kappa_A) \bigl(\lap_\ell-s_0(1-s_0)\bigr)^{-1}\chi}_1<\delta/3
  \end{equation*}
  holds. By \eqref{hsf2} and \eqref{hsf1}
  \begin{equation*}
    \norm{\psi_\epsilon\bigl(\lap_\ell-s(1-s)\bigr)^{-1}(1-\kappa_{-A}) \bigl(\lap_\ell-s_0(1-s_0)\bigr)^{-1}\chi}_1<\delta/3
  \end{equation*}
  holds if $\ell$ is sufficiently small, and then a suitable choice
  of $\epsilon$ gives
  \begin{equation*}
    \norm{\psi_\epsilon\bigl(\lap_\ell-s(1-s)\bigr)^{-1}\kappa_A\bigl(\lap_\ell-s_0(1-s_0)\bigr)^{-1}\chi}_1<\delta/3.
  \end{equation*}
  We proved the proposition except for equations \eqref{hsf1} and
  \eqref{hsf2}. The first of these reduces to the Hilbert-Schmidt
  estimate in lemma \ref{hscusp} because of
  \begin{equation*}
    g_\ell(a_1,a_2,r)\le g_\ell(a_1,-a_2,r)\qquad\text{if}\quad a_1a_2\ge 0
  \end{equation*}
  and by monotonicity in $\ell$. Equation \eqref{hsf2} follows from
  \begin{align*}
    g_\ell(a_1,a_2,r)&\le \Bigl[ {\textstyle\frac 2{\ell^2}}\bigl( \abs{\ell^2+a_1a_2}-(\ell^2+a_1a_2)\bigr)\Bigr]^{1/2-r} \abs{\ell^2+a_1a_2}^{-1/2}\\
    &\le C\ell^{2r-1}\abs{a_1a_2}^{-r},
  \end{align*}
  where we assume that $\abs{a_1a_2}$ is bounded from below by a
  positive number.
\end{proof}
We close the subsection with an additional observation that will
allow surfaces of infinite area to be covered by theorem
\ref{continuity resolvent}. Here the support of $\psi$ may be
unbounded, but $\psi$ and $\chi$ are supported in different components
of $Z$.
\begin{prop}\label{elareainf} Let $\psi$ be bounded and supported in $Z^-\setminus Z^\epsilon$ for some $\epsilon>0$, and let $\chi$ be of compact support in $Z^+$. Then
  \begin{equation*}
    \lim_{\ell\to 0}\norm{\psi\Bigl(\bigl(\lap_\ell-s(1-s)\bigr)^{-1}-\bigl(\lap_\ell-s_0(1-s_0)\bigr)^{-1}\Bigr)\chi}_1 =0
  \end{equation*}
  holds with locally uniform convergence in $s$ and $s_0$.
\end{prop}
\begin{proof}
  For each $A<0$ 
  \begin{equation*}
    \lim_{\ell\to 0}\norm{\psi\bigl(\lap_\ell-s(1-s)\bigr)^{-1}\kappa_A\bigl(\lap_\ell-s_0(1-s_0)\bigr)^{-1}\chi}_1=0
  \end{equation*}
  holds by the proof of proposition \ref{unicauchyprop}. We let
  $A\rightarrow -\infty$ and estimate the remainder. More precisely, we
  show (cf. \eqref{crux})
  \begin{equation*}
    \int_{-\infty}^A \sqrt{f'(a)}\sqrt{g'(a)}\dx a=\Ord\bigl(\ell^{\re(s_0)-1/2}\bigr).
  \end{equation*}
  The restrictions $a<A<0$ and $\epsilon>0$ imply 
  \begin{equation*}
    \int_{-\infty}^{-\epsilon} g_\ell(a_1,a,2\re(s))\,\dx a_1=\Ord(1),\qquad a\rightarrow -\infty,
  \end{equation*}
  and
  \begin{equation*}
    \int_{-\infty}^{-\epsilon} g_\ell(a_1,a,\re(s))^2\,\dx a_1=\Ord\bigl( \abs a^{-1}\bigr),\qquad a\rightarrow -\infty,
  \end{equation*}
  so $\sqrt{f'(a)}$ is bounded. The asymptotic of $g'$ is governed by
  \begin{equation*}
    \int_\epsilon^\infty g_\ell(a,a_1,\re(s_0))^2\,\dx a_1\le C\,\ell^{4\re(s_0)-2}\, \abs{a}^{-2\re(s_0)}\int_{\epsilon}^\infty \abs{a_1}^{-2\re(s_0)} \dx a_1,
  \end{equation*}
  so that $\int_{-\infty}^A \sqrt{f'(a)}\sqrt{g'(a)}\,\dx a$ exists at least provided that $\re(s_0)>1$ holds.
\end{proof}
\subsection{Geometrically finite surfaces}\label{resolvent geometrically finite}
We take up the nomenclature introduced in section \ref{surface
  geometry}. Associated with an admissible graph $\graph G$ and its augmentation $\graph G^*$, there is a
parameter space for hyperbolic structures
\begin{equation*}
  \base=\left\{\lambda=(\ell,\tau)\st \ell\colon\graph G_1^*\rightarrow [0,\infty),\
    \tau\colon\graph G_1^*\rightarrow\IR\right\},
\end{equation*}
and a total space of hyperbolic surfaces
\begin{equation*}
  \total=\bigcup_{\lambda\in\base}\{\lambda\}\times Y_{\graph
  G}(\lambda).
\end{equation*}
We fix $\lambda_0=(\ell_0,\tau_0)\in\base$, and the surfaces
$Y_{\graph G}(\lambda)$ for $\lambda$ near $\lambda_0$ will be
compared to $Y=Y_{\graph G}(\lambda_0)$ via the trivialisation of
$\total$ from section \ref{surface geometry}. Recall the open cover
defined on page \pageref{surface}
\begin{equation*}
  Y=\bigcup_{q\in\graph G_0} P_q\cup\bigcup_{d\in\graph G_1^*} Z_d.
\end{equation*}
The cylindrical parts $Z_d$ of $Y$ can be embedded isometrically into
elementary cylinders $\langle\gamma\rangle\backslash X_{\ell_0(d)}$.
To remove the ambiguity in this statement, we note that $Y$ is defined
in such a way that each orientation $e\in\tilde{\graph G}_1^*$ of
$d\in\graph G_1^*$ gives rise to such an embedding.
If $d$ is a proper edge, replacing $e$ with $\iota(e)$ corresponds to
composition of that embedding with $(x,a)\mapsto(-x,-a)$, so there
is a canonical inclusion of $Z_d$ for each $d\in\graph G_1^*$ into the
cylinder
\begin{equation*}
  \bar Z_d:=\begin{cases}
    \left(\left\{e\right\}\times\langle\gamma\rangle\backslash
      X_{\ell_0(d)} \cup \left\{\iota(e)\right\}\times\langle\gamma\rangle
      \backslash X_{\ell_0(d)}\right)/\sim\,&\text{if $d\in\graph G_1$},\\
    \{d\}\times\langle\gamma\rangle\backslash X_{\ell_0(d)}&\text{if $d\in\graph G_1^*\setminus\graph G_1$},
  \end{cases}
\end{equation*}
where the equivalence relation is generated by
$(e,(x,a))\sim(\iota(e),(-x,-a))$. If $d$ is a proper edge, the image
of this inclusion is a cylindrical subset $\bar Z_d^{A(\ell_0(d))}$
for the interval
\begin{equation*}
  A(t)=
  \begin{cases}
    \bigl( -\tfrac{t}{2 \sinh(t/2)}, \tfrac{t}{2
    \sinh(t/2)}\bigr)& \text{if $t\ne 0$},\\ 
    \left(-1,1\right)& \text{if $t=0$.}
  \end{cases}
\end{equation*}
In case of a phantom edge the image is the set $\bar Z_d^{A(\ell(d))}\cup\bar
Z_d^+$.

We want to use the embeddings to pull back operators acting on $\L^2(\bar Z_d)$. For this purpose, we choose a collection of smooth functions
\begin{equation*}
  \phi_d\colon\ Y\longrightarrow [0,1]\quad\text{for $d\in\graph G_1^*$},\qquad
  \supp\phi_d\subset Z_d,
\end{equation*}
that satisfy for some $\epsilon>0$
\begin{equation*}
  \begin{aligned}
    \phi_d(z)=1\quad&\text{if $d\in\graph
      G_1$ and $z\in Z_d^{A(\ell_0(d))\cdot \epsilon}$,} \\
    \phi_d(z)=1\quad&\text{if $d\in\graph G_1^*\setminus\graph G_1$ and $z\in
       Z_d^+\cup Z_d^{A(\ell_0(d))\cdot \epsilon}$.}
  \end{aligned}
\end{equation*}
Let the smooth functions $\psi_d$, $d\in\graph G_1^*$, fulfil the same support
conditions and in addition $\psi_d\cdot\phi_d=\phi_d$.
If $\re(s)>1/2$, the resolvent operators
\begin{equation*}
  R_d(s):=(\lap_{\bar Z_d}-s(1-s))^{-1}
\end{equation*}
exist. Then $\psi_d R_d(s)\phi_d$ is considered to be 
acting on $\L^2(Y)$, where functions that are supported in $Z_d$ are
identified with functions on $\bar Z_d$ via the embeddings just described.

This construction for $Y=Y_{\graph G}(\lambda_0)$ is carried over to
$Y_{\graph G}(\lambda)$ for all $\lambda$ near $\lambda_0$. There
exists a neighbourhood $\inbase$ of $\lambda_0$ such that, with
respect to a trivialisation of $\total\to\base$ over $\inbase$ as
constructed in section $\ref{surface geometry}$, the cut-off functions
can be chosen simultaneously for all $\lambda\in\inbase$. The
trivialising map of $\total\vert_\inbase$ is denoted by
\begin{equation*}
  \Psi_\lambda\colon\ Y=Y_{\graph G}(\lambda_0)\longrightarrow Y_{\graph G}(\lambda).
\end{equation*}
It maps the complement of all closed geodesics in the
cylinders $Z_d$ diffeomorphically onto its image, and there is an
induced linear homeomorphism $\L^2(Y,\vol)\to \L^2(Y_{\graph
  G}(\lambda),\vol)$. The pull-back of an operator to $\L^2(Y,\vol)$ will be
labelled with a subscript $\lambda$.
\begin{thm}\label{continuity resolvent}
  \begin{enumerate}
  \item Let $\psi$ be a smooth function on $Y_{\graph G}(\lambda_0)$, compactly
    supported in the complement of the closed geodesics in the $Z_d$. If
    $\re(s)>1/2$, the bounded linear maps
    \begin{equation*}
      \psi\cdot(\lap_\lambda-s(1-s))^{-1}\colon\ \L^2(Y_{\graph
        G}(\lambda_0))\longrightarrow \H^2(Y_{\graph G}(\lambda_0))
    \end{equation*}
    depend continuously on $\lambda$ in the strong topology.
  \item Let $\re(s)>1/2$, $\re(s_0)>1$, and let $s(1-s)$ belong to the
    resolvent set. Then the operator on $\L^2(Y_{\graph G}(\lambda))$
    \begin{equation}\label{laplacian}
      T(s):= (\lap-s(1-s))^{-1}-(\lap-s_0(1-s_0))^{-1} -\sum_{d\in\graph
      G_1^*} \psi_d(R_d(s)-R_d(s_0))\phi_d
    \end{equation}
    is of trace class, and the trace depends continuously on $\lambda$.
  \end{enumerate}
\end{thm}
\begin{rem}\label{continuity remark 2}
  One concludes from part 1 of this theorem that the same result holds if the
  resolvent is considered as a map $\H^k_{\mathrm c}(Y')\to
  \H^{k+2}_{\mathrm{loc}}(Y')$, where $Y'$ is the
  complement of the distinguished closed geodesics. The proof is
  analogous to that of proposition \ref{convergence operator
    elementary}. 
\end{rem}
Appearance of the strong operator topology in theorem \ref{continuity
  resolvent} is a consequence of how the trivialisation $\Psi_\lambda$
was chosen. If $Z_d(\lambda)$ and $\bar Z_d(\lambda)$ are the
cylinders defined for $Y_{\graph G}(\lambda)$ in analogy with $Z_d$, $\bar Z_d\subset Y_{\graph G}(\lambda_0)$, there are canonical
diffeomorphisms
\begin{equation}\label{canonical elementary}
  \bar\Psi_{d,\lambda}\colon\ \bar Z_d\longrightarrow \bar Z_d(\lambda)
\end{equation}
induced by the identity of $\IR^2$. It does not coincide with the map induced by the trivialisation $\Psi_\lambda$. Instead, the composition
\begin{equation*}
  \bar Z_d^\epsilon\cong Z_d^\epsilon \xrightarrow{\Psi_\lambda}
  Z_d^\epsilon(\lambda)\cong \bar Z_d^\epsilon(\lambda)
\end{equation*}
is given by
\begin{equation*}
  (x,a)\longmapsto
  \begin{cases} (x+(\tau_0(d)-\tau(d))/2,a)&\text{if $a<0$},\\
    (x-(\tau_0(d)-\tau(d))/2,a)&\text{if $a>0$}.
  \end{cases}
\end{equation*}
If we denote by $\Psi_\lambda$ the obvious extension of this map
to $\bar Z_d\to \bar Z_d(\lambda)$, the automorphism of $\L^2(\bar
Z_d)$ induced by $\bar\Psi_{d,\lambda}^{-1}\circ \Psi_\lambda$ depends
continuously on $\lambda$ in the strong operator topology, but not
necessarily in the norm topology. It is the identity mapping if we
restrict $\lambda$ to a subset of $\base$ where all twist-parameters
are constant. The proof of \ref{continuity resolvent} will show that
under this restriction strong continuity may be replaced with norm
continuity.

\paragraph{Proof of the theorem, part 1.}
  The approach is analogous to that described by Guillop\'e \cite{guillope1},
  where the resolvent is continued analytically across the
  essential spectrum of $\lap$ by using meromorphic Fredholm theory. See the
  book of Reed and Simon \cite{reed-simon4} for a statement of the Fredholm
  theorem.
  
  Let $P$ denote the interior of 
  \begin{equation*}
    Y\setminus\biggl( \bigcup_{d\in\graph G_1^*} Z_d^{A(\ell(d))\cdot\epsilon/2}
    \cup\bigcup_{d\notin \graph G_1} Z_d^+\biggr).
  \end{equation*}
  The family $\phi_d$ of functions on $Y$ is supplemented with an
  additional function $\phi_0$ to form a partition of unity.
  The support of $\phi_0$ is contained in $P$, and there exists a
  smooth function $\psi_0$ of support in $P$ with
  $\psi_0\cdot\phi_0=\phi_0$.
  
  Just as the cylindrical subsets $Z_d$ are embedded into $\bar Z_d$,
  we embed $P$ isometrically into an auxiliary surface $\bar P$.  This
  surface is required to be compact.  Composition of this map with
  $\Psi_{\lambda}^{-1}$ gives an inclusion of the corresponding subset
  of $Y_{\graph G}(\lambda)$ into $\bar P$, and we may assume that
  there is a continuous family $g_\lambda$ of Riemannian metrics given
  on $\bar P$ such that this inclusion is isometric in each case.
  
  Suppressing its dependence on $\lambda$ in our notation, we consider
  the resolvent operators
  \begin{equation*}
    R_0(s):= (\lap_{\bar P}+t)^{-1}\quad\text{and}\quad R_d(s):=\bigl(\lap_{\bar Z_d}-s(1-s)\bigr)^{-1}.
  \end{equation*}
  The parameter $t$ will be chosen suitably below. In the formulae to
  come, the summation index $j$ will run through $\graph G_1^*$ and
  $\{0\}$. On $\L^2(Y_{\graph G}(\lambda))$ there is the equation of
  operators
  \begin{equation}\label{parametrix1}
    (\lap-s(1-s))\circ\sum_{j}\psi_j R_j(s)\phi_j = (1+Q(s)),
  \end{equation}
  with a compact operator $Q(s)$ given by
  \begin{equation*}
    \begin{split}
      Q(s)&= (-t-s(1-s))\psi_0(\lap_{\bar P}+t)^{-1}\phi_0 + \sum_j
      [\lap,\psi_j] R_j(s)\phi_j.
    \end{split}
  \end{equation*}
  If $1+Q(s)$ is invertible and the resolvent operators exist, by
  \eqref{parametrix1} we have
  \begin{equation*}
    \Bigl(\sum_{j} \psi_j R_j(s)\phi_j\Bigr)\circ
    (1+Q(s))^{-1} = (\lap-s(1-s))^{-1}.
  \end{equation*}
  So continuity of the resolvent in $\lambda$ is reduced to that of
  the operators on the left-hand side.  And invertibility of $1+Q(s)$
  indeed holds: If $t=-s(1-s)$, the norm of $Q(s)$ becomes arbitrarily
  small as the real part of $s$ increases.  Thus there exists $t$ such
  that $1+Q(s)$ is invertible for some $s$.  We fix such a $t$, and
  analytic Fredholm theory implies that the map $s\mapsto
  (1+Q(s))^{-1}$ is meromorphic on $\{s\in\IC\st\re(s)>1/2\}$.
  
  We want to see that the pull-back $Q(s)_\lambda$ to $\L^2(Y)$ depends
  continuously on $\lambda$ in the norm topology of $B(\L^2(Y))$, and
  so do the $\psi_jR_j(s)_\lambda\phi_j$ in the strong topology of
  $B(\L^2(Y),\H^2_{\mathrm{loc}}(Y))$.
  
  In both cases, this is obvious for the contribution of the compact
  part $\bar P$, as these operators are associated with the Laplacian
  for a continuous family of metrics on a compact surface. So let us
  consider a contribution $[\lap_\lambda,\psi_d]
  R_d(s)_\lambda\phi_d$ to $Q(s)_\lambda$ for some $d\in \graph
  G_1^*$. In terms of the canonical map $\bar{\Psi}_{d,\lambda}$
  between elementary cylinders, cf. \eqref{canonical elementary},
  this factorises into
  \begin{multline*}
    \Psi^*_\lambda\circ[\lap_{\bar Z_d(\lambda)},\psi_d]\circ (\lap_{\bar
      Z_d(\lambda)}-s(1-s))^{-1}\circ\Psi_\lambda^{-1*}\circ\phi_d\\
    \begin{aligned}
      &=[\lap_\lambda,\psi_d]\circ\Psi_\lambda^*\circ(\lap_{\bar
        Z_d(\lambda)}-s(1-s))^{-1}\circ \Psi_\lambda^{-1*}\circ\phi_d\\
      &=[\lap_\lambda,\psi_d]\circ\Psi_\lambda^*
      \bar\Psi_{d,\lambda}^{-1*} \circ \bar\Psi_{d,\lambda}^* 
      (\lap_{\bar Z_d(\lambda)}-s(1-s))^{-1} \bar\Psi_{d,\lambda}^{-1*}
      \circ \bar\Psi_{d,\lambda}^*\Psi_\lambda^{-1*}\circ\phi_d.
    \end{aligned}
  \end{multline*}
  In proposition \ref{convergence operator elementary} we proved that
  $\chi \bar\Psi_{d,\lambda}^* (\lap_{\bar Z_d(\lambda)}-s(1-s))^{-1}
  \bar\Psi_{d,\lambda}^{-1*}$ is a continuous family of compact linear
  maps $\L^2(\bar Z_d)\to \H^1(\bar Z_d)$ if $\chi$ is compactly
  supported off the closed geodesic. The support of the the
  differential operator on the left also satisfies this condition.
  This and the strong continuity of
  $\bar\Psi_{d,\lambda}^*\circ\Psi_\lambda^{-1*}$ proves norm
  continuity of the composition and that of $Q(s)_\lambda$.
  
  The remaining operators $R_j(s)_\lambda$ are treated in the same
  manner, but we must stick to strong continuity since they are not
  composed with a compact operator.

\paragraph{Proof of the theorem, part 2.} Using the first resolvent
  formula, we may proceed as in the first part to see
    \begin{multline}\label{parametrix2}
    (\lap-s_0(1-s_0))\circ(\lap-s(1-s))\circ\biggl(\sum_{j\ne 0}\psi_j
    (R_j(s)-R_j(s_0))\phi_j +\psi_0 R_0(s)\phi_0\biggr)\\[.5em]
    =(s(1-s)-s_0(1-s_0))\cdot(1+Q(s)),
  \end{multline}
  where $R_0(s)$ is defined by
  \begin{equation*}
    R_0(s):=\frac{s(1-s)-s_0(1-s_0)}{-t-s_0(1-s_0)} \left[(\lap_{\bar
        P}+t)^{-1}-(\lap_{\bar P}-s_0(1-s_0))^{-1}\right],
  \end{equation*}
  and, just for completeness, we have
  \begin{equation*}
    \begin{split}
      Q(s)&= (-t-s(1-s))\psi_0(\lap_{\bar P}+t)^{-1}\phi_0 +
      \frac{-t-s(1-s)}{s(1-s)-s_0(1-s_0)} [\lap,\psi_0] R_0(s)\phi_0\\[.5em]
      &\quad +[\lap,\psi_0] \left((\lap_{\bar P}+t)^{-1} +(\lap_{\bar
          P}-s_0(1-s_0))^{-1}\right)\phi_0 \\[1em]
      &\quad+\sum_{d\in\graph G_1^*} [\lap,\psi_d]\left(R_d(s)+R_d(s_0)\right)\phi_d\\
      &\quad+\frac 1{-t-s_0(1-s_0)} [\lap,[\lap,\psi_0]]
      R_0(s)\phi_0\\
      &\quad+\frac 1{s(1-s)-s_0(1-s_0)} \sum_{d\in\graph G_1^*}
      [\lap,[\lap,\psi_d]] (R_d(s)-R_d(s_0))\phi_d.
    \end{split}
  \end{equation*}
  Again $Q(s)$ is compact and $1+Q(s)$ is meromorphically invertible
  by Fredholm theory, so that
  \begin{multline*}
    \biggl(\sum_{d\in\graph G_1^*} \psi_d (R_d(s)-R_d(s_0))\phi_d+\psi_0
    R_0(s)\phi_0\biggr)\circ
    (1+Q(s))^{-1}\\
    \begin{aligned}
      &= (s_0(1-s_0)-s(1-s))\cdot(\lap-s(1-s))^{-1}\circ (\lap-s_0(1-s_0))^{-1}\\
      &= (\lap-s(1-s))^{-1}-(\lap-s_0(1-s_0))^{-1}.
    \end{aligned}
  \end{multline*}
  According to this equation, the operator $T(s)$ that we want to
  examine is given by
  \begin{multline}\label{truncation}
    (\lap-s(1-s))^{-1}-(\lap-s_0(1-s_0))^{-1}- \sum_{d\in\graph G_1^*}\psi_d
    (R_d(s)-R_d(s_0))\phi_d\\
    =\biggl(\psi_0 R_0(s)\phi_0- \sum_{d\in\graph
      G_1^*}\psi_d\,\bigl(R_d(s)-R_d(s_0)\bigr)\,\phi_d Q(s)\biggr)\circ (1+Q(s))^{-1}.
  \end{multline}
  We want to know that the first factor in \eqref{truncation} is of
  trace class, that its pull-back to $Y=Y_{\graph G}(\lambda_0)$
  depends continuously on $\lambda$ in the trace class topology, and
  that the pull-back of the second factor depends continuously on
  $\lambda$ in the norm topology. Then one can make use of
  Radon-Nikodym derivatives as in corollary \ref{corollary compact} to
  complete the proof.
  
  Continuity of $(Q(s))_\lambda$ in $\lambda$ follows exactly as in
  part 1 of the theorem. To prove trace class property, we observe
  that the image of $\phi_d Q(s)$ is supported in the intersection of
  $\supp\phi_d$ with $\supp\psi_0$. This is a compact subset of the
  complement of the pinching geodesics, and it was shown in section
  \ref{auxiliary surfaces} that this implies trace class property for
  the composition with $\psi_d (R_d(s)-R_d(s_0))$. Since $\bar P$ is
  compact, this property holds for $R_0(s)$ as well, which proves that
  $T(s)$ is of trace class.
  
  To examine continuity of the first factor in \eqref{truncation}, the
  operators are pulled back to $\L^2(Y)$ using $\Psi_\lambda$, and then
  we have in analogy with the elementary contribution in part 1
  \begin{multline*}
    \psi_d(R_d(s)-R_d(s_0))_\lambda \phi_d\,Q(s)_\lambda\\
    \begin{aligned}
      &=\psi_d \Psi_\lambda^*\left[ (\lap_{\bar
      Z_d(\lambda)}-s(1-s))^{-1}-(\lap_{\bar
      Z_d(\lambda)}-s_0(1-s_0))^{-1}\right]\Psi_{\lambda}^{-1*}\phi_d\,
      Q(s)_\lambda\\
      &= \bigl(\Psi_\lambda^*\bar\Psi_{d,\lambda}^{-1*}\bigr)\circ\psi_d\bar \Psi_{d,\lambda}^*
      \left[ (\lap_{\bar Z_d(\lambda)}-s(1-s))^{-1}-(\lap_{\bar
      Z_d(\lambda)}-s_0(1-s_0))^{-1}\right]
      \bar\Psi_{d,\lambda}^{-1*}\\
      &\qquad\circ\left(\bar\Psi_{d,\lambda}^*
      \Psi_\lambda^{-1*}\right)\phi_d\circ Q(s)_\lambda.
    \end{aligned}
  \end{multline*}
  Now
  \begin{equation*}
    \psi_d \bar\Psi_{d,\lambda}^*\left[(\lap_{\bar
    Z_d(\lambda)}-s(1-s))^{-1}-(\lap-s_0(1-s_0))^{-1}\right]
    \bar\Psi_{d,\lambda}^{-1*}\chi
  \end{equation*}
  is a continuous family of trace class operators on $\L^2(Y)$ by
  propositions \ref{convergence operator elementary} and
  \ref{elareainf}. Strong continuity of
  $\bar\Psi_{d,\lambda}^*\Psi_\lambda^{-1*}$ implies that the same
  holds true for $T(s)_\lambda$. 
\section{Applications}\label{applications}
\subsection{Eisenstein functions and the scattering matrix}\label{eisenstein}
We saw that the resolvent of the Laplacians for a degenerating family
converges to that of the limit surface. For surfaces of finite
area, the resolvent is intimately related with the theory of
Eisenstein series. It is therefore natural to ask for functions that
approximate the Eisenstein series during this process. In the first
subsection, we define functions that meet this criterion to some
extend, as well as an appropriate notion of approximate scattering
matrices. The Maass-Selberg relation and the functional equation
generalise to this setting, at least if the surface is of finite area.
In the second subsection, we therefore restrict to the finite area
case and examine how the introduced structure behaves during
degeneration.

\subsubsection{Definitions and fundamental properties}
This part is not concerned with families of degenerating surfaces, but
with a fixed surface $Y=Y_{\graph G}(\lambda)$. We make use of the
meromorphic continuation of the resolvent $(\lap-s(1-s))^{-1}$ in $s$ from $\re(s)>1/2$ to the complex plane \cite{guillope1}. We also use the
conventions from section \ref{surface geometry} that associate with
each oriented edge $j\in\tilde{\graph G}_1^*$ coordinates for an
embedded cylinder $Z_j\subset Y$.

The original definition of Eisenstein series for hyperbolic surfaces
of finite area is based on certain eigenfunctions of the Laplacian on
elementary cusps. We give a definition, in terms of the hypergeometric
function, that serves the same purpose for half-cylinders of arbitrary
circumference.
\begin{dfn}\label{eigenfunction} For each $\ell\ge 0$ and $s\in\IC\setminus(-1/2-\IN_0)$ put
  \begin{equation*}
    h(\ell,s)\colon \IR\longrightarrow\IC,\quad a\longmapsto
    \begin{cases} 
      \abs{a}^{-s}F(s/2,1/2+s/2;1/2+s;-\ell^2/a^2)&\text{if $a< 0$},\\
      0&\text{otherwise.}
    \end{cases}
  \end{equation*}
\end{dfn}
We identify $h(\ell,s)$ in the obvious manner with a
rotationally-symmetric function on the elementary cylinder
$Z_\ell=\langle\gamma\rangle\backslash X_\ell$. 
The Laplacian on $Z_\ell$ is
\begin{equation*}
  \lap_{Z_\ell}= -(\ell^2+a^2)^{-1}\partial_x^2
  -(\ell^2+a^2)\partial_a^2-2a\partial_a.
\end{equation*}
One checks, for example by expanding the hypergeometric function into a power
series, that the restriction of $h(\ell,s)$ to the open
half-cylinder $Z_\ell^-$ belongs to the kernel of
$\lap_{Z_\ell}-s(1-s)$. The choice of these functions is
motivated by the following properties:
\begin{itemize}
\item On $Z_\ell^-$ they converge smoothly to $h(0,s)$
  as $\ell\to 0$, and this is precisely the function used to define
  Eisenstein series associated with a cusp.
\item For generic $s$, the functions $h(\ell,s)$ and $h(\ell,1-s)$ are
  linearly independent.
\item The asymptotic of $h(\ell,s)$ as $a\to -\infty$ is that of $h(0,s)$ for all $\ell$.
\end{itemize}
 As an aside, note that
$h(\ell,s)$ depends quite trivially on $\ell$ in the upper half-plane
model: By means of the isometries introduced in section \ref{surface
  geometry}, we obtain the function
\begin{equation*}
  z\longmapsto \ell^{-s}\,\abs{\tan(\arg z)}^s F(s/2,1/2+s/2;1/2+s;-(\tan(\arg
  z))^2).
\end{equation*}
Being an eigenfunction of a linear ordinary differential equation,
there is of course a continuation of $h(\ell,s)|_{Z_\ell^-}$ to
$Z_\ell$ as an eigenfunction if $\ell>0$, and this continuation will be
used below. But we want to stress out that
$h(\ell,s)$ has a jump discontinuity at $a=0$ by definition.

The asymptotic of $h(\ell,s)$ as $a\to -\infty$ admit to define
approximate Eisenstein functions by summation over elements of a
uniformising group, just as in the classical definition of Eisenstein
series. We analysed the resolvent of the Laplacian in section
\ref{resolvent}, so it is preferable to use the resolvent instead. Fix
a cut-off function $\chi\colon Z_\ell\to [0,1]$ that is smooth,
rotationally symmetric, and supported in a small neighbourhood of
$Z_\ell^+$, such that the support of $1-\chi$ is contained in
$Z_\ell^{(-\infty,-\epsilon)}$ for some $\epsilon>0$.  Then
$[\lap_{Z_\ell},\chi]$ is a differential operator, compactly supported
in $Z_\ell^-$, and $[\lap_{Z_\ell},\chi]\,h(\ell,s)$ is a smooth
function of compact support in $Z_{\ell}^-$.  Recall from section
\ref{surface geometry} how an oriented edge $j\in \tilde{\graph
  G}_1^*$ of a graph gives rise to an identification of
$Z_{\ell(j)}^{A(\ell(j))}$ with a collar in $Y$. We require the
support of $\chi\, h(\ell(j),s)$ to belong to this subset, so that
$\chi\, h(\ell(j),s)$ and $[\lap_{Z_{\ell(j)}},\chi]\,h(\ell(j),s)$
both can be considered as functions on $Y$.
\begin{figure}
\begin{center}
\begin{picture}(0,0)%
\includegraphics{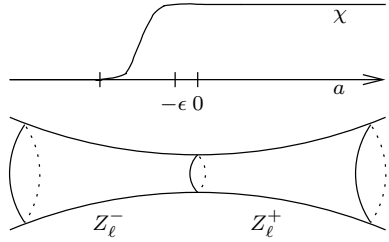}%
\end{picture}%
\setlength{\unitlength}{4144sp}%
\begingroup\makeatletter\ifx\SetFigFont\undefined
\def\x#1#2#3#4#5#6#7\relax{\def\x{#1#2#3#4#5#6}}%
\expandafter\x\fmtname xxxxxx\relax \def\y{splain}%
\ifx\x\y   
\gdef\SetFigFont#1#2#3{%
  \ifnum #1<17\tiny\else \ifnum #1<20\small\else
  \ifnum #1<24\normalsize\else \ifnum #1<29\large\else
  \ifnum #1<34\Large\else \ifnum #1<41\LARGE\else
     \huge\fi\fi\fi\fi\fi\fi
  \csname #3\endcsname}%
\else
\gdef\SetFigFont#1#2#3{\begingroup
  \count@#1\relax \ifnum 25<\count@\count@25\fi
  \def\x{\endgroup\@setsize\SetFigFont{#2pt}}%
  \expandafter\x
    \csname \romannumeral\the\count@ pt\expandafter\endcsname
    \csname @\romannumeral\the\count@ pt\endcsname
  \csname #3\endcsname}%
\fi
\fi\endgroup
\begin{picture}(2274,1474)(439,-1019)
\put(1891,-961){\makebox(0,0)[lb]{\smash{\SetFigFont{12}{14.4}{rm}{$\scriptstyle Z_{\ell}^+$}%
}}}
\put(946,-961){\makebox(0,0)[lb]{\smash{\SetFigFont{12}{14.4}{rm}{$\scriptstyle Z_{\ell}^-$}%
}}}
\put(1351,-241){\makebox(0,0)[lb]{\smash{\SetFigFont{12}{14.4}{rm}{$\scriptstyle -\epsilon$}%
}}}
\put(2386,-151){\makebox(0,0)[lb]{\smash{\SetFigFont{12}{14.4}{rm}{$\scriptstyle a$}%
}}}
\put(2386,299){\makebox(0,0)[lb]{\smash{\SetFigFont{12}{14.4}{rm}{$\scriptstyle\chi$}%
}}}
\put(1543,-244){\makebox(0,0)[lb]{\smash{\SetFigFont{12}{14.4}{rm}{$\scriptstyle 0$}%
}}}
\end{picture}
\end{center}
\caption{The cut-off function $\chi$.}
\end{figure}
\begin{dfn}\label{eisenstein series} If $\re(s)>1$, the \emph{approximate Eisenstein function} of $j\in\tilde{\graph G}_1^*$ is the function on $Y$ given by
  \begin{equation*}
    E_j(s):= (\lap-s(1-s))^{-1}\left([\lap_{Z_{\ell(j)}},\chi]
      h(\ell(j),s)\right)+\chi h(\ell(j),s).
  \end{equation*}
\end{dfn}
From the definition we immediately see that $E_j(s)$ specialises to the
Eisenstein series if $\ell(j)=0$. It is an eigenfunction of the
Laplacian on the subset of $Y$ where it is smooth, that is,
on the complement of a single closed geodesic.

The aforementioned continuation of the resolvent provides a
meromorphic continuation of $E_j(s)$ in $s$ to the complex plane. A
property that these functions inherit from the classical Eisenstein
series is that they satisfy some kind of Maass-Selberg relation, which
involves the notion of scattering matrices. We need a replacement for
the latter in the present situation.

For each function $f\in \L^1_{\mathrm{loc}}(Y)$ and each
$j\in\tilde{\graph G}_1^*$, the canonical coordinates of the
half-cylinder $Z_j^-$ give rise to fibrewise Fourier coefficients of $f$
by
\begin{equation*}
  F_j^nf(a):=\int_0^1 f(x,a) \exp(2 \pi inx)\dx x,\quad\textstyle{-\frac{\ell(j)}{2\sinh(\ell(j)/2)}<a<0}.
\end{equation*}
In particular, if $s-1/2\notin\IZ$ and if the restriction of $f$ to
$Z_j^-$ is annihilated by $\lap-s(1-s)$, then $F_j^0f$ is a linear
combination of $h(\ell(j),s)$ and $h(\ell(j),1-s)$. This observation
is used to define a pair of matrices that are associated with a
distinguished subset $S$ of the edges:
\begin{dfn}\label{scattering matrix}
  Let $S\subset \graph G_1^*$ be a set that contains all edges $j$
  with $\ell(j)=0$ and all phantom edges. Let $\tilde
  S\subset\tilde{\graph G}_1^*$ be the set of oriented representatives
  of edges in $S$. Then there are two analytic families
  $C=(C_{ij})_{i,j\in\tilde S}$ and $D=(D_{ij})_{i,j\in\tilde S}$ of
  matrices defined by the equation
  \begin{equation}\label{coeff} 
    F_j^0 E_i(s)= D_{ij}(s)\, h(\ell,s)+C_{ij}(s)\, h(\ell,1-s)
  \end{equation}
  if $\re(s)>1$ and $s-1/2\notin\IN$. The matrix $C(s)$ is called an
  \emph{approximate scattering matrix}.
\end{dfn}
Again we see that the approximate scattering matrix specialises to the
scattering matrix if $\lambda$ satisfies $\ell(d)=0$ for all $d\in S$.
It has a meromorphic continuation in $s$ to all of $\IC$.  The choice
of $S$ should be imagined as the choice of a certain stratum in the
boundary of the space $\base$ of hyperbolic surfaces constructed from
a given graph. Our objective in the next subsection will be to
approximate the Eisenstein series and the scattering matrix for a
surface $Y_{\graph G}(\lambda')$, where $\lambda'=(\ell',\tau')$
satisfies $\ell'(j)=0$ if and only if $j\in S$.

Some fundamental properties of the approximate scattering matrices are
deduced in lemma \ref{pair calculated} and \ref{D calc} below, before
the Maass-Selberg relation will be given in theorem
\ref{maass-selberg}. Then the functional equations for the approximate
Eisenstein functions and scattering matrices will be discussed as
consequences of the Maass-Selberg-relation.

Recall that $h(\ell,s)$ has simple poles in $s$ if $\ell\ne 0$. To
clarify the behaviour of $C$ at these points, we give an explicit
formula in terms of the constant modes $F_j^0 E_i(s)$.
\begin{lem} \label{pair calculated} Let $s$ be a regular point of the
  approximate Eisenstein function $E_i$. Then
  \begin{multline*}
    \begin{pmatrix} D_{ij}(s)\\ C_{ij}(s)
    \end{pmatrix}
    = \frac{\ell(j)^2+a^2}{1-2s}
    \begin{pmatrix} \partial_a h(\ell(j),1-s)& -h(\ell(j),1-s)\\
      -\partial_a h(\ell(j),s) &h(\ell(j),s)
    \end{pmatrix}
    \begin{pmatrix} F_j^0 E_i(s)\\ \partial_a F_j^0 E_i(s)
    \end{pmatrix}.
  \end{multline*}
\end{lem}
This formula indicates that there might be poles of the $D_{ij}$ in $\{s\in\IC\st \re(s)>1\}$, while the approximate scattering matrix in this
domain is holomorphic.
\begin{proof} As $h(\ell(j),s)$ and $h(\ell(j),1-s)$ are linearly independent solutions of a second-order ordinary differential equation, the coefficients $D_{ij}(s)$ and $C_{ij}(s)$ in \eqref{coeff} are determined by $F_j^0 E_i(s)$ and $\partial_a F_j^0 E_i(s)$ for an arbitrary point $a<0$. The formula then follows from solving the system
  \begin{equation*}
    \begin{pmatrix} F_j^0 E_i(s)\\ \partial_a F_j^0 E_i(s) \end{pmatrix}
    = \begin{pmatrix} h(\ell(j),s)& h(\ell(j),1-s)\\
      \partial_a h(\ell(j),s)&\partial_a h(\ell(j),1-s)\end{pmatrix}\,
    \begin{pmatrix} D_{ij}(s)\\ C_{ij}(s)\end{pmatrix}
  \end{equation*}
  by inversion of the Wronskian matrix. The determinant $\omega$ of the Wronskian solves
  \begin{equation*}
    \partial_a\bigl( (\ell^2+a^2)\,\omega(a)\bigr)=0.
  \end{equation*}
  It can be evaluated at $a=0$ by equation \eqref{eigen center} below, which gives $\omega(0)=(1-2s)\ell^{-2}$. So we obtain
  \begin{equation*}
    \omega(a)=\frac{\ell^2\,\omega(0)}{\ell^2+a^2}= \frac{1-2s}{\ell^2+a^2}.
  \end{equation*}
\end{proof}
If $C$ is the scattering matrix of a finite-area surface, one knows
that $D$ is the identity matrix. This does not hold in the more
general situation here, but the next lemma shows that $D$ can in fact
be calculated from $C$. The reason for this dependency is that the
approximate Eisenstein functions are smooth eigenfunctions up to a
single discontinuity that is known explicitly. So one simply has to
compare the functions $D_{ij}(s)\,h(\ell(j),s)+C_{ij}(s)\,
h(\ell(j),1-s)$, which are associated with $Z_j^-$, with those for
$Z_{\iota(j)}^-=Z_j^+$, to see if they match up correctly.
\begin{lem}\label{D calc}
  For all $s\in\IC$ let $\lambda^s$ and $\Sigma(s)$ denote the
  matrices with the following entries:
  \begin{align*}
    (\lambda^s)_{ij}&=\begin{cases} \ell(j)^s&\text{if $i=j$ and $\ell(j)\ne 0$,}\\
      0&\text{otherwise};\end{cases}\\
    \bigl(\Sigma(s)\bigr)_{ij}&= \begin{cases}
      \cos(\pi s)^{-1}&\text{if $i=j$,}\\
      1&\text{if $i\in\tilde{\graph G}_1$ and $j=\iota(i)$,}\\
      0&\text{otherwise.}\end{cases}
  \end{align*}
  Then $C$ and $D$ satisfy
  \begin{equation*}
    D(s)=1+ (2s-1)\frac{ 4^{-s}\Gamma(s)^2}{\Gamma(1/2+s)^2}\,
    C(s)\cdot\Sigma(s)\cdot\lambda^{2s-1}.
  \end{equation*}
\end{lem}
\begin{proof}
  We already mentioned that this formula is related with the
  continuation of $h(\ell,s)$ as a smooth eigenfunction from
  $Z_\ell^-$ to $Z_\ell$. If $\ell>0$ and $a<0$, we have
  \begin{multline}\label{eigen center}
    h(\ell,s)(a)= \ell^{-s} \frac{\Gamma(1/2)
      \Gamma(1/2+s)}{\Gamma(1/2+s/2)^2} F(s/2,1/2-s/2;1/2;-a^2/\ell^2)\\
    \qquad -\ell^{-(1+s)} \frac{\Gamma(-1/2) \Gamma(1/2+s)}{\Gamma(s/2)^2}a
    F(1/2+s/2,1-s/2;3/2;-a^2/\ell^2),
  \end{multline}
  and each summand on the right-hand side is an eigenfunction on the
  real line. The proof of this equation consists of a functional
  equation for the hypergeometric function, cf. \cite[p.~17,
  eq.~(17)]{erdelyi1}. Now if $a$ is positive and the same functional
  equation is applied to either summand in \eqref{eigen center} again,
  we see that it is equal to
  \begin{equation}\label{eigen right}
    \left[\cos(\pi s)^{-1} h(\ell,s) +\ell^{1-2s} \frac{4^s
        \Gamma(1/2+s)^2}{(2s-1) \Gamma(s)^2} h(\ell,1-s)\right](-a).
  \end{equation}   
  For completeness, let us carry out an intermediate step in the derivation of this formula: The functional equation gives
  \begin{multline*}
    \pi \Gamma(1/2+s)\Gamma(1/2-s) \left( \frac 1{\Gamma(1/2+s/2)^2
        \Gamma(1/2-s/2)^2} +\frac 1{\Gamma(s/2)^2\Gamma(1-s/2)^2} \right) \\[.5em]\cdot a^{-s}
    F(s/2,1/2+s/2;1/2+s,-\ell^2/a^2)\\[.5em]
    +2\pi \ell^{1-2s}
    \frac{\Gamma(1/2+s)\Gamma(-1/2+s)}{\Gamma(1/2+s/2)^2\Gamma(s/2)^2}
    a^{-(1-s)} F(1/2-s/2,1-s/2;3/2-s;-\ell^2/a^2).
  \end{multline*}
  Then \eqref{eigen right} follows from $\cos(\pi
  z)\Gamma(1/2+z)\Gamma(1/2-z)=\pi$ and from Legendre's duplication formula
  for the Gamma function.

  With this at hand, we determine $D_{ij}(s)$ if $j$ is a phantom edge and $\ell(j)\ne 0$. For each $i\in\tilde S$,
  \begin{equation*}
    E_i(s)-\delta_{ij}\,\chi\, h(\ell(i),s)
  \end{equation*}
  is an eigenfunction of $\lap$ in a neighbourhood of the closed
  geodesic in $Z_j$, and its constant term on $Z_j^-$ is
  \begin{equation*}
    (D_{ij}(s)-\delta_{ij}\,h(\ell,s)+C_{ij}(s)\,h(\ell,1-s).
  \end{equation*}
  According to \eqref{eigen right}, the constant term on $Z_j^+$ is 
  \begin{multline}\label{constant term right}
    \begin{aligned}&(D_{ij}(s)-\delta_{ij})\left[ \alpha(s)\, h(\ell,s)
      +\ell^{1-2s} \beta(s)\, h(\ell,1-s)\right](-a)\\
    &+ C_{ij}(s) \bigl[\alpha(1-s)\, h(\ell,1-s)
      +\ell^{1-2(1-s)} \beta(s)\, h(\ell,s)\bigr](-a)\end{aligned}\\
    \begin{aligned}&=\left[ (D_{ij}(s)-\delta_{ij})\alpha(s)+C_{ij}(s)\ell^{2s-1}\beta(1-s)\right] \cdot h(\ell,s)(-a)\\
    &\phantom{=}\ +\left[ (D_{ij}(s)-\delta_{ij}) \ell^{1-2s}\beta(s) +C_{ij}(s)\alpha(1-s)\right]\cdot h(\ell,1-s)(-a),\end{aligned}
  \end{multline}
  where $\alpha(s),\beta(s)$ are abbreviations for the coefficients in \eqref{eigen
    right}. If $\re(s)>1/2$, this function must be square-integrable
  on the infinite half-cylinder $Z_j^+$, so the coefficient of
  $h(\ell,1-s)$ must be zero.  This means
  \begin{equation*}
    D_{ij}(s)=\delta_{ij} +(2s-1) \frac{4^{-s} \Gamma(s)^2}{\Gamma(1/2+s)^2}
    C_{ij}(s) \cos(\pi s)^{-1} \ell^{2s-1}.
  \end{equation*}
  Now suppose that $j$ is a proper edge and $\ell(j)\ne 0$. The
  procedure to compute $D_{ij}(s)$ in this case is similar: The
  constant term of $E_i(s)-\delta_{ij}\chi h(\ell(i),s)$ on $Z_j^-$ is
  \begin{equation*}
    (D_{ij}(s)-\delta_{ij})\, h(\ell(j),s)+ C_{ij}(s)\, h(\ell(j),1-s),
  \end{equation*}
  and the constant term on $Z_{\iota(j)}^-\cong Z_j^+$ is again given
  by \eqref{constant term right}. On the other hand, the latter is equal to
  \begin{equation*}
    \bigl(D_{i\iota(j)}(s)-\delta_{i\iota(j)}\bigr)\,h(\ell(j),s)(-a) + C_{i\iota(j)}(s)\, h(\ell(j),1-s)(-a)
  \end{equation*}
  by definition. Equating the coefficients we see
  \begin{multline*}
    \begin{pmatrix} D_{ij}(s)-\delta_{ij} \\
      D_{i\iota(j)}-\delta_{i\iota(j)}
    \end{pmatrix}\\ 
    =(2s-1)\frac{4^{-s}\Gamma(s)^2}{\Gamma(1/2+s)^2} \ell(j)^{2s-1}
    \begin{pmatrix} \cos(\pi s)^{-1}&1\\ 1&\cos(\pi s)^{-1}
    \end{pmatrix}\cdot
    \begin{pmatrix} C_{ij}(s)\\ C_{i\iota(j)}(s)
    \end{pmatrix}.
  \end{multline*}
  We have determined the value of $D_{ij}(s)$ for all $j$ with
  $\ell(j)\ne 0$. If $\ell(j)=0$, then $D_{ij}(s)-\delta_{ij}=0$
  immediately follows from the square-integrability of $E_i(s)-\chi
  h(\ell(i),s)$ on a cusp if $\re(s)>1/2$.
\end{proof}
To state the Maass-Selberg relation, we need one more definition.
\begin{dfn}\label{truncted eisenstein}
  If $A$ is a sufficiently small positive number, the \emph{truncated
    Eisenstein functions} are
  \begin{equation*}
    E_i^A(s)\colon\ z\longmapsto 
    \begin{cases}
      E_i(s)(z)-(F_j^0 E_i(s))(a)&\text{if $z=(x,a)\in Z_j^-$, $j\in \tilde S$ and $\abs
        a<A$,}\\
      E_i(s)(z)-(F_j^0 E_i(s))(a)&\text{if $z=(x,a)\in Z_j^+$ and $j\notin
        \tilde{\graph G}_1$},\\
      E_i(s)(z)&\text{otherwise}.
    \end{cases}
  \end{equation*}
\end{dfn}
\begin{thm}\label{maass-selberg}
  If both $s$ and $s'$ are regular points of the approximate
  Eisenstein functions, and if the truncated functions are
  square-integrable, then the following relation holds:
  \begin{multline*}
    (s(1-s)-s'(1-s'))\cdot \scal{E_i^A(s)}{E_j^A(\bar s')}\\
    \begin{aligned}
      = \sum_{e\in\tilde S}\bigl(\ell(e)^2+A^2\bigr) \Bigl[ &D_{ie}(s)\, D_{je}(s')\cdot\omega_{\ell(e)}(s,s';-A) \\[-3mm]
      &+ D_{ie}(s)\,C_{je}(s')\cdot\omega_{\ell(e)}(s,1-s';-A)\\[1mm]
      &+C_{ie}(s)\,D_{je}(s')\cdot\omega_{\ell(e)}(1-s,s';-A) \\
      &+C_{ie}(s)\,C_{je}(s')\cdot\omega_{\ell(e)}(1-s,1-s';-A)\Bigr],
    \end{aligned}
  \end{multline*}
  where
  \begin{equation*}
    \omega_\ell(s_1,s_2;a)=
    \det\begin{pmatrix} h(\ell,s_1)(a)& h(\ell,s_2)(a)\\
      \partial_a h(\ell,s_1)(a)& \partial_a h(\ell,s_2)(a)\end{pmatrix}.
  \end{equation*}
\end{thm}
\begin{proof} We proceed exactly as Kubota \cite{kubota}. Let $f=E_i(s)$,
  $g=E_j(s')$, and let $V\subset Y$ be the subset where
  the approximate Eisenstein series are not truncated, i.e.
  \begin{equation*}
    V=Y\setminus\Bigl(\bigcup_{e\in\tilde
      S}Z_e^-\cup\bigcup_{ e\in\tilde{\graph G_1^*}\setminus \tilde{\graph
        G}_1} Z_e^{[0,\infty)}\Bigr).
  \end{equation*}
  The outer unit normal of $V$ at $Z_e^-$ is
  $\nu=\sqrt{\ell(e)^2+A^2}\,\partial_a$, and the one-dimensional volume form of $\partial V\cap Z_e^-$ is $\vol_{\partial V}=\sqrt{\ell(e)^2+A^2}\,\dx x$. Green's theorem states that
  \begin{equation}\label{MS1}
    -\int_V\bigl(f\lap g-g\lap f\bigr)\vol =\int_{\partial
      V}\bigl(f\cdot\nu(g)-g\cdot\nu(f)\bigr)\vol_{\partial V}.
  \end{equation}
  Insertion of the Fourier decomposition of $f$ and $g$ on the right leads to
  \begin{multline}\label{MS2}
      \!\sum_{e\in\tilde S}\int_{0}^1 \biggl[ \Bigl(\sum_{n\in\IZ} F_e^n f(-A)
      \exp(2\pi inx)\Bigr)\!\Bigl( \sum_{m\in\IZ}\nu (F_e^mg)(-A) \exp(2\pi imx)\Bigr)\\
      -\Bigl(\sum_{n\in\IZ} F_e^n g(-A) \exp(2\pi inx)\Bigr)\!\Bigl(
      \sum_{m\in\IZ} \nu (F_e^m f)(-A) \exp(2\pi
      imx)\Bigr)\biggr] \vol_{\partial}(x)\\
      =-\sum_{e\in\tilde S} \biggl[\bigl(\ell(e)^2+A^2\bigr)
      \cdot\sum_{n\in\IZ} \Bigl( (F_e^nf)'(-A)\,
      F_e^{-n}g(-A)-F_e^nf(-A)\,(F_e^{-n}g)'(-A)\Bigr) \biggl].
  \end{multline}
  For each $e\in\tilde{\graph G}_1^*\setminus\tilde{\graph G}_1$ and
  $n\ne 0$, we assume that the Fourier coefficients $F_e^nf$, $F_e^ng$
  are square-integrable eigenfunctions on the infinite half-cylinder
  $Z_e^{(-A,\infty)}$. The following shows that their respective
  contribution to the previous expression is given by an $\L^2$ scalar
  product over this cylinder:
  \begin{multline}\label{MS3}
    -\bigl(\ell(e)^2+A^2\bigr)\Bigl( (F_e^nf)'(-A)\, F_e^{-n}g(-A) -F_e^nf(-A)\,
    (F_e^{-n}g)' (-A)\Bigr)\\
    \begin{aligned}
      &=\int_{-A}^\infty \partial_a \Bigl[ \bigl(\ell(e)^2+a^2\bigr)\cdot \bigl( (F_e^nf)'(a)\,
      F_e^{-n}g(a) - F_e^nf(a)\, (F_e^{-n}g)'(a)\bigr)\Bigr]\dx a\\
      &= \int_{-A}^\infty \Bigl[ (\ell(e)^2+a^2) (F_e^nf)''(a)+2a
      (F_e^nf)'(a)\Bigr] F_e^{-n}g(a)\,\dx a\\
      &\phantom{=}-\int_{-A}^\infty F_e^nf(a) \Bigl[ \bigl(\ell(e)^2+a^2\bigr)
      (F_e^{-n}g)''(a)+2a (F_e^{-n}g)'(a)\Bigr]\dx a\\
      &= (-s(1-s)+s'(1-s'))\int_{-A}^\infty F_e^nf(a)\cdot F_e^{-n}g(a)\,\dx
      a.
    \end{aligned}
  \end{multline}
  This scalar product is part of the final formula. The same calculation gives for each proper edge $e\in\tilde{\graph G}_1$ and $n\ne 0$
  \begin{multline}\label{MS4}
    -\bigl(\ell(e)^2+A^2\bigr)\Bigl( (F_e^nf)'(-A)\, F_e^{-n}g(-A) -F_e^nf(-A)
    \,(F_e^{-n}g)'(-A)\Bigr)\\
    \begin{aligned}
      &= (-s(1-s)+s'(1-s'))\int_{-A}^0 F_e^nf(a)\cdot F_e^{-n}g(a)\,\dx
      a\\
      &\quad +\ell(e)^2 \Bigl( (F_e^nf)'(0)\, F_e^{-n}g(0)-
      F_e^nf(0)\,(F_e^{-n}g)'(0)\Bigr),
    \end{aligned}
  \end{multline}
  and we see that the additional summand here cancels with that of
  $\iota(e)$ in the sum over all $e\in\tilde S$ in \eqref{MS2}.
  Combining \eqref{MS1}, \eqref{MS2} and \eqref{MS3} with the
  assumption that $f=E_i^A(s)$ and $g=E_j^A(s')$ are eigenfunctions on
  the interior of $V$, we get
  \begin{multline*}
    \sum_{e\in\tilde S} \bigl(\ell(e)^2+A^2\bigr) \Bigl[ F_e^0f(-A)\, (F_e^0g)'(-A)
    -(F_e^0f)'(-A)\, F_e^0g(-A)\Bigr]\\
    \begin{aligned}
      &= (s(1-s)-s'(1-s')) \biggl[\int_S fg\vol
      +\sum_{e\in\tilde{\graph G}_1}\sum_{n\ne 0}
      \int_{-A}^0 F_e^nf(a)\,F_e^{-n}g(a)\dx a\\
      &\phantom{=} +\sum_{e\notin\tilde{\graph G}_1} \sum_{n\ne 0}
      \int_{-A}^\infty F_e^nf(a)\,F_e^{-n}g(a)\,\dx a\biggr]\\
      &=(s(1-s)-s'(1-s'))\cdot\scal{E_i^A(s)}{E_j^A(\bar s')}.
    \end{aligned}
  \end{multline*}
  The left-hand side is
  \begin{equation*}
    \sum_{e\in\tilde S}\bigl(\ell(e)^2+A^2\bigr)\,
    \det \begin{pmatrix} F_e^0f(-A)& F_e^0g(-A)\\ (F_e^0f)'(-A)& (F_e^0g)'(-A) \end{pmatrix}.
  \end{equation*}
  What remains is to substitute the definition of $C(s)$ and $D(s)$,
  and to expand the result using the linearity of the determinant
  in each column.
\end{proof}
Of particular interest is the Mass-Selberg relation in the cases $s=s'$ and
$s=1-s'$. The first yields a relation between $C$ and $D$ that is equivalent to the symmetry of $C$, the second case implies a functional equation for surfaces of finite area, which generalises the functional equation of the classical scattering matrix.
\begin{lem} \label{hyper}The function $\omega_\ell$ of theorem \ref{maass-selberg} satisfies 
  \begin{equation*}
    \bigl(\ell^2+a^2\bigr)\cdot\omega_\ell(s,1-s;a)=1-2s.
  \end{equation*}
\end{lem}
\begin{proof}
  Here $\omega_\ell$ specialises to the Wronskian determinant of $h(\ell,s)$ and $h(\ell,1-s)$. We already noted in the proof of lemma \ref{pair calculated} that
  \begin{equation*}
    a\longmapsto \bigl(\ell^2+a^2\bigr)\cdot\omega_\ell(s,1-s;a)
  \end{equation*}
  is constant on $(-\infty,0)$.
\end{proof}
\begin{cor} \label{c symmetric} The approximate scattering matrix is symmetric, and it satisfies
  \begin{equation*}
    C(s)\cdot D(s)^t= D(s)\cdot C(s).
  \end{equation*}
\end{cor}
\begin{proof} We begin with proving $C(s)\cdot D(s)^t=D(s)\cdot C(s)^t$, and this will imply the symmetry of $C(s)$. It is sufficient to
  consider $\re(s)>1$. There are no poles of $C$ or of the approximate
  Eisenstein functions in this domain, and the truncated Eisenstein
  series are square-integrable by definition. The Mass-Selberg
  relation for $s=s'$ gives via the preceding lemma
  \begin{align*} 0 &=(1-2s)\sum_{e\in\tilde S} (D_{ie}(s)\, C_{je}(s)-
  C_{ie}(s)\,D_{je}(s))\\
  &= (1-2s) \left( (D(s)\,C(s)^t)_{ij}-(C(s)\,D(s)^t)_{ij}\right),
  \end{align*}
  so we have $C(s)\cdot D(s)^t=D(s)\cdot C(s)^t$. Lemma \ref{D calc}
  implies that 
  \begin{equation*}
    C(s)(1+C(s)X)^t=(1+C(s)X)C(s)^t
  \end{equation*}
  holds for some symmetric matrix $X$, and in consequence
  \begin{equation*}
    C(s)+ C(s)XC(s)^t= C(s)^t+ C(s)XC(s)^t.
  \end{equation*}
\end{proof}
Apparently, the matrices $D(s)$ need not be symmetric in general.
Rather, lemma \ref{D calc} gives
\begin{align}\label{D symmetric}
  D(s)-D(s)^t &= 1+(2s-1) \frac{4^{-s}\Gamma(s)^2}{\Gamma(1/2+s)^2} C(s)\Sigma(s)\lambda^{2s-1}\notag\\
  &\quad -\biggl(1+(2s-1) \frac{4^{-s}\Gamma(s)^2}{\Gamma(1/2+s)^2}\lambda^{2s-1}\Sigma(s)C(s)\biggr)\notag\\
  &= (2s-1) \frac{4^{-s}\Gamma(s)^2}{\Gamma(1/2+s)^2}\ \left[ C(s),\Sigma(s)\lambda^{2s-1}\right].
\end{align}
In certain special cases, for example if all edges in the
distinguished set $S$ are phantom edges and the length $\ell(i)$ is
independent of $i\in S$, the commutator in equation \eqref{D
  symmetric} vanishes.

As a second consequence of the Maass-Selberg relation, we obtain a
functional equation for the approximate scattering matrices if the
surface is of finite area. In this situation, there is also a
functional equation for the approximate Eisenstein functions, corollary \ref{functional eisenstein}.

\begin{cor}\label{functional equation} The following functional equation holds if $Y$ is of finite area:
  \begin{equation*}
    D(s)\cdot D(1-s)^t= C(s)\cdot C(1-s).
  \end{equation*}
\end{cor}
\begin{proof}
  It suffices to prove the equation for all regular points $s\ne 1/2$,
  $\re(s)=1/2$ of the approximate Eisenstein functions. Let $(s_n)_n$ and
  $(t_n)_n$ be sequences in $\{w\in\IC\st \re(w)>1/2\}$ that converge to $s$
  and $1-s$, respectively. We may assume that the Eisenstein functions
  $E_i(s_n)$, $E_i(t_n)$ exist. The truncated ones $E_i^A(s_n)$ and
  $E_i^A(t_n)$ are square-integrable, and the matrices $C(s_n)$, $D(s_n)$,
  $C(t_n)$ and $D(t_n)$ converge. It is well-known in the classical theory
  of Eisenstein series that this implies the existence of a uniform bound for
  the $\L^2$-norms of $E_i^A(s_n)$ and $E_i^A(t_n)$ (cf.~Kubota \cite[theorem~4.1.2]{kubota}). Therefore
  \begin{equation*}
    \lim_{n\to\infty} (s_n(1-s_n)-t_n(1-t_n))\cdot \scal {E_i^A(s_n)}{E_j^A(\bar t_n)}=0.
  \end{equation*}
  The right-hand side of the Maass-Selberg relation converges to
  \begin{multline*}
    (1-2s)\, \sum_{e\in\tilde S} (D_{ie}(s)\ 
    D_{je}(1-s)-C_{ie}(s)\ C_{je}(1-s))\\
    =(1-2s) \bigl(D(s)\cdot D(1-s)^t -C(s)\cdot C(1-s)^t\bigr)_{ij}.
  \end{multline*}
\end{proof}
\begin{cor}\label{functional eisenstein}
  Assume that $Y$ is of finite area and that $D(s)$ is meromorphically invertible. Let $E(s)$ be the column vector that has one entry $E_j(s)$ for each $j\in\tilde S$. Then the functional equation
  \begin{equation*}
    E(s)=C(s)\,D(1-s)^{-1}\cdot E(1-s)
  \end{equation*}
  holds.
\end{cor}
  \begin{proof}
    We define a square matrix $Q(s)=(q_{ij}(s))$ from the constant
    terms of the Eisenstein series, namely by
    \begin{equation*}
      q_{ij}(s):=F_j^0 E_i(s).
    \end{equation*}
    If $h(s)$ denotes the diagonal matrix with entries $h(\ell(j),s)$, then
    \begin{equation*}
      Q(s)=D(s)\,h(s)+C(s),h(1-s).
    \end{equation*}
    The functional equation \ref{functional equation} and the relation
    in \ref{c symmetric} give
    \begin{align*}
      0&= D(s)\,h(s)+C(s)\,h(1-s)\\
      &\phantom{=\ }-C(s)\,D(1-s)^{-1}\cdot\bigl( D(1-s)\,h(1-s)+ C(1-s)\,h(s)\bigr)\\
      &= Q(s)- C(s)\,D(1-s)^{-1}\cdot Q(1-s).
    \end{align*}
    The latter is the constant term matrix of the column $E(s)-C(s)\,D(1-s)^{-1}\, E(1-s)$. This is zero for the following reasons:
    \begin{itemize}
    \item Each entry is a linear combination of Eisenstein functions
      such that the constant term is 0 on each cylinder. Thus it is a
      meromorphic family of smooth eigenfunctions of the Laplacian.
      But if $Y$ is compact, the spectrum is discrete.
    \item If $Y$ is non-compact, this follows from vanishing of the
      constant term on the cusps as in Kubota \cite[thm~4.4.2]{kubota}.
    \end{itemize}
\end{proof}
At the end of this subsection, we want to say a few words on the finite-area assumption in the functional equations.

Uniform boundedness of $\L^2$-norms, which allows to take the limit
$\re(s)\to 1/2$ of the Maass-Selberg relation, is essential in the
proof. In the general case of geometrically finite surface with
possibly infinite area, one knows from the definition of the
approximate Eisenstein functions in terms of the resolvent
\begin{equation*}
  \norm{E_i^A(s)}=\Ord\bigl((\re(s)-1/2)^{-1}\bigr),\ \re(s)\to 1/2.
\end{equation*}
The Maass-Selberg relation improves this bound:
\begin{cor}\label{eisenstein bound} Let $s\in\IC\setminus\{1/2\}$ with
  $\re(s)=1/2$ be a regular point of the approximate Eisenstein functions,
  and let $A$ be a positive number such that the truncated functions below are
  defined. If $(s_n)$ is a sequence in $\{w\in\IC\st \re(w)>1/2\}$ that converges to $s$, then the norms $\norm{E_i^A(s_n)}_{\L^2(Y)}$
  satisfy
  \begin{equation*}
    \norm{E_i^A(s_n)}_{\L^2(Y)} =\Ord\bigl((\re(s_n)-1/2)^{-1/2}\bigr),\ n\to\infty.
  \end{equation*}
\end{cor}
\begin{proof}
  The Maass-Selberg relation for $s=s_n$ and $s'=\bar s_n$ reads
  \begin{align*}
    (s_n-\bar s_n)& (1-2\re(s_n)) \norm{E_i^A(s)}^2\\
    &\begin{aligned}
      =\sum_{e\in\tilde S}\bigl(\ell(e)^2+A^2\bigr)& \Bigl[ \abs{D_{ie}(s_n)}^2 \,\omega_{\ell(e)}(s_n,\bar s_n;-A)\\[-.5em]
      &+ D_{ie}(s_n)\, C_{ie}(\bar s_n) \,\omega_{\ell(e)}(s_n,1-\bar s_n;-A)\\[.5em]
      &+C_{ie}(s_n)\, D_{ie}(\bar s_n)\,\omega_{\ell(e)}(1-s_n,\bar s_n;-A)\\[.2em]
      &+\abs{C_{ie}(s_n)}^2\,\omega_{\ell(e)}(1-s_n,1-\bar s_n;-A)\Bigr].
    \end{aligned}
  \end{align*}
  The assumptions now imply that the right-hand side of this formula converges
  to
  \begin{equation*}
    (1-2s) \sum_{e\in\tilde S} \bigl( \abs{D_{ie}(s)}^2-\abs{C_{ie}(s)}^2\bigr).
  \end{equation*}
\end{proof}
One might expect the $E_i^A(s_n)$ in corollary \ref{eisenstein bound} to
satisfy a stronger estimate like $\ord((\re(s_n)-1/2)^{-1/2})$, so that the previous proof would lead to the same functional equation as in corollary
\ref{functional equation}. But it turns out that the limit
\begin{equation*}
  \lim_{n\to\infty} (2\re(s_n)-1)\norm{E_i^A(s_n)}^2 =\sum_{e\in\tilde
  S} \bigl(\abs{D_{ie}(s)}^2-\abs{C_{ie}(s)}^2\bigr).
\end{equation*}
does not vanish.

Up to now, only the \emph{constant} term in the Fourier decomposition
of an eigenfunction was given explicitly, namely as a linear
combination of $h(\ell,s)$ and $h(\ell,1-s)$. The ordinary
differential equation that is solved by arbitrary Fourier modes $F_e^n f$ of
an eigenfunction $f$ is
\begin{equation*}
  \bigl(\ell^2+a^2\bigr) u''(a)+2au'(a) +\Bigl(s(1-s)-\frac{4\pi^2
  n^2}{\ell^2+a^2}\Bigr) u(a)=0.
\end{equation*}
It is solved (on the complement of $a=0$) by
\begin{equation*}
  \textstyle h^n(\ell,s)\colon\ a\longmapsto \abs{a}^{-s}(1+\ell^2/a^2)^{-i\frac{\pi n}\ell} F(s/2-i\frac{\pi n}\ell,
  1/2+s/2-i\frac{\pi n}\ell;1/2+s;-\ell^2/a^2).
\end{equation*}
If the surface $Y$ is of infinite area, and if $e$ is a phantom edge with
$\ell=\ell(e)\ne 0$, then $F_e^n(E_j(s))$ is a linear combination of
$h^n(\ell,s)$ and $h(\ell,1-s)$ on the infinite half-cylinder $Z_e^+$. It
follows from square-integrability that the coefficient of $h^n(\ell,1-s)$
vanishes if $\re(s)>1/2$ and $n\ne 0$, so
\begin{equation*}
  F_e^n(E_j(s))(a)= c(s)\, h^n(\ell,s)(a),\ a\in(0,\infty).
\end{equation*}
This shows
\begin{equation*}
  \int_{A}^{\infty} \abs{(F_e^n E_j(s))(a)}^2\dx a \sim\abs{c(s)}^2 \frac{
  A^{1-2\re(s)}}{2\re(s)-1},\ \re(s)\to 1/2.
\end{equation*}

Now suppose that $E_j^A(s)=\ord((\re(s)-1/2)^{-1/2})$ holds. The
asymptotic shows $\lim_{n}\abs{c(s_n)}^2=0$ for each sequence $(s_n)$
that converges from the right to a point $s$ with $\re(s)=1/2$. Thus
the coefficient $c(s)$ must vanish on the critical line. But this
means that $c$ is the zero function, and therefore all higher Fourier
modes $F_e^n E_j(s)$ vanish identically. This is not the case in
general, because the approximate Eisenstein functions do approximate
the Eisenstein functions at least if $\re(s)>1/2$.

\subsubsection{Approximating the scattering theory of finite area surfaces}
We proved in section \ref{resolvent} that the resolvent of the
Laplacian depends continuously on $\lambda\in\base$, and this implies
that the approximate Eisenstein functions $E_j(s)$ and scattering
matrices $C(s)$ depend continuously on $\lambda$ if $\re(s)>1/2$. In
particular, if $\lambda_0=(\ell_0,\tau_0)\in\partial B$ satisfies
$\ell_0(j)=0$ for all $j\in S$, and if $(\lambda_n)$ is a sequence in
$\base$ that converges to $\lambda_0$, then the approximate Eisenstein
functions and scattering matrices converge to Eisenstein series and
scattering matrix of the limit surface in the right half-plane. The
following theorem summarises what was proved so far.
\begin{thm}\label{approximation 1}
  Let $\lambda_0\in\base$ and $\re(s_0)>1/2$ such that
  $s_0(1-s_0)$ belongs to the resolvent set of the Laplacian on
  $Y_{\graph G}(\lambda_0)$.
  \begin{enumerate}
  \item There exists a neighbourhood $\inbase\subset\base$ of
    $\lambda_0$ and a neighbourhood $V$ of $s_0$ such that each $s\in
    V$ is not in the spectrum of $Y_{\graph G}(\lambda)$ for
    each $\lambda\in\inbase$. In particular, the approximate
    Eisenstein series $E_i(\lambda,s)$ and the matrices $C(\lambda,s)$
    are always defined, and so is $D(\lambda,s)$ if $s-1/2\notin\IN$.
  \item Let $\Phi^* E_i(\lambda,s)$ denote the pull-back of $E_i(\lambda,s)$
  to $Y_{\graph G}(\lambda_0)$ via the trivialisation maps from
  section \ref{surface geometry}. Then $\Phi^* E_i(\lambda,s)$ depends
  continuously on $\lambda\in\inbase$ and $s\in V$ in the sense of
  locally uniform convergence on $Y_{\graph G}(\lambda_0)$, the
  approximate scattering matrices $C(\lambda,s)$ and the
  $D(\lambda_k,s)$ depend continuously on $\lambda$.
  \end{enumerate}
\end{thm}
\begin{proof}
  Part 1 is only a reminder of theorem \ref{continuity
    resolvent}. Additional poles of $D(\lambda,s)$ at $s-1/2\in\IN$ might
  appear due to the poles of
  $s\mapsto h(\ell,1-s)$ (cf.~lemma \ref{pair calculated}).
  
  Theorem \ref{continuity resolvent} and remark \ref{continuity remark
    2} also show that $\Phi E_i(\lambda,s)$ depends smoothly on $s$,
  at least up to the twist caused by a variation of $\tau$. Now to
  determine the matrices $C(\lambda,s)$ and $D(\lambda,s)$, we only
  need the constant term of approximate Eisenstein functions and its
  first derivative in the complement of all closed geodesics that are
  associated with edges of the graph.
\end{proof}
Behaviour of the approximate scattering data on the left of
$\{s\in\IC\st\re(s)=1/2\}$ is completely different. We examine this
under the assumption that all surfaces are of finite
area, i.e.~we fix an admissible graph $\graph G$ of type $(p,0)$ for
the rest of this section. There are no phantom edges, and $Y_{\graph
  G}(\lambda)$ is compact for each $\lambda\in\interior\base$. If
$\lambda\in\partial\base$ then $Y_{\graph G}(\lambda)$ is of finite
area.  We also fix a set $S\subset\graph G_1$ as in definition
\ref{scattering matrix}, and therefore the parameter space for
surfaces must be restricted to
\begin{equation*}
  \base_S:=\{(\ell,\tau)\in\base\st \ell(d)\ne 0\text{ if $d\notin S$}\}.
\end{equation*}
In this situation, we can use the matrices $D(\lambda, s)$ in
conjunction with the approximate Eisenstein functions and scattering
matrices to achieve convergence on the complement of the critical
axis, see the corollary below. Its proof is based on the following
asymptotic.
\begin{prop}\label{D left}
  Let $\lambda_0=(\ell_0,\tau_0)\in\partial\base$ satisfy $\ell(d)=0$ for all
  $d\in S$.  Let $(\lambda_k)$ be a sequence in $\interior\base$ that
  converges to $\lambda_0$. If $\re(s)<1/2$, then
  \begin{equation*}
    C(\lambda_k,s)\sim \lambda_k^{1-2s} (1-2s) \frac{ 4^{-(1-s)}
      \Gamma(1-s)^2}{\Gamma(3/2-s)^2}\,\Sigma(1-s),\quad k\to\infty,
  \end{equation*}
  where the notation is that of lemma \ref{D calc}.
\end{prop}
\begin{proof}
  If $\lambda=(\ell,\tau)\in\interior \base$, the resolvent $s\mapsto
  (\lap_\lambda-s(1-s))^{-1}$ on $Y_{\graph G}(\lambda)$ is meromorphic on the
  complex plane for the spectrum is purely discrete. Thus
  the approximate Eisenstein functions are given by
  \begin{equation*}
    E_i(\lambda,s)=(\lap_\lambda-s(1-s))^{-1} \bigl([\lap_\lambda,\chi]
    h(\ell(i),s)\bigr) +\chi h(\ell(i),s)
  \end{equation*}
  even if $\re(s)<1/2$. Theorem \ref{continuity resolvent} implies that the
  pull-back of $E_i(\lambda_k,s)$ on $Y_{\graph G}(\lambda_0)$ converges to
  \begin{equation*}
    \tilde E_i(s):=(\lap_{\lambda_0}-s(1-s))^{-1} \bigl( [\lap_{\lambda_0},\chi]
    h(\ell_0(i),s)\bigr) +\chi h(\ell_0(i),s)\quad\text{if $\re(s)<1/2$.}
  \end{equation*}
  Note that $\tilde E_i(s)$ is defined in terms of the resolvent of the
  Laplacian, we do not need a meromorphic continuation across the
  critical axis. So $\tilde E_i(s) -\chi h(\ell_0(i),s)$ is
  square-integrable, and the constant term of this function on $Z_j^-$
  is
  \begin{equation*}
    F_j^0\tilde E_i(s)=\alpha_{ij}(s)\, h(\ell_0(j),s)+\beta_{ij}(s) \,h(\ell_0(j),1-s)
  \end{equation*}
  with $\beta_{ij}(s)=0$ if $\re(s)<1/2$. On the other hand we know that $\beta_{ij}(s)$ is the limit of
  $C_{ij}(\lambda_k,s)$, and this gives
  \begin{equation}\label{approx C 1}
    \lim_{k\to\infty} C(\lambda_k,s)=0,\ \re(s)<1/2.
  \end{equation}
  The same can be concluded for $D(\lambda_k,s)$ by the functional
  equation \ref{functional equation}: We know that $D(\lambda_0,s)$ is
  the identity matrix, and that it is the limit of $D(\lambda_k,s)$ if
  $\re(s)>1/2$. This implies that $s\mapsto D(\lambda_k,s)$ is
  meromorphically invertible if $k$ is sufficiently large, and the
  functional equation gives
  \begin{equation*}
    D(\lambda_k,s)= C(\lambda_k,s)\cdot C(\lambda_k,1-s)\cdot (D(\lambda_k,1-s)^t)^{-1}.
  \end{equation*}
  If $\re(s)<1/2$ the rightmost factor converges to the identity and
  $C(\lambda_k,1-s)$ converges to the scattering matrix, so $\lim_k
  D(\lambda_k,s)=0$ follows from \eqref{approx C 1}.

   Now the claim is proved by lemma \ref{D calc}:
  \begin{align*}
    C(\lambda_k,s)&= (D(\lambda_k,s)-1) \left( (2s-1) \frac{4^{-s}\Gamma(s)^2}
      {\Gamma(1/2+s)^2} \Sigma(s)\lambda^{2s-1}\right)^{-1}\\
    &=(D(\lambda_k,s)-1) \left( -\lambda^{1-2s}(1-2s)
      \frac{4^{-(1-s)}\Gamma(1-s)^2}{\Gamma(3/2-s)^2} \Sigma(1-s)\right).
  \end{align*}
\end{proof}
Proposition \ref{D left} suggests that the approximate scattering matrix should be replaced with the matrix defined now.
\begin{dfn}\label{scattering matrix 2} Let $\base_S'$ be the set of all $\lambda\in\base_S$ such that $s\mapsto D(\lambda,s)$ is meromorphically invertible. Then we define $C'\!(\lambda,s)$ for each $\lambda\in\base_S'$ to be the meromorphic family in $s$ given by
  \begin{equation*}
    C'\!(\lambda,s):=D(\lambda,s)^{-1}\cdot C(\lambda,s).
  \end{equation*}
\end{dfn}

Of course it is sufficient to require $D(\lambda,s)$ to be invertible for some $s$ in the definition of $\base_S'$. This always holds true if $Y_{\graph G}(\lambda)$ is compact, so we have $\interior\base\subset\base'_S$ (note that $D(\lambda,s)$ is the identity matrix for each $s\in -1/2-\IN$ by lemma \ref{D calc}).

Recall that $D(\lambda,s)=1$ if $\lambda(j)=0$ holds for each $j\in S$, so $C'\!(\lambda,s)$ again coincides with the scattering matrix. It turns out that there are additional properties that $C'\!(\lambda,s)$ shares with scattering matrices.
\begin{lem} The $C'\!(\lambda,s)$ are symmetric, they are unitary if $\re(s)=1/2$, and they satisfy the functional equation
  \begin{equation*}
    C'\!(\lambda,1-s)\cdot C'\!(\lambda,s)=1.
  \end{equation*}
\end{lem}
\begin{proof} We leave out $\lambda$ in the notation. Symmetry is a consequence of both statements in corollary \ref{c symmetric}:
  \begin{equation*}
    \bigl(C'\!(s)\bigr)^t=C(s)\cdot \bigl(D(s)^t\bigr)^{-1} = D(s)^{-1}\cdot C(s)=C'\!(s).
  \end{equation*}
  By definition of $C'$ this implies
  \begin{equation*}
    C'\!(1-s)\cdot C'\!(s)=C'\!(1-s)^t\cdot C'\!(s) =C(1-s)\cdot \bigl(D(s)\,D(1-s)^t\bigr)^{-1}\cdot C(s),
  \end{equation*}
  and this is the identity matrix by \ref{functional equation}. So if $\re(s)=1/2$
  \begin{equation*}
    C'\!(s)^*=C'\!(\bar s)^t=C'\!(\bar s)=C'\!(1-s)=C'\!(s)^{-1}.
  \end{equation*}
\end{proof}
We conclude this section with a corollary of theorem \ref{approximation 1}, of the functional equation for $C'$, and of corollary \ref{functional eisenstein}.
\begin{cor}\label{approximation 2}
  Let $E(\lambda,s)$ be the column vector of approximate Eisenstein functions as in corollary \ref{functional eisenstein}, and let $E'(\lambda,s):=D(\lambda,s)^{-1}\cdot E(\lambda,s)$. Then
  \begin{equation*}
    s\mapsto \Phi^*E'(\lambda,s)\qquad\text{and}\qquad s\mapsto C'\!(\lambda,s)
  \end{equation*}
  depend continuously on $\lambda\in\base_S'$ as meromorphic families on $\{s\in\IC\st\re(s)\ne 1/2\}$.
\end{cor}

\subsection{The Selberg Zeta function}\label{zeta function}
We extend proposition \ref{prop zeta} to the domain $\re(s)>1/2$, thus
proving Wolpert's second conjecture. As in section \ref{zetaconv},
$\Zeta_\lambda$ is the Selberg Zeta function of $Y_{\graph
  G}(\lambda)$, and the contribution of a particular edge $d\in\graph
G_1^*$ to the Zeta function is denoted by $\Zeta_{d,\lambda}$.
\begin{thm}\label{main zeta}  With respect to the topology of locally uniform
  convergence of functions on $\{s\in\IC\st \re(s)> 1/2\}$, the map
  $\lambda\mapsto\Zeta_\lambda/\prod_{d\in\graph G_1^*}\Zeta_{d,\lambda}$ is
  continuous on $\base$. 
\end{thm}
Before we come to the proof, let us state an immediate corollary of
this theorem and of the functional equation for $\Zeta_\lambda$ for
surfaces of finite area. It can be used to show that there is no
straight-forward extension of theorem \ref{main zeta} to the left of
the critical axis \cite{wolpert}. At the end of this section, we will
use this corollary to prove theorem \ref{zeta everywhere}, which
covers \ref{main zeta} but adds convergence on the left half-plane if
the non-degenerate surfaces are compact.
\begin{cor}\label{zeta left} Let $\lambda_0=(\ell_0,\tau_0)\in\partial B$ be a point where $Y(\lambda_0)$ is the union of elementary cusps and a surface of finite
  area. The area of the latter surface is denoted by $\abs F$. Let
  $(\lambda_m)$ be a sequence in $\base$ that converges to
  $\lambda_0$. Then
  \begin{align*}
    \lim_{m\to\infty}\frac{\Zeta_{\lambda_m}(s)}{\textstyle{\prod_{d\in S}}
    \Zeta_{d,\lambda_m}(s)}&=\pm\Zeta_{\lambda_0}(1-s)\,\det( C(\lambda_0,1-s))\cdot\biggl(\frac{\Gamma(1/2+s)}{\Gamma(3/2-s)}\biggr)^{\!\!k}\\[2mm]
  &\phantom{=}\cdot \exp
    \biggl( \abs F \int _0^{s-1/2} t \tanh(\pi t)\dx t +(2s-1) k\log 2\biggr)
  \end{align*}
  holds in $\{s\in\IC\st \re(s)>1/2\}$. Here $C(\lambda_0,s)$ is the
  scattering matrix of $Y(\lambda_0)$, and $k$ is the number of cusps.
  In particular, if the $Y(\lambda_n)$ are compact for all $n\ne 0$,
  we have for $\re(s)>1/2$
  \begin{equation*}
   \lim_{m\to\infty}\frac{\Zeta_{\lambda_m}(1-s)}{\textstyle{\prod_{d\in S}}\Zeta_{d,\lambda_m}(s)}= \Zeta_{\lambda_0}(1-s)\, \det( C(\lambda_0,1-s))\cdot 2^{(2s-1)k}\biggl(\frac{\Gamma(1/2+s)}{\Gamma(3/2-s)}\biggr)^{\!\!k}
    .
  \end{equation*}
\end{cor}
The sign in the first formula above depends on the determinant of the
scattering matrix for $Y(\lambda_0)$ at $s=1/2$. According to Lax and
Phillips \cite[Prop.~8.14]{lpscatt}, this determinant equals $(-1)^k$.
It vanishes in the second formula as the number of cusps is even and
there are no elementary components.

Now we prove Theorem \ref{main zeta}. The logarithmic derivative of
$\Zeta_\lambda$ is related with the resolvent of the Laplacian by the
trace formula, as stated in appendix \ref{trace formula}. More
precisely, if $\re(s)>1/2$ the functions
\begin{equation*}
  h_s(\xi)=\bigl(1/4+\xi^2-s(1-s)\bigr)^{-1},\quad g_s(u)=(2s-1)^{-1} e^{-(s-1/2)\abs
  u}
\end{equation*}
correspond via the Selberg transform to
\begin{equation*}
  \begin{split}
    k_s\colon\ t&\longmapsto \frac{4^{s-1}}\pi \int_{0}^{\infty} \frac{
      \left(\sqrt{w+t}+\sqrt{w+t+4}\right)^{1-2s}}{\sqrt w
      \sqrt{w+t}\sqrt{w+t+4}}\dx w\\
    &\qquad =\frac {4^{s-1}}\pi \int_0^1 (x(1-x))^{s-1}(4x+t)^{-s}\dx x.
  \end{split}
\end{equation*}
The trace formula is applicable to the difference $h_s-h_{s_0}$ if $\re(s)>1$
and $\re(s_0)>1$, and this yields for $Y=Y(\lambda)$
\begin{equation*}
  \begin{split}
    \lim_{A\to 0}\int_{Y^+}\chi_A K^0\vol &= \frac 1{2s-1} \sum_{c\in\mathcal
      C}\sum_{n=1}^{\infty} \frac{\ell(c) e^{-(s-1/2)
        n\ell(c)}}{e^{n\ell(c)/2}-e^{- n\ell(c)/2}}\\
    &\quad-\frac 1{2 s_0-1} \sum_{c\in\mathcal C}\sum_{n=1}^\infty
    \frac{\ell(c) e^{-(s_0-1/2)n\ell(c)}}{e^{n\ell(c)/2}-e^{-n\ell(c)/2}}\\
    &=\frac 1{2s-1}\, \frac{\Zeta'_\lambda(s)}{\Zeta_\lambda(s)}-\frac 1{
      2s_0-1}\,\frac{\Zeta'_\lambda(s_0)}{\Zeta_{\lambda}(s_0)}.
  \end{split}
\end{equation*}

The surface $Y^+$ is a disjoint union of $Y(\lambda)$ and one pair
$\bar Z_d$ of elementary cusps for each $d\in\graph G_1^*$ with
$\ell(d)=0$.  We add to $Y^+$ an elementary cylinder $\bar Z_d$ for
each $d\in\graph G_1^*$ with $\ell(d)\ne 0$. This gives 
\begin{equation*}
  Y^*:=Y^+\cup\bigcup_{\substack{d\in\graph G_1^*\\ \ell(d)\ne 0}}\bar Z_d
  =Y\cup\bigcup_{d\in\graph G_1^*}\bar Z_d,
\end{equation*}
where each $\bar Z_d$ is an isometric copy of the elementary cylinder
$\langle(x,a)\mapsto (x+1,a)\rangle\backslash X_{\ell(d)}$. In analogy
with the trace formula above, the logarithmic derivative of the
quotient $\Zeta_\lambda/\prod_d\Zeta_{d,\lambda}$ will be realised as
an integral over $Y^*$, so that the newly added $\bar Z_d$ yield a
logarithmic derivative of
$\Zeta_{d,\lambda}(s_0)/\Zeta_{d,\lambda}(s)$ each. On the other hand,
the surface $Y^*$ already appeared in section \ref{resolvent
  geometrically finite}, where a family of trace class operators on
$\L^2(Y)$ was defined by
\begin{equation*}
  T(s)=(\lap-s(1-s))^{-1}-(\lap-s_0(1-s_0))^{-1}-\sum_{d\in\graph
      G_1^*} \psi_d(R_d(s)-R_d(s_0))\phi_d.
\end{equation*}
The sum on the right-hand side consists of operators on the $\bar Z_d$
that are pulled back to $Y$.  The following proposition computes the
trace by integration of the operator kernel. The function $K^0$ on
$Y^*$ that appears in the statement is defined, as in appendix
\ref{trace formula}, by the truncated symmetrisation of $k_s-k_{s_0}$
over a uniformising group for each component of $Y$, and by that of
$-(k_s-k_{s_0})$ for each component of $Y^*\setminus Y$, and
$\phi_0=1-\sum_{d\in\graph G_1^*}\phi_d$.
\begin{prop}\label{trace of laplacian}
  If $\re(s)>1$ and $\re(s_0)>1$, the trace of
  \begin{equation*}
    (\lap-s(1-s))^{-1}-(\lap-s_0(1-s_0))^{-1}-\sum_{d\in\graph G_1^*}\psi_d
    (R_d(s)-R_d(s_0))\phi_d
  \end{equation*}
  is equal to the following expression:
  \begin{multline}
    \frac 1{2s-1}\biggl( \frac{\Zeta'_\lambda(s)}{\Zeta_\lambda(s)}-
    \sum_{d\in\graph G_1^*}
    \frac{\Zeta'_{d,\lambda}(s)}{\Zeta_{d,\lambda}(s)}\biggr) -\frac 1{
    2s_0-1} \biggl(\frac{\Zeta'_\lambda(s_0)}{\Zeta_\lambda(s_0)}
    -\sum_{d\in\graph G_1^*}
    \frac{\Zeta'_{d,\lambda}(s_0)}{\Zeta_{d,\lambda}(s_0)}\biggr)\\
    -\sum_{d\in\graph G_1^*}\int_{\bar Z_d} (1-\bar\phi_d) K^0\,\vol
    +\int_Y\phi_0\vol\cdot\frac 1{4\pi} \int_{-\infty}^\infty \xi\cdot
    (h_s(\xi)- h_{s_0}(\xi))\tanh(\pi\xi)\dx\xi.\label{correction terms}
  \end{multline}
\end{prop}
\begin{proof}
  The operator is known to be of trace class, and the trace can be
  computed by integrating its integral kernel over the diagonal in
  $Y\times Y$. If $\re(s)>1$ and $\re(s_0)>1$ the kernel for each
  operator in the definition of $T$ is given as the symmetrisation of
  $\pm(k_s-k_{s_0})$ for all components of $Y^*$. The overall contribution of
  the identity is
  \begin{multline*}
    \int_Y\biggl( 1-\sum_{d\in\graph G_1^*}\psi_d\phi_d\biggr)\vol\cdot
    (k_s(0)-k_{s_0}(0))\\
    =\int_Y\phi_0\vol\cdot \frac 1{4\pi}\int_{-\infty}^{\infty} \xi\cdot
    (h_s(\xi)-h_{s_0}(\xi))\tanh(\pi\xi)\dx\xi.
  \end{multline*}
  If $\ell(d)\ne 0$, the remaining contribution of $\psi_d(R_d(s)-R_d(s_0))\phi_d$ is
  \begin{equation*}
    \int_{\bar Z_d}\bar\phi_d K^0\vol
    =\frac 1{2s-1}\, \frac{\Zeta'_{d,\lambda}(s)}{\Zeta_{d,\lambda}(s)} -\frac
    1{2s_0-1}\, \frac{\Zeta'_{d,\lambda}(s_0)}{\Zeta_{d,\lambda}(s_0)}
    -\int_{\bar Z_d}(1-\bar\phi_d) K^0\vol,
  \end{equation*}
  and the trace formula, as applied at the beginning of this paragraph,
  completes the proof.
\end{proof}
We proved in theorem \ref{continuity resolvent} and in proposition
\ref{elareainf} that the trace depends continuously on
$\lambda\in\base$, and so do the correction terms in equation
\eqref{correction terms} by lemma \ref{convergence remainder trace}.
Hence we see that the logarithmic derivative of
$\Zeta_\lambda/\prod_{d\in\graph G_1^*}\Zeta_{d,\lambda}$ is
continuous in $\lambda$.  In view of equation \eqref{truncation},
proposition \ref{trace of laplacian} provides an analytic continuation
of $\Zeta_\lambda$ from its domain of convergence to the half plane
$\{s\in\IC\st \re(s)>1/2\}$. The theorem follows from proposition
\ref{prop zeta} by integration of the logarithmic derivative.

In a similar way we obtain the following observation, Wolpert's \emph{Conjecture 1}.
\begin{thm} 
  If $\lambda\in\base$ let $\mathcal N
  (\lambda,s):=\prod_{t<1/4}(t-s(1-s))$, where $t$ runs through
  the eigenvalues of the Laplacian below $1/4$. Let
  $K\subset\{s\in\IC\st \re(s)>1/2\}$ and $\inbase\subset\base$ be
  relatively compact subsets. Then there exist positive numbers
  $\alpha,\beta$ such that
  \begin{equation*}
    \alpha\le \abs{\mathcal N(\lambda,s)^{-1}\cdot\Zeta_\lambda(s)/\textstyle\prod_{d\in\graph G_1^*}Z_{d,\lambda}(s)}\le\beta
  \end{equation*}
  holds for all $s\in K$ and $\lambda\in\inbase$.
\end{thm}
\begin{proof} 
  We only need to subtract the singular part of the operator in
  proposition \ref{trace of laplacian}, to get a holomorphic family of
  operators on $\{s\in\IC\st\re(s)>1/2\}$ that is continuous in
  $\lambda$. If $t\in[0,1/4)$ is an eigenvalue of the Laplacian, the
  corresponding singular part of $s\mapsto\bigl(\lap-s(1-s)\bigr)^{-1}$ is
  \begin{equation*}
    \frac 1{t-s(1-s)}\cdot \mathrm{pr}_{t}=\frac 1{2s-1}\cdot\frac {2s-1}{t-s(1-s)}\cdot\mathrm{pr}_t\,,
  \end{equation*}
  where $\mathrm{pr}_{t}$ denotes projection onto the eigenspace. Now
  $\frac{2s-1}{t-s(1-s)}$ is the contribution of a particular eigenvalue to the
  logarithmic derivative of $\mathcal N(\lambda,s)$.
\end{proof}
In order to examine the Zeta function on the left, we need the
precise asymptotics of the contribution of a single geodesic as its
length decreases.
\begin{lem}\label{asymptotic} Define $\Zeta_l(s):=\prod_{k=0}^\infty \bigl(1-e^{-(s+k)l}\bigr)^2$ for all $l>0$. Then
  \begin{equation*}
    \lim_{l\to 0} \Gamma(s)^2\Zeta_l(s) e^{\pi^2/3l} l^{2s-1}=2\pi.
  \end{equation*}
\end{lem}
\begin{proof} It is sufficient to prove the lemma if $\re(s)>0$, for it can be continued iteratively to the left using the functional equation of $\Gamma$. Like Wolpert and Hejhal we begin with  
  \begin{align*}
    \log\Zeta_l(s)&=2\sum_{k=0}^\infty\log\bigl(1-e^{-(s+k)l}\bigr) =-2\sum_{k=0}^{\infty}\sum_{n=1}^\infty\frac{ e^{-(s+k)nl}}n\\
    &=-2\sum_{n=1}^\infty \frac{e^{-snl}}n\bigl(1-e^{-nl}\bigr)^{-1}.
  \end{align*}
  This expression is split into three part according to
  \begin{align*}
    \bigl(1-e^{-nl}\bigr)^{-1} &=\frac 12+\frac 1{nl}+ \Bigl[\bigl( 1-e^{-nl}\bigr)^{-1}-\frac 1{nl}-\frac 12\Bigr]\\
    &= \frac 12+\frac 1{nl} +\Bigl[\bigl(e^{nl}-1\bigr)^{-1}-\frac 1{nl}+\frac 12\Bigr].
  \end{align*}
  The contribution of the first term to $\log\Zeta_l(s)$ is
  \begin{equation}\label{log}
    -\sum_{n=1}^\infty \frac{e^{-snl}}n= \log\bigl(1-e^{-sl}\bigr).
  \end{equation}
  The second summand yields Euler's dilogarithm:
  \begin{equation*}
    -2l^{-1}\sum_{n=1}^\infty \frac{e^{-snl}}{n^2}=-2l^{-1}\Li\bigl(e^{-sl}\bigr),
  \end{equation*}
  and the dilogarithm satisfies $\Li(z)=\pi^2/6-\log (z)\cdot\log(1-z)-\Li(1-z)$, so
  \begin{align}\label{Euler}
    -2l^{-1}\sum_{n=1}^\infty \frac{e^{-snl}}{n^2}&= -\frac{\pi^2}{3l}-2s\log\bigl(1-e^{-sl}\bigr)+2l^{-1}\Li\bigl(1-e^{-sl}\bigr)\\\nonumber
    &\sim -\frac{\pi^2}{3l}-2s\log\bigl(1-e^{-sl}\bigr)+2s,\quad l\to 0.
  \end{align}
  Finally, as $l\to 0$, we have \cite[p.~21, eq.~(4)]{erdelyi1}
  \begin{multline}\label{gamma}
    -2l\sum_{n=1}^\infty \frac{e^{-snl}}{nl}\Bigl[\bigl(e^{nl}-1\bigr)^{-1}-\frac 1{nl}+\frac 12\Bigr]\\
  \begin{aligned}\longrightarrow &
  -2\int_0^\infty t^{-1} e^{-st}\Bigl[\bigl( e^t-1\bigr)^{-1}-\frac 1t+\frac 12\bigr]\dx t\\[1em]
  &=-2\log\Gamma(s)+(2s-1)\log s-2s+\log (2\pi).
  \end{aligned}
\end{multline}
  The sum of \eqref{log},\eqref{Euler},\eqref{gamma} gives the result.
\end{proof}
Now we want to consider the Zeta function on the left of the critical
axis. The theorem below relies on the asymptotic of the approximate
scattering matrix given in proposition \ref{D left}. We must therefore
restrict to a degenerating family of compact surfaces.
\begin{thm}\label{zeta everywhere} Let $\graph G$ be a graph of type $(p,0)$ and $\lambda_0=(\ell_0,\tau_0)\in\partial\base$. Let $S$ be the set of edges $d$ with $\ell_0(d)=0$. If $(\lambda_m)$ is a sequence in $\interior\base$ that converges to $\lambda_0$, then
  \begin{equation*}
    \lim_{m\to\infty} \frac{\det D(\lambda_m,s)\cdot\Zeta_{\lambda_m}(s)}{\prod_{d\in S}\Zeta_{d,\lambda_m}(s)}=\Zeta_{\lambda_0}(s)
  \end{equation*}
  holds on $\{s\in\IC\st\re(s)\ne 1/2\}$.
\end{thm}
\begin{proof}
  The statement is an extension of theorem \ref{main zeta} to the left
  of the critical axis: the matrices $D(\lambda_m,s)$ converge to the
  identity matrix if $\re(s)>1/2$, so multiplication of the quotient
  $\Zeta(\lambda_m,s)/\prod_d\Zeta_{d,\lambda_m}$ with $\det
  D(\lambda_m,s)$ does not alter the limit in this domain. 
  
  Then the proof consists in the application of a number of formulas
  stated earlier.  Let $k$ denote the number of cusps in $Y_{\graph
    G}(\lambda_0)$. By corollary \ref{zeta left} we have for all $s$
  with $\re(s)>1/2$
  \begin{equation*}
    \lim_{m\to\infty}\frac{\Zeta_{\lambda_m}(1-s)}{\textstyle{\prod_{d\in S}}\Zeta_{d,\lambda_m}(s)}=  \Zeta_{\lambda_0}(1-s)\, \det( C(\lambda_0,1-s))\cdot 2^{(2s-1)k}\biggl(\frac{\Gamma(1/2+s)}{\Gamma(3/2-s)}\biggr)^{\!\!k}.
  \end{equation*}
  In corollary \ref{approximation 2} we proved $\lim_m
  D(\lambda_m,1-s)^{-1}C(\lambda_m,1-s)=C(\lambda_0,1-s)$, so
  \begin{equation}\label{X}
    \lim_{m\to\infty} \frac{\det D(\lambda_m,1-s)}{\det C(\lambda_m,1-s)}\, \frac {\Zeta_{\lambda_m}(1-s)}{\prod_d\Zeta_{d,\lambda_m}(s)}=
     \Zeta_{\lambda_0}(1-s)\,\frac{\Gamma(1/2+s)^k}{\Gamma(3/2-s)^k}\, 2^{(2s-1)k}.
  \end{equation}
  Proposition \ref{D left} gives
  \begin{equation*}
    \det C(\lambda_m,1-s)\sim \biggl(\prod_{d\in\tilde S}\ell(d)^{2s-1}\biggr) (2s-1)^k\frac {4^{-ks}\,\Gamma(s)^{2k}}{\Gamma(1/2+s)^{2k}}\det\Sigma(s),
  \end{equation*}
  and
  \begin{equation*}
    \det\Sigma(s)=\Bigl(\frac 1{\cos(\pi s)}-1\Bigr)^{k/2}=\frac {\Gamma(1/2+s)^k\, \Gamma(1/2-s)^k}{\Gamma(s)^k\, \Gamma(1-s)^k}
  \end{equation*}
  implies
  \begin{equation*}
    \det C(\lambda_m,1-s)\sim \biggl(\prod_{d\in\tilde S}\ell(d)^{2s-1}\biggr) (2s-1)^k \frac{4^{-ks}\Gamma(s)^k}{\Gamma(1/2+s)^k}\, \frac{\Gamma(1/2-s)^k}{\Gamma(1-s)^k}.
  \end{equation*}
  We substitute this into \eqref{X} to see
  \begin{equation*}
    \lim_{m\to\infty} \frac{\det D(\lambda_m,1-s) \Zeta_{\lambda_m}(1-s)\Gamma(1-s)^k}{\prod_{d\in S}\bigl(\Zeta_{d,
        \lambda_m}(s)\,\ell(d)^{4s-2}\bigr)\Gamma(s)^k} = \Zeta_{\lambda_0}(1-s).
  \end{equation*}
  Thus the proof is reduced to
  \begin{equation*}
    \frac{ \Zeta_{d,\lambda_m}(s)\,\ell(d)^{4s-2}\,\Gamma(s)^2}{\Zeta_{d,\lambda_m}(1-s)\,\Gamma(1-s)^2}\longrightarrow 1,
  \end{equation*}
  which follows from lemma \ref{asymptotic}.
\end{proof}

\appendix
\section{Families of Fuchsian groups}\label{families}
Throughout this appendix, the hyperbolic plane is identified with its
unit-disc model $D$. By $\dist{z_1}{z_2}$ we mean the hyperbolic
distance between two points. Let $\Gamma$ denote a finitely generated
group without torsion that admits a continuous family $\phi\colon
\base\to\hom(\Gamma,\isom (D))$ of discrete inclusions into the
orientation-preserving isometries of $D$.   The image of
$\phi(b)$ is denoted by $\Gamma_b$, and the image of any
$\gamma\in\Gamma$ under $\phi(b)$ by $\gamma_b$. By continuity
we mean that $b\mapsto\gamma_b$ is a continuous map for each
$\gamma\in\Gamma$.  This is equivalent to the continuity of
$b\mapsto\Gamma_b$, where the set of closed subgroups of $\isom(D)$
carries the Chabauty topology. We assume $\base$ to be a
path-connected, locally compact and metrisable space.

We assume that the $\Gamma_b$ are non-elementary. The only element of
finite order in $\Gamma_b$ is the identity. An action of $\Gamma$ on
$\base\times D$ is defined by $\gamma(b,z):=(b,\gamma_b z)$. Our first
observation is that this action is freely discontinuous, i.e.~every
element of $\base\times D$ has a neighbourhood $U$ such that $\gamma
U\cap U\ne\emptyset$ implies $\gamma=1$.
\begin{lem}\label{discont} The induced action of $\Gamma$ on $\base\times D$
  is freely discontinuous.
\end{lem}
\begin{proof}  We use the fact that the inequality
  \begin{equation*} \sinh(\dist z{\gamma_1 z}/2)\,\sinh (\dist z{\gamma_2 z}/2)\ge 1\end{equation*}
  is satisfied for any $z\in D$ if the group generated by
  $\gamma_1,\gamma_2\in\isom(D)$ is torsion-free and non-elementary
  (Beardon \cite[\S 8.3, p.~198]{beardon}). 
  
  Let $(b_0,z_0)\in \base\times D$, the aim is to prove that $\Gamma$
  acts freely discontinuously in $(b_0,z_0)$. Choose $\tau\in
  \Gamma\setminus\{1\}$. Since $\base$ is locally compact, there exist
  neighbourhoods $V\subset \base$ of $b_0$ and $W\subset D$ of $z_0$
  such that $\dist{\phi(b)(\tau) z}{z}>\epsilon >0$ for all $(b,z)\in
  V\times W$. The cited inequality implies $\dist{\phi(b)(\gamma)
    z}{z}> 2 \arsinh(\epsilon^{-1})$ for $(b,z)\in V\times W$ if
  $\gamma$ is not contained in the centraliser $Z(\tau)$. Now the
  images of $Z(\tau)$ under the homomorphisms considered are
  elementary and discrete, so we may choose a neighbourhood $U'$ of
  $(b_0,z_0)$ with the desired property for $Z(\tau)$. Then
  \begin{equation*}
    U'\cap \left( V\times\left\{ z\in D\st \dist{z}{z_0}<\arsinh (\epsilon^{-1})\right\}\cap W\right)
  \end{equation*}
  is a suitable neighbourhood of $(b_0,z_0)$.
\end{proof}
The previous lemma allows for a definition of $\Gamma$-invariant structures on
$\base\times D$, e.g.~Riemannian metrics or vector fields if the
action of $\Gamma$ is differentiable. We will stick to the latter in
order to define locally a trivialisation of the quotient
$\Gamma\backslash \total\to \base$, where $\total$ denotes the fibre-wise Nielsen domain for the group action. We must restrict to the Nielsen domain as we explicitly want to allow the type of an isometry $\gamma_b$ to change from parabolic to hyperbolic and vice versa. 

Of particular importance in this respect are those primitive elements
of a Fuchsian group that are related to the infinite parts of the
corresponding quotient of the hyperbolic plane. They are specified in
the following definition.
\begin{dfn}\label{boundary pairing}
  An element $\gamma\in\Gamma\setminus\{1\}$ is called a \emph{boundary pairing} if it is not a
  proper multiple of any other element of $\Gamma$, and if there exists $b\in
  \base$ such that $\gamma_b$ either is parabolic or it leaves invariant a
  component of the set of discontinuity in $\partial D$.
\end{dfn}
The definition refers to a property of $\gamma_b$ for some particular
$b\in \base$, but the next lemma shows that it is in fact independent
of $b$.
\begin{lem}
  If $\gamma$ is a boundary pairing, then the condition in definition
  \ref{boundary pairing} is satisfied by $\gamma_b$ for all $b\in
  \base$.
\end{lem}
\begin{proof} Suppose that the condition is satisfied by $\gamma_{b_0}$ but
  not by $\gamma_{b_1}$. There exists a curve $t\mapsto b_t$ from $b_0$ to
  $b_1$, and two curves $c_1,c_2$ in $\partial D$ such that $c_1(t)$ and
  $c_2(t)$ are the fixed points of $\gamma_{b_t}$ for all $t$. The fixed
  points of non-trivial elements of $\Gamma_{b_1}$ are dense in the limit
  set. By assumption, $\gamma_{b_1}$ is hyperbolic and each component of
  $\partial D\setminus \{c_1(1), c_2(1)\}$ contains a fixed point of an
  element $\eta_{b_1}$ or $\eta'_{b_1}$. Now $\partial D\setminus \{c_1(0),
  c_2(0)\}$ either consists of one component, or one of its components is a
  subset of the set of discontinuity. This implies that the intersection of
  $\{c_1(t), c_2(t)\}$ with the set of fixed points of $\eta_{b(t)}$ or
  $\eta'_{b(t)}$ is not empty for some $t$. This is a contradiction to discreteness
  of $\Gamma_{b(t)}$.
\end{proof}
Now we want to construct a trivialisation of $\Gamma\backslash
\total\to \base$ from an invariant vector field $s$ on $\base\times
D$. For convenience we restrict to a smooth family of Fuchsian groups
parametrised by a compact interval $\base$. Let $\Gamma'\subset\Gamma$
be the set of boundary pairings.

\begin{figure}
  \begin{center}
\begin{picture}(0,0)%
\includegraphics{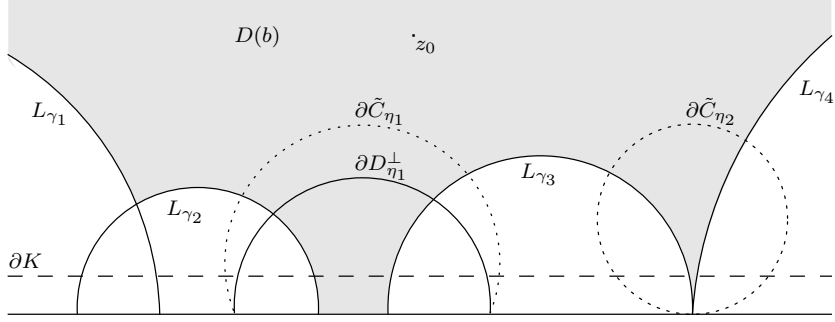}%
\end{picture}%
\setlength{\unitlength}{4144sp}%
\begingroup\makeatletter\ifx\SetFigFont\undefined
\def\x#1#2#3#4#5#6#7\relax{\def\x{#1#2#3#4#5#6}}%
\expandafter\x\fmtname xxxxxx\relax \def\y{splain}%
\ifx\x\y   
\gdef\SetFigFont#1#2#3{%
  \ifnum #1<17\tiny\else \ifnum #1<20\small\else
  \ifnum #1<24\normalsize\else \ifnum #1<29\large\else
  \ifnum #1<34\Large\else \ifnum #1<41\LARGE\else
     \huge\fi\fi\fi\fi\fi\fi
  \csname #3\endcsname}%
\else
\gdef\SetFigFont#1#2#3{\begingroup
  \count@#1\relax \ifnum 25<\count@\count@25\fi
  \def\x{\endgroup\@setsize\SetFigFont{#2pt}}%
  \expandafter\x
    \csname \romannumeral\the\count@ pt\expandafter\endcsname
    \csname @\romannumeral\the\count@ pt\endcsname
  \csname #3\endcsname}%
\fi
\fi\endgroup
\begin{picture}(4956,1921)(432,-3320)
\put(2521,-2131){\makebox(0,0)[lb]{\smash{\SetFigFont{12}{14.4}{rm}{$\scriptstyle \partial\tilde{C}_{\eta_1}$}%
}}}
\put(1801,-1681){\makebox(0,0)[lb]{\smash{\SetFigFont{12}{14.4}{rm}{$\scriptstyle D(b)$}%
}}}
\put(586,-2131){\makebox(0,0)[lb]{\smash{\SetFigFont{12}{14.4}{rm}{$\scriptstyle L_{\gamma_1}$}%
}}}
\put(5176,-1996){\makebox(0,0)[lb]{\smash{\SetFigFont{12}{14.4}{rm}{$\scriptstyle L_{\gamma_4}$}%
}}}
\put(451,-3031){\makebox(0,0)[lb]{\smash{\SetFigFont{12}{14.4}{rm}{$\scriptstyle \partial K$}%
}}}
\put(2881,-1726){\makebox(0,0)[lb]{\smash{\SetFigFont{12}{14.4}{rm}{$\scriptstyle z_0$}%
}}}
\put(1396,-2716){\makebox(0,0)[lb]{\smash{\SetFigFont{12}{14.4}{rm}{$\scriptstyle L_{\gamma_2}$}%
}}}
\put(2521,-2446){\makebox(0,0)[lb]{\smash{\SetFigFont{12}{14.4}{rm}{$\scriptstyle \partial D_{\eta_1}^\perp$}%
}}}
\put(3511,-2491){\makebox(0,0)[lb]{\smash{\SetFigFont{12}{14.4}{rm}{$\scriptstyle L_{\gamma_3}$}%
}}}
\put(4501,-2131){\makebox(0,0)[lb]{\smash{\SetFigFont{12}{14.4}{rm}{$\scriptstyle \partial\tilde{C}_{\eta_2}$}%
}}}
\end{picture}

\caption{Trimming of a Dirichlet domain}\label{flow}
\end{center}
\end{figure}

Let us recall the definition of $\total$. If $\gamma\in\Gamma'$ and
$\gamma_b$ is hyperbolic, its axis separates $D$ into two components.
Precisely one of them, denoted by $D_\gamma^\perp (b)$, is not bounded
at $\partial D$ by a component of the set of discontinuity. If
$\gamma_b$ is parabolic, we define $D_\gamma^\perp (b)$ to be empty.
Then the fibre-wise Nielsen domain is the following $\Gamma$-invariant subset of
$\base\times D$:
\begin{equation*}
  \total:= \bigcup_{b\in\base} \bigcap_{\gamma\in\Gamma'} \{b\}\times
  D_\gamma^\perp (b).
\end{equation*}
The projection $\Gamma\backslash\total\to\base$ is not a proper map if
any of the Fuchsian groups contain parabolic elements. If we want to
apply the flow of a vector field to trivialise this map, we must
impose additional conditions on the vector field. These ensure that
for each compact neighbourhood $\inbase\subset\interior\base$ of a point $b_0$
there exists a relatively compact set $K\subset D$ such that the
following holds:
\begin{itemize}
\item If $\Phi_b\colon \inbase\times K\to \base\times D$, $\abs
  {b-b_0}<\epsilon$, denotes the flow of the vector field, then $\total$
  is invariant under $\Phi_t$.
\item The map $\Phi_b\vert_{\{b_0\}\times K}$ extends to a
  diffeomorphism $\total_{b_0}\to\total_b$ for all $\abs
  {b-b_0}<\epsilon$.
\end{itemize}
Essentially this means that we must require the vector field to be
well-behaved near the fibre-wise boundary of
$\Gamma\backslash\total\subset\Gamma\backslash(\base\times D)$. We
give a precise definition of a neighbourhood of this boundary in terms
of covering subsets $\tilde C_\gamma$ of collars (figure~\ref{flow}):

There is an open subset $\tilde C_\gamma$ of $\total$ associated with
each $\gamma\in\Gamma'$. If $\ell(\gamma_b)$ denotes the hyperbolic
translation of $\gamma_b$, resp.~$\ell(\gamma_b)=0$ if $\gamma_b$ is
parabolic, there exists an orientation preserving isometry
$X_{\ell(\gamma_b)}^-\to D_\gamma^\perp(b)$ such that the isometry
$(x,a)\mapsto (x\pm 1,a)$ of $X_{\ell(\gamma_b)}^-\to
D_\gamma^\perp(b)$ corresponds to $\gamma_b$ (cf.~p.~\pageref{model}).
The isometry is uniquely determined up to composition on the left with
maps $(x,a)\mapsto (x+\delta,a)$. So we can define the fibre of
$\tilde C_\gamma$ over $b$ to be the the image of
$X_{\ell(\gamma_b)}^{A(\ell(\gamma_b))}$, where
\begin{equation*}
  A(\ell)=
  \begin{cases}
    \bigl( -\frac{\ell}{2\sinh (\ell/2)},0\bigr) &\text{if $\ell\ne 0$}\\
    (-1,0)&\text{if $\ell=0$}.
  \end{cases}
\end{equation*}

If $\ell_0$ satisfies $\ell_0<\ell(\gamma_b)$ for all
$b\in\base$, we define $\tilde C_\gamma^{\ell_0}\subset\tilde
C_\gamma$ as the fibre-wise image of $X_{\ell(\gamma_b)}^{A(\ell_0)}$.
The Fuchsian groups are required to be geometrically finite,
so there are only finitely many conjugacy classes of boundary
pairings. This implies that such a number $\ell_0$ can be chosen
independently of $\gamma\in\Gamma'$.
\begin{lem}
  Let $s$ be a smooth, $\Gamma$-invariant lift of the canonical vector
  field $\partial_b$ on the interval $\base$. Provided that the
  integral flow of $s$ defines a diffeomorphism of the fibres of \/
  $\bigcup_{\gamma}\tilde C_\gamma^{\ell_0}\to \base$ for some
  $\ell_0>0$, it extends these diffeomorphisms to a local
  trivialisation of $\total\to \base$.
\end{lem}

\begin{proof}
  We verify that the two properties above are satisfied.
  Invariance of $\total$ is immediate from the condition on $s$.
  
  The set $K$ will be defined in terms of Dirichlet domains for the
  Fuchsian groups. Choose $z_0\in D$. If $\gamma\in\Gamma$ and
  $b\in\base$, the open half-plane
  \begin{equation*}
    D_\gamma(b) := \{z\in D\st \dist{z}{z_0}< \dist{z}{\gamma_b z_0}\}
  \end{equation*}
  is bounded by the geodesic
  \begin{equation*}
    L_\gamma(b):= \{z\in D\st \dist{z}{z_0}= \dist{z}{\gamma_b z_0}\}.
  \end{equation*}
  The Dirichlet domain $D(b)$ of $\Gamma_b$ with centre $z_0$ is
  $\bigcap_\gamma D_\gamma(b)$. Due to finite geometry, for fixed
  $b\in B$ this intersection can by replaced with an intersection of
  finitely many half-planes, say of those associated with
  $\gamma_1,\ldots,\gamma_n\in \Gamma$.
  
  We need to have a closer look at those points where the Euclidean
  closure of $D(b)$ meets $\partial D$. This only happens in a
  component of the set of discontinuity or in a parabolic fixed point
  for the following reason: The limit set of a geometrically finite
  group only consists of conical limit points and parabolic fixed
  points. But a convex fundamental polygon of such a group cannot meet
  $\partial D$ in a point of approximation \cite[thm 10.2.3 and thm
  10.2.5]{beardon}.  This observation implies the existence of
  finitely many boundary pairings $\eta_1,\ldots,\eta_j$ and of a
  compact neighbourhood $K$ of $z_0$ such that
  \begin{equation}\label{dirichlet compact}
    D(b)\setminus K\subset \bigcup_{k=1}^j \left( D_{\eta_k}^\perp(b) \cup
      \tilde C_{\eta_k}^{\ell_0}\right)
  \end{equation}
  (this is illustrated in figure \ref{flow}). Moreover, the Dirichlet
  domain $D(b')$ is always contained in the intersection
  $\bigcap_{k=1}^n D_{\gamma_k}(b')$, so the compact set $K$ may be
  chosen such that \eqref{dirichlet compact} holds if $b$ is replaced
  with any element $b'$ of a suitable neighbourhood of $b$. If
  $\inbase\subset\base$ is a compact neighbourhood with this property, we
  only need to apply the integral flow to $\inbase\times K$ and to each
  $\tilde C_\gamma^{\ell_0}$ separately.
\end{proof}
\section{The Selberg trace formula}\label{trace formula}
The trace formula as formulated by Selberg \cite{selberg1} expresses
the trace of certain operators $\tilde h(\lap)$ in terms of the closed
geodesics on a geometrically finite surface.  Difficulties may arise
since the spectrum of $\lap$ is not discrete for non-compact surfaces.
In this appendix, we will only state a rudimentary version of this
trace formula in the sense that it computes the distributional trace
of an operator by integration of its integral kernel over the
diagonal, but the result is not known to be related to the trace of a
trace class operator a priori. Such a relation is established in a
special case only, in section \ref{zeta function}, where the operator
is some kind of relative resolvent of the Laplacian (cf.
p.~\pageref{laplacian} for the definition).

For suitable functions $\tilde h$, the Selberg transform provides an
explicit formula for the integral kernel of $\tilde h(\lap)$. The
kernel on the hyperbolic plane is given by
\begin{equation*}
  (z_1,z_2)\ \longmapsto\ k(4 \sinh^2(\dist{z_1}{z_2}/2)),
\end{equation*}
where $k$ is related to the function $h\colon\xi\mapsto \tilde h(1/4+\xi^2)$ by
\begin{gather*}\label{selberg transform}
  k(t)=-\frac 1\pi \int_t^{\infty} \frac{\dx Q(w)}{\sqrt{w-t}},\qquad
  \int_w^{\infty} \frac{k(t)}{\sqrt{t-w}}\,\dx t=Q(w),\\[1em]
  Q(e^u+e^{-u}-2)=g(u),\\[.5em]
  g(u)=\frac 1{2\pi} \int_{-\infty}^{\infty} h(\xi) e^{-i\xi u}\dx \xi,
  \qquad  \int_{-\infty}^{\infty} g(u) e^{i\xi u}\dx u=h(\xi).
\end{gather*}
For the time being, we only consider those operators with an integral
kernel that can be derived from $k$ by summation over a uniformising
group. To specify a class of functions $h$ that meet this condition,
we prove the following lemma. It is certainly well-known, but to our
knowledge the complete proof is only implicit in the literature, so we
give it here for the sake of completeness (c.f. Hejhal \cite{hej1}).
\begin{lem}\label{bound}
  Let $h$ be an analytic function in $\left\{\xi\in\IC \st \abs{\im
  (\xi)}<1/2+\delta\right\}$ satisfying $h(\xi)=h(-\xi)$ and $\abs{h(\xi)} \le
  M\, (1+\abs{\re(\xi)}^2)^{-(1+\alpha)}$. Then $k$ is continuous on
  $\left[0,\infty\right)$ and there exist numbers $C_\rho$ for all
  $\rho<\delta$ such that
  \begin{equation}
    \abs{g(u)}\le C_\rho e^{-(1/2+\rho)}\abs u\quad\text{and}\quad \abs{k(t)} \le
    C_\rho (1+t)^{-(1+\rho)}.
  \end{equation}
  Conversely, if we start with a continuous function $k$ on
  $\left[0,\infty\right)$ with $\abs{k(t)}\le C_\rho
  (1+t)^{-(1+\rho)}$ for $\rho<\delta$ and we apply the inverse
  transformations, then $g$ satisfies estimates as above and $h$ is
  analytic in the strip.
\end{lem}
\begin{proof}
  The relation between $g$ and $h$ is an immediate consequence of the
  Cauchy integral formula. The final assertion is also easy to see:
  Given $k$, the function $Q$ is continuous and the estimate for $g$
  follows from
  \begin{equation*}
    \int_w^{\infty} \frac{(1+t)^{-(1+\rho)}}{\sqrt{t-w}}\dx t\le
    (1+w)^{-(1/2+\rho-\epsilon)} \int_0^{\infty} (1+t+w)^{-(1/2+\epsilon)}
    t^{-1/2}\dx t.
  \end{equation*}
  So we proceed with the main part. A function $h$ that meets the
  assumptions gives rise to an element $g$ of the Sobolev space
  $\H^{3/2+\alpha}(\IR)$, because the following function is
  square-integrable:
  \begin{equation*}
    \xi\longmapsto (1+\xi^2)^{1/2 (3/2+\alpha)} h(\xi).
  \end{equation*}
  Now $3/2+\alpha > 1+(\dim\IR)/2$ holds, so the Sobolev embedding
  theorem implies $g\in\C^{1,\gamma}(\IR)$ for $\gamma<\alpha$.
  H\"older continuity of the derivative $g'$ implies continuity of $k$
  at 0:
  \begin{equation*}
    k(t)= -\frac 1\pi \int_0^\infty
    \frac{g'(\arcosh(1+{\textstyle\frac{w+t}2}))}{\sqrt{w(w+t)(w+t+4)}}\dx w,
  \end{equation*}
  and using $g'(0)=0$, the integrand can be estimated near 0 by
  \begin{equation*}
    \frac{\left(\arcosh(1+{\textstyle\frac{w+t}2})\right)^\gamma}
    {\sqrt{w(w+t)(w+t+4)}}\le \frac{\sqrt{w+t}^\gamma}{\sqrt{w(w+t)(w+t+4)}} \le
    w^{-1+\gamma/2}.
  \end{equation*}
  To derive the bounds for $k$, we use the fact that
  \begin{equation*}
    f\colon\ \left\{\xi\in\IC\st \abs{\im\xi}<
      1/2+\delta\right\}\longrightarrow\IC,\quad \xi\longmapsto -i\xi h(\xi)
  \end{equation*}
  is an analytic function which satisfies
  \begin{equation*}
    \sup_{\abs{\eta}<1/2+\delta} \norm {f(\cdot+ i\eta)}_{\L^2(\IR)}<\infty.
  \end{equation*}
  The Fourier transform of $f$ is $g'$. Therefore (cf.~Reed-Simon
  \cite[p.~18]{reed-simon2}) the function $u\mapsto e^{b\abs u} g'(u)$ is
  square-integrable for any $b<1/2+\delta$, and we see
  \begin{equation*}
    \begin{split}
      k(t)&= -\frac 1\pi \int_t^\infty \frac{Q'(w)\dx w}{\sqrt{w-t}}\\
      &= -\frac 1\pi \biggl( \int_t^{t+t^{-2\delta}}\frac{ Q'(w)\dx
        w}{\sqrt{w-t}} +\int_{\arcosh(1+(t+t^{-2\delta})/2)}^\infty
      \frac{g'(u)\dx u}{\sqrt{2(\cosh(u)-1)-t}}\biggr).
    \end{split}
  \end{equation*}
  We use boundedness of $g'$ for the first integral, and the second is
  estimated with the Cauchy-Schwarz inequality for some
  $b=b_1+b_2<1/2+\delta$:
  \begin{equation*}
    \begin{split}
      \abs{k(t)}&\le  \frac{\norm{g'}_\infty}\pi \int_{t}^{t+t^{-2\delta}}
      ((w-t)w(w+4))^{-1/2}\dx w\\
      &\quad +\frac 1\pi \biggl( \int_{\arcosh(1+(t+t^{-2\delta})/2)}^\infty e^{-(1-2
        b_1)u} \abs{g'(u)}^2\dx u\biggr)^{1/2}\\
      &\qquad\cdot \biggl( \int_{\arcosh
        (1+(t+t^{-2\delta})/2)}^\infty \frac{e^{(1-2 b_1)u}\dx
        u}{2(\cosh(u)-1)-t}\biggr)^{1/2}\\
      &\le \frac{\norm {g'}_\infty}{\pi} t^{-1} \int_0^{t^{-2\delta}} w^{-1/2}\dx
      w\\
      &\quad + \frac 1\pi e^{-(1/2+b_2)\arcosh (1+(t+t^{-2\delta})/2)} \biggl(
      \int_{\arcosh(1+(t+t^{-2\delta})/2)}^\infty e^{2bu}\abs{g'(u)}^2\dx
      u\biggr)^{1/2}\\
      &\qquad\cdot \biggl(\int_{t+t^{-2\delta}}^\infty ((w-t)\sqrt w\sqrt{w+4})^{-1}
      e^{(1-2b_1) \arcosh (1+w/2)}\dx w\biggr)^{1/2}\\
      &\le \frac{2\norm {g'}_\infty}\pi t^{-(1+\delta)} +C (1+t)^{-(1/2+
        b_2)}.
    \end{split}
  \end{equation*}
  The last line follows from $e^{\arcosh v}=v+\sqrt{v^2-1}$. This proves
  the bound for $k$ as $b_2$ may be chosen arbitrarily close to $1/2+\delta$.
\end{proof}
 
Now assume that $k\in\C([0,\infty))$ satisfies $\abs{k(t)}\le C_\rho
(1+t)^{-(1+\rho)}$ for all $\rho<\delta$. Let $Y'$ be a connected,
geometrically finite surface. We choose a uniformising group
$\Gamma'\subset\isom(D)$ and identify $Y'$ with
$\Gamma'\backslash D$. The critical exponent of a Fuchsian group is
less than or equal to one, and therefore the series
\begin{equation*}
  K_{\Gamma'}(z_1,z_2):=\sum_{\gamma\in\Gamma'} k(4 \sinh^2(\dist{z_1}{\gamma
  z_2}/2))
\end{equation*}
is absolutely convergent for $z_1,z_2\in D$. It defines a
$\Gamma'\times\Gamma'$-invariant function, so let $K_{Y'}\colon
Y'\times Y'\to\IC$ denote the function induced on the quotient. We
would like to integrate $K_{Y'}$ over the diagonal in $Y'\times Y'$,
but the integral might diverge for two reasons, either due to the
growth of the kernel (if the surface has cusps) or because of the
surface having infinite area (in the presence of funnels). We will
simply subtract the divergent parts.

Consider $Y'$ to be a component of $Y=Y_{\graph G}(\lambda)$ as defined
in section \ref{surface geometry}.  An edge $d\in\graph G_1^*$ with
$\ell(d)=0$ represents a pair of cusps $Z_d\subset Y$. We adjoin an
extra pair of cusps $\bar Z_d$, which is by definition an isometric
copy of the elementary quotient $\langle (x,a)\mapsto
(x+1,a)\rangle\backslash X_0$. The resulting surface is a disjoint
union
\begin{equation*}
  Y^+:= Y\cup \bigcup_{\substack{d\in\graph G_1^*\\\ell(d)=0}}\bar Z_d,
\end{equation*}
and we will integrate over $Y^+$ the function $K^0$ defined component-wise by
\begin{equation*}\label{selberg kerneldef}
  K^0(z):= 
  \begin{cases} K_{Y'}(z,z)-k(0)&\text{if $z\in Y'\subset Y$,}\\
    -(K_{Y'}(z,z)-k(0))&\text{if $z\in\bar Z_d\subset Y^+\setminus Y$.}
  \end{cases}
\end{equation*}
The constant $k(0)$ is subtracted to avoid difficulties with
integration over infinite-area surfaces. Finally, we need a family of
cut-off functions $\chi_A\colon Y^+\to\{0,1\}$ that vanish in the
cuspidal ends of area $A$:
\begin{equation*}
  \chi_A(z)=0\quad\Longleftrightarrow\quad
    \displaystyle z=(x,a)\in\bigcup_{\ell(d)=0} Z_d\cup \bar Z_d\quad
    \text{and $\abs a<A$.}
\end{equation*}
\begin{thm}[Selberg]
  Let $\mathcal C$ denote the set of oriented, simple, closed geodesics in
  $Y$, and $\ell(c)$ the length of a curve $c\in\cal C$. Then the following
  formula holds:
  \begin{equation}\label{selberg}
    \lim_{A\rightarrow 0}\int_{Y^+}\chi_A K^0\vol = \sum_{c\in\cal C}
    \sum_{n=1}^{\infty}\frac{ \ell(c) g(n\ell(c))}{ e^{n\ell(c)/2}-
    e^{-n\ell(c)/2}}.
  \end{equation}
\end{thm}
The proof is given in Selberg's G\"ottingen lectures, in Hejhal's
monograph~\cite{hej2} or by Venkov~\cite{venkov}, one only needs to
pay attention to the following facts.

Assume the function $k$ to be positive first. For every component $Y'$ of $Y$ we fix a
complete set $B'$ of representatives for the conjugacy classes of
hyperbolic elements in $\Gamma'$. The famous calculation of Selberg's
\cite{selberg1,hej2}, in connection with Beppo Levi's theorem, 
shows that
\begin{multline}\label{sef1}
  \sum_{Y'\subset Y}\int_{F} \sum_{\substack{\gamma\in
    \Gamma'\\\gamma\ne 1}} \tilde\chi_A(z) k(\sig(z,\gamma z))\vol(z) -\sum_{Y'\subset Y}
  \int_{F} \sum_{\substack{\gamma\in\Gamma'\\
      \text{parabolic}}}\tilde\chi_A(z) k(\sig(z,\gamma z))\vol(z)\\
  \longrightarrow\ \sum_{Y'\subset Y}\sum_{\gamma\in B'}
  \frac{\ell(\gamma_0)g(\ell(\gamma))}{
    e^{\ell(\gamma)/2}-e^{-\ell(\gamma)/2}},\quad A\rightarrow 0,
\end{multline}
where $\gamma_0$ denotes the unique primitive element of $\Gamma'$ such that
$\gamma=\gamma_0^n$ holds for some positive $n$, and $\ell(\gamma)$ is the
hyperbolic translation of $\gamma$. The domain $F$ of integration is a suitable
fundamental domain for $\Gamma'$ and $\sig(z_1,z_2)=
4\sinh^2(\dist{z_1}{z_2}/2)$. The infinite series on the right-hand side converges since
\begin{equation*}
  \#\{\gamma\in B'\st e^{\ell(\gamma)}\le x\}=\Ord(x),\
   x\rightarrow\infty
\end{equation*}
and $g$ satisfies $g(u)\le C e^{-(1/2+\rho)\abs u}$ by lemma \ref{bound}.

If $k$ is not positive, we see that Lebesgue's theorem on dominated convergence implies that the positivity
condition on $k$ may be dropped in equation \eqref{sef1}. On the other hand
(see Hejhal \cite[pp.~198, 311]{hej2})
\begin{equation}\label{sef2}
  \lim_{A\rightarrow 0} \biggl( \sum_{\ell(d)=0} \int_{\bar Z_d}
  \chi_AK^0\vol + \sum_{Y'\subset Y}\int_{F}
  \sum_{\substack{\gamma\in\Gamma'\\ \text{parabolic}}}\tilde\chi_A(z) k(\sig(z,\gamma
  z))\vol(z) \biggr)=0.
\end{equation}
The primitive elements $\gamma_0\in B'$ correspond bijectively to the
simple closed geodesics of $Y'$, so the sum of \eqref{sef1} and \eqref{sef2}
yields the claim.

\bibliographystyle{plain}

\end{document}